\pdfoutput=1

\documentclass[11pt,a4paper, twoside]{amsart}
\usepackage{geometry}
\usepackage[utf8]{inputenc}
\usepackage[english]{babel}
\usepackage[T1]{fontenc}
\usepackage{lmodern}

\usepackage{amsmath}
\usepackage{amssymb}
\usepackage{amsthm}

\usepackage{mathtools}
\usepackage{stmaryrd}

\usepackage{wrapfig}
\usepackage{subcaption}

\usepackage{enumitem}
\usepackage{verbatim}
\usepackage{wrapfig}
\usepackage{makecell}
\usepackage{xcolor}

\usepackage{hyperref}
\hypersetup{colorlinks,linkcolor={blue!75!black},citecolor={blue!75!black},urlcolor={blue!75!black}}

\usepackage[backend=biber]{biblatex}
\usepackage{tikz}

\DeclareMathOperator*{\id}{id}

\DeclareMathOperator*{\supp}{supp}

\DeclareMathOperator*{\linh}{span}
\DeclareMathOperator*{\sgn}{sgn}

\DeclareMathOperator*{\conv}{conv}
\DeclareMathOperator*{\aff}{aff}
\DeclareMathOperator*{\pos}{pos}
\DeclareMathOperator*{\bd}{bd}
\DeclareMathOperator*{\cl}{cl}
\DeclareMathOperator*{\relbd}{relbd}

\DeclareMathOperator*{\diam}{diam}
\DeclareMathOperator*{\TV}{TV}
\DeclareMathOperator*{\Wt}{W}

\DeclareMathOperator*{\dist}{dist}

\newcommand{\MID}{\,\middle|\,}
\DeclarePairedDelimiter\scpr{\langle}{\rangle}

\newcommand{\norm}[1]{\|#1\|}
\newcommand{\NORM}[1]{\left\|#1\right\|}
\newcommand{\SCPR}[1]{\left\langle #1\right\rangle}

\newcommand{\R}{\mathbb{R}}

\newcommand{\N}{\mathbb{N}}

\mathtoolsset{showonlyrefs}

\theoremstyle{theorem}

\newtheorem{theorem}{Theorem}[section] %
\newtheorem{corollary}[theorem]{Corollary}
\newtheorem{proposition}[theorem]{Proposition}
\newtheorem{lemma}[theorem]{Lemma}
\newtheorem{claim}[theorem]{Claim}

\theoremstyle{definition}
\newtheorem{question}[theorem]{Question}
\newtheorem{definition}[theorem]{Definition}
\newtheorem{construction}[theorem]{Construction}

\theoremstyle{remark}
\newtheorem{remark}[theorem]{Remark}
\newtheorem{example}[theorem]{Example}

\numberwithin{equation}{section}

\DeclareMathOperator*{\vertices}{vert}
\DeclareMathOperator*{\interior}{int}
\DeclareMathOperator*{\relint}{relint}
\DeclareMathOperator*{\vol}{vol}
\DeclareMathOperator*{\wlim}{wlim}
\newcommand{\volt}[1]{{\textstyle\vol_{#1}}}

\DeclareMathOperator*{\cvx}{cvx}

\DeclareRobustCommand{\One}{\text{\usefont{U}{bbold}{m}{n}1}}

\newcommand{\transp}{\mathsf{T}}
\newcommand{\Kn}{\mathcal{K}^n}
\newcommand{\Knn}{\mathcal{K}_n^n}
\newcommand{\Pn}{\mathcal{P}^n}
\newcommand{\Pnn}{\mathcal{P}_n^n}

\newcommand{\E}{\mathbb{E}} 
\newcommand{\M}{\mathsf{M}}
\newcommand{\W}{\mathsf{W}}
\newcommand{\sphere}{\mathbb{S}^{n-1}}

\newcommand{\dint}{\,\mathrm{d}}
\newcommand{\ddt}{\frac{\mathrm{d}}{\mathrm{d}t}}

\newcommand{\eqand}{\quad \text{and} \quad}
\newcommand{\dH}{d_{\text{H}}}
\newcommand{\dS}{d_{\text{S}}}
\newcommand{\goestoinfty}[1]{\underset{#1 \rightarrow \infty}{\longrightarrow}}
\newcommand{\goestosth}[2]{\underset{#1 \rightarrow #2}{\longrightarrow}}

\newcommand{\warrow}{\overset{w}{\rightarrow}}

\newcommand{\cdotbox}{\hspace{0.4mm}\cdot\hspace{0.4mm}}

\let\emptyset\varnothing

\title{Small Perturbations of Polytopes}
\author{Christian Kipp}

\address{Technische Universität Berlin, Institut für Mathematik, Sekr.~MA4-1, Straße des 17.~Juni 136, 10623 Berlin, Germany}
\email{kipp@math.tu-berlin.de}
\begin{document}
	
	\begin{abstract}
		Motivated by first-order conditions for extremal bodies of geometric functionals, we study a functional analytic notion of infinitesimal perturbations of convex bodies and give a full characterization of the set of realizable perturbations if the perturbed body is a polytope. As an application, we derive a necessary condition for polytopal maximizers of the isotropic constant.
	\end{abstract}
	
	\maketitle
	
	\setcounter{tocdepth}{1}
	
	\tableofcontents
	
	\section{Introduction} \label{sect_introduction}
	
	In this article, a convex body is a non-empty compact convex subset of $\R^n$. We equip $\R^n$ with the standard scalar product $\scpr{\cdot,\cdot}$. The space of convex bodies in $\R^n$ is denoted by $\Kn$ and the subspace of full-dimensional convex bodies by $\Knn$.
	
	In convex geometry, one often encounters functionals $\phi \colon \mathcal{A} \rightarrow \R$ of the form
	\begin{equation} \label{eq_functional_phi}
		\phi(K) = g\left(\int_K f_1(x) \dint x, \dots, \int_K f_m(x) \dint x\right),
	\end{equation}
	where $\mathcal{A} \subset \Kn$, $f_1, \dots, f_m$ are continuous functions $\R^n \rightarrow \R$ and $g \colon \R^m \rightarrow \R$ is differentiable. Examples of such functionals include the volume, components of the centroid, the moment of inertia, %
	the difference between the left-hand and right-hand sides in the $B$-theorem, and the isotropic constant. If one is interested in extremal questions related to such functionals, it is natural to ask how $\phi$ changes if $K$ is slightly ``perturbed'', i.e., replaced by another convex body that is close to $K$. Instead of considering a single perturbation, we consider a whole family of convex bodies, absolving us temporarily of the question what we mean exactly by ``close''.  Let $(K_t)_{t \in [0,1]}$ be a one-parameter family of convex bodies in $\mathcal{A}$ for which the derivatives $\left.\ddt \int_{K_t} f_i(x) \dint x\right|_{t=0}$ exist. Then the derivative of $\phi$ ``along'' $(K_t)_{t \in [0,1]}$ is given by 
	\begin{equation} \label{eq_derivative_of_phi}
		\left.\ddt \phi(K_t) \right|_{t=0}=  \sum_{i=1}^m \frac{\partial g}{\partial x_i} \left(\int_{K_0} f_1(x) \dint x, \dots, \int_{K_0} f_m(x) \dint x\right) \cdot \left.\ddt \int_{K_t} f_i(x) \dint x\right|_{t=0}.
	\end{equation}
	As an example, we consider the isotropic constant $L_K$ of a full-dimensional convex body $K\in \Knn$, which is given by
	\begin{equation} \label{eq_def_l_k}
		L_K^{2n} = \frac{\det \left[\int_{K-c_K}x_i x_j \dint x\right]}{\vol(K)^{n+2}}, \quad \text{where } c_K=\frac{1}{\vol K}\int_K x \dint x.
	\end{equation}
	Because the centroid $c_K$ modifies the domain of integration, the functional $L_K^{2n }$ is not of the form \eqref{eq_functional_phi}, but this turns out to be a minor technicality. Let $\mathcal{A}\subset \Knn$ be the class of centered convex bodies, i.e., those $K\in \Knn$ with $c_K=0$. Defining $\phi \colon \Knn \rightarrow \R$ by
	\begin{equation}
		\phi(K) = \frac{\det \left[\int_{K}x_i x_j \dint x\right]}{\vol(K)^{n+2}},
	\end{equation}
	we clearly have $\phi(K)=L_K^{2n}$ for $K \in \mathcal{A}$. Let $(K_t)_{t \in [0,1]}$ be a one-parameter family in $\Knn$ (not necessarily in $\mathcal{A}$). It was shown by Rademacher \cite{rad16} that if $K_0$ is \textit{isotropic}, i.e., it satisfies
	\begin{equation} \label{eq_def_isotropic}
		c_{K_0} = 0 \quad \text{and} \quad \frac{1}{\vol K_0} \int_{K_0}x_i x_j \dint x = \delta_{ij},
	\end{equation}
	then the derivatives of $t\mapsto \phi(K_t)$ and $t\mapsto L_{K_t}^{2n}$ at $t=0$ coincide and are given by
	\begin{equation} \label{eq_rademacher}
		\left.\ddt \phi(K_t) \right|_{t=0}=\left.\ddt L_{K_t}^{2n} \right|_{t=0}=\frac{1}{\vol(K_0)^3}\left.\left(\ddt\int_{K_t}[\norm{x}_2^2-n-2] \dint x\right)\right|_{t=0},
	\end{equation}
	if the derivative on the right-hand side exists. Using certain one-parameter families $(K_t)_{t \in [-1,1]}$, which can intuitively be described as ``hinging the facets'', Rademacher used the first-order conditions that follow from \eqref{eq_rademacher} to derive the remarkable conclusion that a simplicial polytope $P$ that maximizes $K \mapsto L_K$ must be a simplex \cite{rad16}.
	
	\begin{figure}[ht]
		\begin{subfigure}[b]{.33\linewidth}
			\centering
			\begin{tikzpicture}%
	[x={(0.6cm, -0.4cm)},
	y={(0.827112cm, 0.224669cm)},
	z={(0.000024cm, 0.916614cm)},
	scale=1.700000,
	back/.style={dotted, thin},
	edge/.style={color=black},
	facet/.style={fill=white,fill opacity=0.200000},
	vertex/.style={},
	back2/.style={dotted},
	edge2/.style={color=black},
	facet2/.style={fill=red,fill opacity=0.2000000},
	vertex2/.style={},
	back3/.style={dotted, thin},
	edge3/.style={color=black,very thick},
	facet3/.style={fill=cyan,fill opacity=0.200000},
	vertex3/.style={}]
	
	\coordinate (0.00000, 0.50000, 0.80902) at (0.00000, 0.50000, 0.80902);
	\coordinate (0.00000, -0.50000, 0.80902) at (0.00000, -0.50000, 0.80902);
	\coordinate (0.50000, 0.80902, 0.00000) at (0.50000, 0.80902, 0.00000);
	\coordinate (0.50000, -0.80902, 0.00000) at (0.50000, -0.80902, 0.00000);
	\coordinate (0.80902, 0.00000, 0.50000) at (0.80902, 0.00000, 0.50000);
	\coordinate (0.80902, 0.00000, -0.50000) at (0.80902, 0.00000, -0.50000);
	\coordinate (-0.50000, 0.80902, 0.00000) at (-0.50000, 0.80902, 0.00000);
	\coordinate (-0.50000, -0.80902, 0.00000) at (-0.50000, -0.80902, 0.00000);
	\coordinate (-0.80902, 0.00000, 0.50000) at (-0.80902, 0.00000, 0.50000);
	\coordinate (0.00000, 0.50000, -0.80902) at (0.00000, 0.50000, -0.80902);
	\coordinate (0.00000, -0.50000, -0.80902) at (0.00000, -0.50000, -0.80902);
	\coordinate (-0.80902, 0.00000, -0.50000) at (-0.80902, 0.00000, -0.50000);
	\draw[edge,back] (0.00000, 0.50000, 0.80902) -- (-0.50000, 0.80902, 0.00000);
	\draw[edge,back] (0.50000, 0.80902, 0.00000) -- (-0.50000, 0.80902, 0.00000);
	\draw[edge,back] (0.50000, 0.80902, 0.00000) -- (0.00000, 0.50000, -0.80902);
	\draw[edge,back] (0.80902, 0.00000, -0.50000) -- (0.00000, 0.50000, -0.80902);
	\draw[edge,back] (-0.50000, 0.80902, 0.00000) -- (-0.80902, 0.00000, 0.50000);
	\draw[edge,back] (-0.50000, 0.80902, 0.00000) -- (0.00000, 0.50000, -0.80902);
	\draw[edge,back] (-0.50000, 0.80902, 0.00000) -- (-0.80902, 0.00000, -0.50000);
	\draw[edge,back] (-0.50000, -0.80902, 0.00000) -- (-0.80902, 0.00000, -0.50000);
	\draw[edge,back] (-0.80902, 0.00000, 0.50000) -- (-0.80902, 0.00000, -0.50000);
	\draw[edge,back] (0.00000, 0.50000, -0.80902) -- (0.00000, -0.50000, -0.80902);
	\draw[edge,back] (0.00000, 0.50000, -0.80902) -- (-0.80902, 0.00000, -0.50000);
	\draw[edge,back] (0.00000, -0.50000, -0.80902) -- (-0.80902, 0.00000, -0.50000);
	\fill[facet] (-0.50000, -0.80902, 0.00000) -- (0.00000, -0.50000, 0.80902) -- (0.50000, -0.80902, 0.00000) -- cycle {};
	\fill[facet] (0.00000, -0.50000, -0.80902) -- (0.50000, -0.80902, 0.00000) -- (-0.50000, -0.80902, 0.00000) -- cycle {};
	\fill[facet] (0.80902, 0.00000, 0.50000) -- (0.00000, 0.50000, 0.80902) -- (0.50000, 0.80902, 0.00000) -- cycle {};
	\fill[facet,gray] (0.80902, 0.00000, 0.50000) -- (0.00000, 0.50000, 0.80902) -- (0.00000, -0.50000, 0.80902) -- cycle {};
	\fill[facet] (0.80902, 0.00000, 0.50000) -- (0.00000, -0.50000, 0.80902) -- (0.50000, -0.80902, 0.00000) -- cycle {};
	\fill[facet] (-0.80902, 0.00000, 0.50000) -- (0.00000, -0.50000, 0.80902) -- (-0.50000, -0.80902, 0.00000) -- cycle {};
	\fill[facet] (-0.80902, 0.00000, 0.50000) -- (0.00000, 0.50000, 0.80902) -- (0.00000, -0.50000, 0.80902) -- cycle {};
	\fill[facet] (0.00000, -0.50000, -0.80902) -- (0.50000, -0.80902, 0.00000) -- (0.80902, 0.00000, -0.50000) -- cycle {};
	\fill[facet] (0.80902, 0.00000, -0.50000) -- (0.50000, 0.80902, 0.00000) -- (0.80902, 0.00000, 0.50000) -- cycle {};
	\fill[facet] (0.80902, 0.00000, -0.50000) -- (0.50000, -0.80902, 0.00000) -- (0.80902, 0.00000, 0.50000) -- cycle {};
	\draw[edge3] (0.00000, 0.50000, 0.80902) -- (0.00000, -0.50000, 0.80902);
	\draw[edge] (0.00000, 0.50000, 0.80902) -- (0.50000, 0.80902, 0.00000);
	\draw[edge] (0.00000, 0.50000, 0.80902) -- (0.80902, 0.00000, 0.50000);
	\draw[edge] (0.00000, 0.50000, 0.80902) -- (-0.80902, 0.00000, 0.50000);
	\draw[edge] (0.00000, -0.50000, 0.80902) -- (0.50000, -0.80902, 0.00000);
	\draw[edge] (0.00000, -0.50000, 0.80902) -- (0.80902, 0.00000, 0.50000);
	\draw[edge] (0.00000, -0.50000, 0.80902) -- (-0.50000, -0.80902, 0.00000);
	\draw[edge] (0.00000, -0.50000, 0.80902) -- (-0.80902, 0.00000, 0.50000);
	\draw[edge] (0.50000, 0.80902, 0.00000) -- (0.80902, 0.00000, 0.50000);
	\draw[edge] (0.50000, 0.80902, 0.00000) -- (0.80902, 0.00000, -0.50000);
	\draw[edge] (0.50000, -0.80902, 0.00000) -- (0.80902, 0.00000, 0.50000);
	\draw[edge] (0.50000, -0.80902, 0.00000) -- (0.80902, 0.00000, -0.50000);
	\draw[edge] (0.50000, -0.80902, 0.00000) -- (-0.50000, -0.80902, 0.00000);
	\draw[edge] (0.50000, -0.80902, 0.00000) -- (0.00000, -0.50000, -0.80902);
	\draw[edge] (0.80902, 0.00000, 0.50000) -- (0.80902, 0.00000, -0.50000);
	\draw[edge] (0.80902, 0.00000, -0.50000) -- (0.00000, -0.50000, -0.80902);
	\draw[edge] (-0.50000, -0.80902, 0.00000) -- (-0.80902, 0.00000, 0.50000);
	\draw[edge] (-0.50000, -0.80902, 0.00000) -- (0.00000, -0.50000, -0.80902);
\end{tikzpicture}
			
			\caption{simplicial polytope $P$}
		\end{subfigure}
		\begin{subfigure}[b]{.32\linewidth}
			\centering
			\begin{tikzpicture}%
	[x={(0.6cm, -0.4cm)},
	y={(0.827112cm, 0.224669cm)},
	z={(0.000024cm, 0.916614cm)},
	scale=1.700000,
	back/.style={dotted, thin},
	edge/.style={color=black},
	facet/.style={fill=white,fill opacity=0.000000},
	vertex/.style={},
	back2/.style={dotted},
	edge2/.style={color=black},
	facet2/.style={fill=orange,fill opacity=0.20000},
	vertex2/.style={},
	back3/.style={dotted, thin},
	edge3/.style={color=black,very thick},
	facet3/.style={fill=green,fill opacity=0.20000},
	vertex3/.style={}]
	
	\coordinate (0.00000, 0.50000, 0.80902) at (0.00000, 0.50000, 0.80902);
	\coordinate (0.00000, -0.50000, 0.80902) at (0.00000, -0.50000, 0.80902);
	\coordinate (0.50000, 0.80902, 0.00000) at (0.50000, 0.80902, 0.00000);
	\coordinate (0.50000, -0.80902, 0.00000) at (0.50000, -0.80902, 0.00000);
	\coordinate (0.80902, 0.00000, 0.50000) at (0.80902, 0.00000, 0.50000);
	\coordinate (0.80902, 0.00000, -0.50000) at (0.80902, 0.00000, -0.50000);
	\coordinate (-0.50000, 0.80902, 0.00000) at (-0.50000, 0.80902, 0.00000);
	\coordinate (-0.50000, -0.80902, 0.00000) at (-0.50000, -0.80902, 0.00000);
	\coordinate (-0.80902, 0.00000, 0.50000) at (-0.80902, 0.00000, 0.50000);
	\coordinate (0.00000, 0.50000, -0.80902) at (0.00000, 0.50000, -0.80902);
	\coordinate (0.00000, -0.50000, -0.80902) at (0.00000, -0.50000, -0.80902);
	\coordinate (-0.80902, 0.00000, -0.50000) at (-0.80902, 0.00000, -0.50000);
	\draw[edge,back] (0.00000, 0.50000, 0.80902) -- (-0.50000, 0.80902, 0.00000);
	\draw[edge,back] (0.50000, 0.80902, 0.00000) -- (-0.50000, 0.80902, 0.00000);
	\draw[edge,back] (0.50000, 0.80902, 0.00000) -- (0.00000, 0.50000, -0.80902);
	\draw[edge,back] (0.80902, 0.00000, -0.50000) -- (0.00000, 0.50000, -0.80902);
	\draw[edge,back] (-0.50000, 0.80902, 0.00000) -- (-0.80902, 0.00000, 0.50000);
	\draw[edge,back] (-0.50000, 0.80902, 0.00000) -- (0.00000, 0.50000, -0.80902);
	\draw[edge,back] (-0.50000, 0.80902, 0.00000) -- (-0.80902, 0.00000, -0.50000);
	\draw[edge,back] (-0.50000, -0.80902, 0.00000) -- (-0.80902, 0.00000, -0.50000);
	\draw[edge,back] (-0.80902, 0.00000, 0.50000) -- (-0.80902, 0.00000, -0.50000);
	\draw[edge,back] (0.00000, 0.50000, -0.80902) -- (0.00000, -0.50000, -0.80902);
	\draw[edge,back] (0.00000, 0.50000, -0.80902) -- (-0.80902, 0.00000, -0.50000);
	\draw[edge,back] (0.00000, -0.50000, -0.80902) -- (-0.80902, 0.00000, -0.50000);
	\fill[facet] (-0.50000, -0.80902, 0.00000) -- (0.00000, -0.50000, 0.80902) -- (0.50000, -0.80902, 0.00000) -- cycle {};
	\fill[facet] (0.00000, -0.50000, -0.80902) -- (0.50000, -0.80902, 0.00000) -- (-0.50000, -0.80902, 0.00000) -- cycle {};
	\fill[facet] (0.80902, 0.00000, 0.50000) -- (0.00000, 0.50000, 0.80902) -- (0.50000, 0.80902, 0.00000) -- cycle {};
	\fill[facet] (0.80902, 0.00000, 0.50000) -- (0.00000, 0.50000, 0.80902) -- (0.00000, -0.50000, 0.80902) -- cycle {};
	\fill[facet] (0.80902, 0.00000, 0.50000) -- (0.00000, -0.50000, 0.80902) -- (0.50000, -0.80902, 0.00000) -- cycle {};
	\fill[facet] (-0.80902, 0.00000, 0.50000) -- (0.00000, -0.50000, 0.80902) -- (-0.50000, -0.80902, 0.00000) -- cycle {};
	\fill[facet] (-0.80902, 0.00000, 0.50000) -- (0.00000, 0.50000, 0.80902) -- (0.00000, -0.50000, 0.80902) -- cycle {};
	\fill[facet] (0.00000, -0.50000, -0.80902) -- (0.50000, -0.80902, 0.00000) -- (0.80902, 0.00000, -0.50000) -- cycle {};
	\fill[facet] (0.80902, 0.00000, -0.50000) -- (0.50000, 0.80902, 0.00000) -- (0.80902, 0.00000, 0.50000) -- cycle {};
	\fill[facet] (0.80902, 0.00000, -0.50000) -- (0.50000, -0.80902, 0.00000) -- (0.80902, 0.00000, 0.50000) -- cycle {};
	\draw[edge] (0.00000, 0.50000, 0.80902) -- (0.00000, -0.50000, 0.80902);
	\draw[edge] (0.00000, 0.50000, 0.80902) -- (0.50000, 0.80902, 0.00000);
	\draw[edge] (0.00000, 0.50000, 0.80902) -- (0.80902, 0.00000, 0.50000);
	\draw[edge] (0.00000, 0.50000, 0.80902) -- (-0.80902, 0.00000, 0.50000);
	\draw[edge] (0.00000, -0.50000, 0.80902) -- (0.50000, -0.80902, 0.00000);
	\draw[edge] (0.00000, -0.50000, 0.80902) -- (0.80902, 0.00000, 0.50000);
	\draw[edge] (0.00000, -0.50000, 0.80902) -- (-0.50000, -0.80902, 0.00000);
	\draw[edge] (0.00000, -0.50000, 0.80902) -- (-0.80902, 0.00000, 0.50000);
	\draw[edge] (0.50000, 0.80902, 0.00000) -- (0.80902, 0.00000, 0.50000);
	\draw[edge] (0.50000, 0.80902, 0.00000) -- (0.80902, 0.00000, -0.50000);
	\draw[edge] (0.50000, -0.80902, 0.00000) -- (0.80902, 0.00000, 0.50000);
	\draw[edge] (0.50000, -0.80902, 0.00000) -- (0.80902, 0.00000, -0.50000);
	\draw[edge] (0.50000, -0.80902, 0.00000) -- (-0.50000, -0.80902, 0.00000);
	\draw[edge] (0.50000, -0.80902, 0.00000) -- (0.00000, -0.50000, -0.80902);
	\draw[edge] (0.80902, 0.00000, 0.50000) -- (0.80902, 0.00000, -0.50000);
	\draw[edge] (0.80902, 0.00000, -0.50000) -- (0.00000, -0.50000, -0.80902);
	\draw[edge] (-0.50000, -0.80902, 0.00000) -- (-0.80902, 0.00000, 0.50000);
	\draw[edge] (-0.50000, -0.80902, 0.00000) -- (0.00000, -0.50000, -0.80902);
	
	\coordinate (0.77049, 0.10086, 0.43767) at (0.77049, 0.10086, 0.43767);
	\coordinate (0.80902, 0.00000, 0.41910) at (0.80902, 0.00000, 0.41910);
	\coordinate (0.77049, -0.10086, 0.43767) at (0.77049, -0.10086, 0.43767);
	\coordinate (0.00000, 0.50000, 0.80902) at (0.00000, 0.50000, 0.80902);
	\coordinate (0.00000, -0.50000, 0.80902) at (0.00000, -0.50000, 0.80902);
	\coordinate (0.50000, 0.80902, 0.00000) at (0.50000, 0.80902, 0.00000);
	\coordinate (0.50000, -0.80902, 0.00000) at (0.50000, -0.80902, 0.00000);
	\coordinate (0.80902, 0.00000, -0.50000) at (0.80902, 0.00000, -0.50000);
	\coordinate (-0.80902, 0.00000, 0.50000) at (-0.80902, 0.00000, 0.50000);
	\coordinate (-0.50000, 0.80902, 0.00000) at (-0.50000, 0.80902, 0.00000);
	\coordinate (-0.50000, -0.80902, 0.00000) at (-0.50000, -0.80902, 0.00000);
	\coordinate (0.00000, 0.50000, -0.80902) at (0.00000, 0.50000, -0.80902);
	\coordinate (0.00000, -0.50000, -0.80902) at (0.00000, -0.50000, -0.80902);
	\coordinate (-0.80902, 0.00000, -0.50000) at (-0.80902, 0.00000, -0.50000);
	\fill[facet] (0.50000, 0.80902, 0.00000) -- (0.77049, 0.10086, 0.43767) -- (0.00000, 0.50000, 0.80902) -- cycle {};
	\fill[facet2] (0.00000, -0.50000, 0.80902) -- (0.77049, -0.10086, 0.43767) -- (0.80902, 0.00000, 0.41910) -- (0.77049, 0.10086, 0.43767) -- (0.00000, 0.50000, 0.80902) -- cycle {};
	\fill[facet] (0.50000, -0.80902, 0.00000) -- (0.77049, -0.10086, 0.43767) -- (0.00000, -0.50000, 0.80902) -- cycle {};00000
	\fill[facet] (0.80902, 0.00000, -0.50000) -- (0.80902, 0.00000, 0.41910) -- (0.77049, 0.10086, 0.43767) -- (0.50000, 0.80902, 0.00000) -- cycle {};
	\fill[facet] (0.80902, 0.00000, -0.50000) -- (0.80902, 0.00000, 0.41910) -- (0.77049, -0.10086, 0.43767) -- (0.50000, -0.80902, 0.00000) -- cycle {};
	\fill[facet] (-0.80902, 0.00000, 0.50000) -- (0.00000, 0.50000, 0.80902) -- (0.00000, -0.50000, 0.80902) -- cycle {};
	\fill[facet] (-0.50000, -0.80902, 0.00000) -- (0.00000, -0.50000, 0.80902) -- (0.50000, -0.80902, 0.00000) -- cycle {};
	\fill[facet] (-0.50000, -0.80902, 0.00000) -- (0.00000, -0.50000, 0.80902) -- (-0.80902, 0.00000, 0.50000) -- cycle {};
	\fill[facet] (0.00000, -0.50000, -0.80902) -- (0.50000, -0.80902, 0.00000) -- (0.80902, 0.00000, -0.50000) -- cycle {};
	\fill[facet] (0.00000, -0.50000, -0.80902) -- (0.50000, -0.80902, 0.00000) -- (-0.50000, -0.80902, 0.00000) -- cycle {};
	\draw[edge,gray] (0.77049, 0.10086, 0.43767) -- (0.80902, 0.00000, 0.41910);
	\draw[edge,gray] (0.77049, 0.10086, 0.43767) -- (0.00000, 0.50000, 0.80902);
	\draw[edge,gray] (0.80902, 0.00000, 0.41910) -- (0.77049, -0.10086, 0.43767);
	\draw[edge,gray] (0.77049, -0.10086, 0.43767) -- (0.00000, -0.50000, 0.80902);
	\draw[edge3] (0.00000, 0.50000, 0.80902) -- (0.00000, -0.50000, 0.80902);
\end{tikzpicture}
			
			\caption{facet hinges inward}
		\end{subfigure}
		\begin{subfigure}[b]{.32\linewidth}
			\centering
			\begin{tikzpicture}%
	[x={(0.6cm, -0.4cm)},
	y={(0.827112cm, 0.224669cm)},
	z={(0.000024cm, 0.916614cm)},
	scale=1.700000,
	back/.style={dotted, thin},
	edge/.style={color=black},
	facet/.style={fill=white,fill opacity=0.000000},
	vertex/.style={},
	back2/.style={dotted},
	edge2/.style={color=black},
	facet2/.style={fill=red,fill opacity=0.2500000},
	vertex2/.style={},
	back3/.style={dotted, thin},
	edge3/.style={color=black,very thick},
	facet3/.style={fill=cyan,fill opacity=0.20000},
	vertex3/.style={}]
	
	\draw[edge,back] (0.00000, 0.50000, 0.80902) -- (-0.50000, 0.80902, 0.00000);
	\draw[edge,back] (0.50000, 0.80902, 0.00000) -- (-0.50000, 0.80902, 0.00000);
	\draw[edge,back] (0.50000, 0.80902, 0.00000) -- (0.00000, 0.50000, -0.80902);
	\draw[edge,back] (0.80902, 0.00000, -0.50000) -- (0.00000, 0.50000, -0.80902);
	\draw[edge,back] (-0.50000, 0.80902, 0.00000) -- (-0.80902, 0.00000, 0.50000);
	\draw[edge,back] (-0.50000, 0.80902, 0.00000) -- (0.00000, 0.50000, -0.80902);
	\draw[edge,back] (-0.50000, 0.80902, 0.00000) -- (-0.80902, 0.00000, -0.50000);
	\draw[edge,back] (-0.50000, -0.80902, 0.00000) -- (-0.80902, 0.00000, -0.50000);
	\draw[edge,back] (-0.80902, 0.00000, 0.50000) -- (-0.80902, 0.00000, -0.50000);
	\draw[edge,back] (0.00000, 0.50000, -0.80902) -- (0.00000, -0.50000, -0.80902);
	\draw[edge,back] (0.00000, 0.50000, -0.80902) -- (-0.80902, 0.00000, -0.50000);
	\draw[edge,back] (0.00000, -0.50000, -0.80902) -- (-0.80902, 0.00000, -0.50000);
	\fill[facet] (-0.50000, -0.80902, 0.00000) -- (0.00000, -0.50000, 0.80902) -- (0.50000, -0.80902, 0.00000) -- cycle {};
	\fill[facet] (0.00000, -0.50000, -0.80902) -- (0.50000, -0.80902, 0.00000) -- (-0.50000, -0.80902, 0.00000) -- cycle {};
	\fill[facet] (0.80902, 0.00000, 0.50000) -- (0.00000, 0.50000, 0.80902) -- (0.50000, 0.80902, 0.00000) -- cycle {};
	\fill[facet] (0.80902, 0.00000, 0.50000) -- (0.00000, 0.50000, 0.80902) -- (0.00000, -0.50000, 0.80902) -- cycle {};
	\fill[facet] (0.80902, 0.00000, 0.50000) -- (0.00000, -0.50000, 0.80902) -- (0.50000, -0.80902, 0.00000) -- cycle {};
	\fill[facet] (-0.80902, 0.00000, 0.50000) -- (0.00000, -0.50000, 0.80902) -- (-0.50000, -0.80902, 0.00000) -- cycle {};
	\fill[facet] (-0.80902, 0.00000, 0.50000) -- (0.00000, 0.50000, 0.80902) -- (0.00000, -0.50000, 0.80902) -- cycle {};
	\fill[facet] (0.00000, -0.50000, -0.80902) -- (0.50000, -0.80902, 0.00000) -- (0.80902, 0.00000, -0.50000) -- cycle {};
	\fill[facet] (0.80902, 0.00000, -0.50000) -- (0.50000, 0.80902, 0.00000) -- (0.80902, 0.00000, 0.50000) -- cycle {};
	\fill[facet] (0.80902, 0.00000, -0.50000) -- (0.50000, -0.80902, 0.00000) -- (0.80902, 0.00000, 0.50000) -- cycle {};
	\draw[edge] (0.00000, 0.50000, 0.80902) -- (0.00000, -0.50000, 0.80902);
	\draw[edge] (0.00000, 0.50000, 0.80902) -- (0.50000, 0.80902, 0.00000);
	\draw[edge,gray] (0.00000, 0.50000, 0.80902) -- (0.80902, 0.00000, 0.50000);
	\draw[edge] (0.00000, 0.50000, 0.80902) -- (-0.80902, 0.00000, 0.50000);
	\draw[edge] (0.00000, -0.50000, 0.80902) -- (0.50000, -0.80902, 0.00000);
	\draw[edge,gray] (0.00000, -0.50000, 0.80902) -- (0.80902, 0.00000, 0.50000);
	\draw[edge] (0.00000, -0.50000, 0.80902) -- (-0.50000, -0.80902, 0.00000);
	\draw[edge] (0.00000, -0.50000, 0.80902) -- (-0.80902, 0.00000, 0.50000);
	\draw[edge] (0.50000, 0.80902, 0.00000) -- (0.80902, 0.00000, 0.50000);
	\draw[edge] (0.50000, 0.80902, 0.00000) -- (0.80902, 0.00000, -0.50000);
	\draw[edge] (0.50000, -0.80902, 0.00000) -- (0.80902, 0.00000, 0.50000);
	\draw[edge] (0.50000, -0.80902, 0.00000) -- (0.80902, 0.00000, -0.50000);
	\draw[edge] (0.50000, -0.80902, 0.00000) -- (-0.50000, -0.80902, 0.00000);
	\draw[edge] (0.50000, -0.80902, 0.00000) -- (0.00000, -0.50000, -0.80902);
	\draw[edge] (0.80902, 0.00000, 0.50000) -- (0.80902, 0.00000, -0.50000);
	\draw[edge] (0.80902, 0.00000, -0.50000) -- (0.00000, -0.50000, -0.80902);
	\draw[edge] (-0.50000, -0.80902, 0.00000) -- (-0.80902, 0.00000, 0.50000);
	\draw[edge] (-0.50000, -0.80902, 0.00000) -- (0.00000, -0.50000, -0.80902);
	
	\fill[facet3] (0.00000, -0.50000, 0.80902) -- (0.69635, 0.00000, 0.61267) -- (0.00000, 0.50000, 0.80902) -- cycle {};
	\fill[facet3] (0.00000, -0.50000, 0.80902) -- (0.69635, 0.00000, 0.61267) -- (0.80902, 0.00000, 0.50000) -- cycle {};
	\fill[facet3] (0.00000, 0.50000, 0.80902) -- (0.69635, 0.00000, 0.61267) -- (0.80902, 0.00000, 0.50000) -- cycle {};
	\draw[edge] (0.69635, 0.00000, 0.61267) -- (0.80902, 0.00000, 0.50000);
	\draw[edge] (0.69635, 0.00000, 0.61267) -- (0.00000, 0.50000, 0.80902);
	\draw[edge] (0.69635, 0.00000, 0.61267) -- (0.00000, -0.50000, 0.80902);
	\draw[edge3] (0.00000, 0.50000, 0.80902) -- (0.00000, -0.50000, 0.80902);
\end{tikzpicture}
			
			\caption{facet hinges outward}
		\end{subfigure}
		\caption{Rademacher's one-parameter family from \cite{rad16}. The gray facet hinges around the thickly drawn edge.}
	\end{figure}
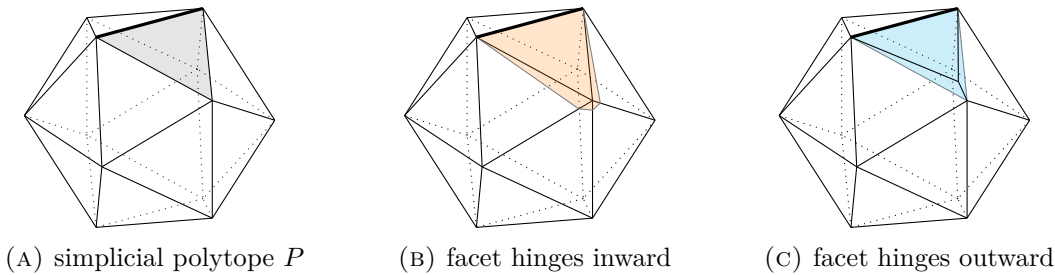
	
	The present article is motivated by the question whether similar perturbation arguments might be used to derive further information about the extremal bodies of the isotropic constant (similar arguments are used in \cite{rad12} and \cite{mr16}). Phrased somewhat more generally, one might ask which derivatives ``can appear'' on the right-hand side of \eqref{eq_rademacher}. Clearly, in order to make this question meaningful, some geometric restrictions on the family $(K_t)_{t \in [0,1]}$ are needed to ensure that it actually describes a ``local'' perturbation. To make this notion precise, we return to the more general setting of the derivative given in \eqref{eq_derivative_of_phi}. The exact geometric restriction that we impose on the family $(K_t)_{t \in [0,1]}$ is that \textit{all derivatives} of the form \eqref{eq_derivative_of_phi} exist, i.e., for all functionals $\phi$ of the form \eqref{eq_functional_phi}. This leads us to the following definition, where we denote by $C(\R^n)$ the vector space of continuous functions $\R^n \rightarrow \R$.
	\begin{definition} \label{def_weak_derivative}
		Let $(K_t)_{t \in [0,1]}$ be a family of full-dimensional convex bodies. For every $f\in C(\R^n)$, we obtain a function $\phi_f(t)\colon [0,1] \rightarrow \R$ via
		\begin{equation}
			\phi_f(t) \coloneqq \int_{K_t} f(x) \dint x.
		\end{equation} 
		We say that $(K_t)_{t \in [0,1]}$ is \textit{weakly differentiable} if $\phi_f(t)$ is differentiable (from the right) at $t=0$ for every $f\in C(\R^n)$. The linear map
		\begin{equation}
			K'_0 \colon C(\R^n) \rightarrow \R, \quad f \mapsto \left.\ddt \phi_f(t) \right|_{t=0}
		\end{equation}
		is called the \textit{weak derivative} of $(K_t)_{t \in [0,1]}$.
	\end{definition}
	Using an analogy from differential geometry, the weak derivative $K_0'$ can be thought of as a ``tangent vector'' at $K_0$ to the space $\Knn$, which describes the directional derivative along the curve $(K_t)_{t \in [0,1]}$ at $t=0$. Based on Definition \ref{def_weak_derivative}, we are now able to state our problem more precisely.
	\begin{question} \label{qu_first_question}
		Which linear maps $C(\R^n)\rightarrow \R$ can be realized (i.e., can appear) as weak derivatives of families $(K_t)_{t \in [0,1]}$?
	\end{question}

	A finite signed Borel measure on %
	$\R^n$ is a set function of the form $\alpha_1 \mu_1 +\alpha_2 \mu_2$, where $\mu_1, \mu_2$ are probability measures on the Borel $\sigma$-algebra and $\alpha_1,\alpha_2 \in \R$ are arbitrary scalars. In the following, we will refer to such set functions simply as \textit{signed measures}. Defining the support of a (signed) measure is somewhat complicated; for our purposes, it suffices to define that ``$\mu$ is supported on $A$'' means $\mu(B)=0$ for all measurable sets $B \subset X \setminus A$. 
	Let $X\subset \R^n$ be a Borel set. We denote the vector space of signed measures supported on $X$ by $\M(X)$ and the cone of positive (signed) measures supported on $X$ by $\M^+(X)$. Moreover, for $\mu \in \M(\R^n)$, we use the notation $\mu|_X$ to denote the restricted measure $\mu(\cdotbox \cap X)$.  %
	
	Using the Riesz-Markov-Kakutani representation theorem, we will show in Section \ref{sect_weak_derivatives_signed_measures} that weak derivatives can be represented by signed measures on the boundary of $K_0$. 
	Formally, this statement reads as follows.
	\begin{proposition} \label{prop_weak_derivatives_are_signed_measures}
		A family $(K_t)_{t \in [0,1]}$ is weakly differentiable if and only if there exists a signed measure $\mu \in \M(\bd K_0)$ such that
		\begin{equation} \label{eq_prop_weak_derivatives_are_signed_measures}
			\left.\ddt \phi_f(t) \right|_{t=0}=\lim_{t \rightarrow 0^+} \frac{1}{t}\left(\int_{K_t} f(x) \dint x - \int_{K_0} f(x) \dint x\right) = \int_{\bd K_0} f \dint \mu
		\end{equation}
		for all $f \in C(\R^n)$.
	\end{proposition}

	\begin{figure}[ht]
		\begin{subfigure}[b]{.33\linewidth}
			\centering
			\begin{tikzpicture}%
	[x={(0.562037cm, -0.330669cm)},
	y={(0.827112cm, 0.224669cm)},
	z={(0.000024cm, 0.916614cm)},
	scale=1.050000,
	back/.style={dotted, thin},
	edge/.style={color=black},
	facet/.style={fill=white,fill opacity=0.25},
	vertex/.style={}]
	
	\coordinate (1.00000, -1.00000, -1.00000) at (1.00000, -1.00000, -1.00000);
	\coordinate (1.00000, 1.00000, -1.00000) at (1.00000, 1.00000, -1.00000);
	\coordinate (1.00000, 1.00000, 1.00000) at (1.00000, 1.00000, 1.00000);
	\coordinate (1.00000, -1.00000, 1.00000) at (1.00000, -1.00000, 1.00000);
	\coordinate (-1.00000, -1.00000, 1.00000) at (-1.00000, -1.00000, 1.00000);
	\coordinate (-1.00000, -1.00000, -1.00000) at (-1.00000, -1.00000, -1.00000);
	\coordinate (-1.00000, 1.00000, -1.00000) at (-1.00000, 1.00000, -1.00000);
	\coordinate (-1.00000, 1.00000, 1.00000) at (-1.00000, 1.00000, 1.00000);
	\draw[edge,back] (1.00000, 1.00000, -1.00000) -- (-1.00000, 1.00000, -1.00000);
	\draw[edge,back] (-1.00000, -1.00000, -1.00000) -- (-1.00000, 1.00000, -1.00000);
	\draw[edge,back] (-1.00000, 1.00000, -1.00000) -- (-1.00000, 1.00000, 1.00000);
	\fill[facet] (1.00000, -1.00000, 1.00000) -- (1.00000, -1.00000, -1.00000) -- (1.00000, 1.00000, -1.00000) -- (1.00000, 1.00000, 1.00000) -- cycle {};
	\fill[facet] (-1.00000, 1.00000, 1.00000) -- (1.00000, 1.00000, 1.00000) -- (1.00000, -1.00000, 1.00000) -- (-1.00000, -1.00000, 1.00000) -- cycle {};
	\fill[facet] (-1.00000, -1.00000, -1.00000) -- (1.00000, -1.00000, -1.00000) -- (1.00000, -1.00000, 1.00000) -- (-1.00000, -1.00000, 1.00000) -- cycle {};
	\draw[edge] (1.00000, -1.00000, -1.00000) -- (1.00000, 1.00000, -1.00000);
	\draw[edge] (1.00000, -1.00000, -1.00000) -- (1.00000, -1.00000, 1.00000);
	\draw[edge] (1.00000, -1.00000, -1.00000) -- (-1.00000, -1.00000, -1.00000);
	\draw[edge] (1.00000, 1.00000, -1.00000) -- (1.00000, 1.00000, 1.00000);
	\draw[edge] (1.00000, 1.00000, 1.00000) -- (1.00000, -1.00000, 1.00000);
	\draw[edge] (1.00000, 1.00000, 1.00000) -- (-1.00000, 1.00000, 1.00000);
	\draw[edge] (1.00000, -1.00000, 1.00000) -- (-1.00000, -1.00000, 1.00000);
	\draw[edge] (-1.00000, -1.00000, 1.00000) -- (-1.00000, -1.00000, -1.00000);
	\draw[edge] (-1.00000, -1.00000, 1.00000) -- (-1.00000, 1.00000, 1.00000);
\end{tikzpicture}
			
			\caption{unit cube $C_3=K_0$}
		\end{subfigure}
		\begin{subfigure}[b]{.33\linewidth}
			\centering
			\begin{tikzpicture}%
	[x={(0.562037cm, -0.330669cm)},
	y={(0.827112cm, 0.224669cm)},
	z={(0.000024cm, 0.916614cm)},
	scale=1.050000,
	back/.style={dotted, thin},
	edge/.style={color=black},
	facet/.style={fill=white,fill opacity=0.000000},
	vertex/.style={},
	back2/.style={dotted},
	edge2/.style={color=black},
	facet2/.style={fill=orange,fill opacity=0.200000},
	vertex2/.style={},
	back3/.style={dotted, thin},
	edge3/.style={color=black},
	facet3/.style={fill=cyan,fill opacity=0.20000},
	vertex3/.style={}]
	
	\draw[edge,back] (1.00000, 1.00000, -1.00000) -- (-1.00000, 1.00000, -1.00000);
	\draw[edge,back] (-1.00000, -1.00000, -1.00000) -- (-1.00000, 1.00000, -1.00000);
	\draw[edge,back] (-1.00000, 1.00000, -1.00000) -- (-1.00000, 1.00000, 1.00000);
	\node[vertex] at (-1.00000, 1.00000, -1.00000)     {};
	\fill[facet] (1.00000, -1.00000, 1.00000) -- (1.00000, -1.00000, -1.00000) -- (1.00000, 1.00000, -1.00000) -- (1.00000, 1.00000, 1.00000) -- cycle {};
	\fill[facet] (-1.00000, 1.00000, 1.00000) -- (1.00000, 1.00000, 1.00000) -- (1.00000, -1.00000, 1.00000) -- (-1.00000, -1.00000, 1.00000) -- cycle {};
	\fill[facet] (-1.00000, -1.00000, -1.00000) -- (1.00000, -1.00000, -1.00000) -- (1.00000, -1.00000, 1.00000) -- (-1.00000, -1.00000, 1.00000) -- cycle {};
	\draw[edge] (1.00000, -1.00000, -1.00000) -- (1.00000, 1.00000, -1.00000);
	\draw[edge] (1.00000, -1.00000, -1.00000) -- (1.00000, -1.00000, 1.00000);
	\draw[edge] (1.00000, -1.00000, -1.00000) -- (-1.00000, -1.00000, -1.00000);
	\draw[edge] (1.00000, 1.00000, -1.00000) -- (1.00000, 1.00000, 1.00000);
	\draw[edge] (1.00000, 1.00000, 1.00000) -- (1.00000, -1.00000, 1.00000);
	\draw[edge] (1.00000, 1.00000, 1.00000) -- (-1.00000, 1.00000, 1.00000);
	\draw[edge] (1.00000, -1.00000, 1.00000) -- (-1.00000, -1.00000, 1.00000);
	\draw[edge] (-1.00000, -1.00000, 1.00000) -- (-1.00000, -1.00000, -1.00000);
	\draw[edge] (-1.00000, -1.00000, 1.00000) -- (-1.00000, 1.00000, 1.00000);
	
	\fill[facet2] (0.72727, -0.95455, 0.77273) -- (1.00000, -1.00000, -1.00000) -- (1.00000, 1.00000, 1.00000) -- cycle {};
	\fill[facet2] (0.50000, -0.54545, 0.97727) -- (1.00000, 1.00000, 1.00000) -- (0.72727, -0.95455, 0.77273) -- cycle {};
	\fill[facet2] (0.50000, -0.54545, 0.97727) -- (-1.00000, -1.00000, 1.00000) -- (0.72727, -0.95455, 0.77273) -- cycle {};
	\fill[facet2] (0.50000, -0.54545, 0.97727) -- (-1.00000, -1.00000, 1.00000) -- (1.00000, 1.00000, 1.00000) -- cycle {};
	\fill[facet2] (0.52273, -1.00000, 0.58636) -- (-1.00000, -1.00000, 1.00000) -- (0.72727, -0.95455, 0.77273) -- cycle {};
	\fill[facet2] (0.52273, -1.00000, 0.58636) -- (1.00000, -1.00000, -1.00000) -- (0.72727, -0.95455, 0.77273) -- cycle {};
	\draw[edge,gray] (-1.00000, -1.00000, 1.00000) -- (1.00000, 1.00000, 1.00000);
	\draw[edge,back2] (-1.00000, -1.00000, 1.00000) -- (0.72727, -0.95455, 0.77273);
	\draw[edge,back2] (-1.00000, -1.00000, 1.00000) -- (0.50000, -0.54545, 0.97727);
	\draw[edge,back2] (-1.00000, -1.00000, 1.00000) -- (0.52273, -1.00000, 0.58636);
	\draw[edge,gray] (1.00000, -1.00000, -1.00000) -- (1.00000, 1.00000, 1.00000);
	\draw[edge,back2] (1.00000, -1.00000, -1.00000) -- (0.72727, -0.95455, 0.77273);
	\draw[edge,back2] (1.00000, -1.00000, -1.00000) -- (0.52273, -1.00000, 0.58636);
	\draw[edge,back2] (1.00000, 1.00000, 1.00000) -- (0.72727, -0.95455, 0.77273);
	\draw[edge,back2] (1.00000, 1.00000, 1.00000) -- (0.50000, -0.54545, 0.97727);
	\draw[edge,back2] (0.72727, -0.95455, 0.77273) -- (0.50000, -0.54545, 0.97727);
	\draw[edge,back2] (0.72727, -0.95455, 0.77273) -- (0.52273, -1.00000, 0.58636);
	
	\coordinate (-0.29545, -1.18182, -0.15909) at (-0.29545, -1.18182, -0.15909);
	\coordinate (-1.00000, -1.00000, -1.00000) at (-1.00000, -1.00000, -1.00000);
	\coordinate (-1.00000, -1.00000, 1.00000) at (-1.00000, -1.00000, 1.00000);
	\coordinate (1.00000, -1.00000, -1.00000) at (1.00000, -1.00000, -1.00000);
	\coordinate (0.52273, -1.00000, 0.58636) at (0.52273, -1.00000, 0.58636);
	\fill[facet3] (1.00000, -1.00000, -1.00000) -- (-0.29545, -1.18182, -0.15909) -- (-1.00000, -1.00000, -1.00000) -- cycle {};
	\fill[facet3] (-1.00000, -1.00000, 1.00000) -- (-0.29545, -1.18182, -0.15909) -- (-1.00000, -1.00000, -1.00000) -- cycle {};
	\fill[facet3] (0.52273, -1.00000, 0.58636) -- (-0.29545, -1.18182, -0.15909) -- (-1.00000, -1.00000, 1.00000) -- cycle {};
	\fill[facet3] (0.52273, -1.00000, 0.58636) -- (-0.29545, -1.18182, -0.15909) -- (1.00000, -1.00000, -1.00000) -- cycle {};
	\draw[edge] (-0.29545, -1.18182, -0.15909) -- (-1.00000, -1.00000, -1.00000);
	\draw[edge] (-0.29545, -1.18182, -0.15909) -- (-1.00000, -1.00000, 1.00000);
	\draw[edge] (-0.29545, -1.18182, -0.15909) -- (1.00000, -1.00000, -1.00000);
	\draw[edge] (-0.29545, -1.18182, -0.15909) -- (0.52273, -1.00000, 0.58636);
	\draw[edge] (-1.00000, -1.00000, 1.00000) -- (0.52273, -1.00000, 0.58636);
	\draw[edge] (1.00000, -1.00000, -1.00000) -- (0.52273, -1.00000, 0.58636);
\end{tikzpicture}
			
			\caption{perturbation of $C_3$}%
		\end{subfigure}
		\caption{Proposition \ref{prop_weak_derivatives_are_signed_measures} asserts that the family $(K_t)_{t \in [0,1]}$ is weakly differentiable if and only if the family of signed measures with densities $\tfrac{1}{t}[\One_{K_{t}}-\One_{K_0}]$, $t \in [0,1]$, \textit{converges weakly} (see Definition \ref{def_weak_convergence} below) to a signed measure on the boundary of $K_0$ as $t \rightarrow 0^+$.}
	\end{figure}
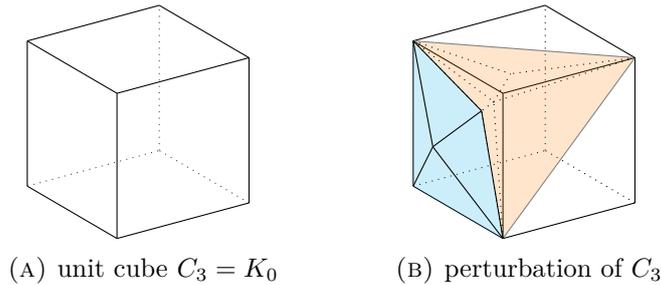

	Throughout the article, we will identify the weak derivative $K'_0$ with its representing measure $\mu \in \M(\bd K_0)$. %
	Proposition \ref{prop_weak_derivatives_are_signed_measures} shows that the weak derivative from Definition \ref{def_weak_derivative} coincides with a notion of weak derivatives of probability measures due to Pflug \cite{pfl96}. Such weak derivatives were studied for convex bodies by Weisshaupt in \cite{wei01} and Weisshaupt and Pflug in \cite{pw05}. For a broader discussion of different notions of derivatives of set-valued maps, we refer to \cite[Sect.~1]{khm07}.
	
	The following terminology reflects the fact that in Question \ref{qu_first_question} the weak derivatives themselves are the objects of primary interest, rather than the families that realize them.
	\begin{definition} \label{def_small_perturbation}
		We call a signed measure $\mu \in \M(\bd K)$ a \textit{small perturbation} of $K \in \Knn$ if there is a weakly differentiable family $(K_t)_{t \in [0,1]}$ with $K_0 = K$ and $K'_0(f)=\int_{\bd K} f \dint\mu$ for all $f \in C(\R^n)$. The small perturbation $\mu$ is called reversible if $-\mu$ is also a small perturbation of $K$. We denote the set of small perturbations of $K$ by $\W(K)$ and the set of reversible small perturbations of $K$ by $\W^\pm(K)$.
	\end{definition}
	
	Evidently, Definition \ref{def_small_perturbation} customizes the convenient but somewhat casual term ``small perturbation'' to our specific purpose of finding first-order conditions -- in other contexts, it might be appropriate to consider higher-order representatives of weakly differentiable families, or to apply a completely different notion of ``small perturbations''.
	
	In light of this terminology, Question \ref{qu_first_question} can be rephrased as follows:
	\begin{question} \label{qu_second_question}
		Given $K \in \Knn$, what are the sets $\W(K)$ and $\W^\pm(K)$?
	\end{question} 
	A signed measure $\nu\in \M(X)$ is called \textit{absolutely continuous} with respect to another signed measure $\mu\in \M(X)$ if one of the following two equivalent conditions (see \cite[Sect.~7.5]{kle20}) holds:
	\begin{enumerate}[label=\roman*)]
		\item for every Borel set $A \subset X$ with $\mu(A)=0$, we have $\nu(A)=0$;
		\item there exists a \textit{density}, i.e., a function $f \in L^1(X)$ with $\nu(A) = \int_A f \dint \mu$ for all  Borel sets $A \subset X$.
	\end{enumerate}
	Let $\mathcal{H}^{n-1} \in \M(\bd K)$ be the $(n-1)$-dimensional Hausdorff measure (which assigns to each Borel set $A \subset \bd K$ the corresponding surface area). For a function $f \colon K \rightarrow \R$, let $\mu_f\in \M(\bd K)$ be the signed measure that is given by $\mu_f(A) = \int_A f \dint \mathcal{H}^{n-1}$ for all Borel sets $A\subset \bd K$. In \cite{wei01}, Weisshaupt showed that
	\begin{equation}
		\W(K) \supset \{\mu_f \in \M(\bd K) \mid f\colon K \rightarrow \R \text{ is Lipschitz continuous and concave}\}.
	\end{equation}
	We will extend Weisshaupt's result in the special case that $K$ is a polytope, i.e., the convex hull of a finite set $V \subset \R^n$. We denote the spaces of polytopes and full-dimensional polytopes in $\R^n$ by $\Pn$ and $\Pnn$, respectively. For a comprehensive overview of the combinatorial theory of polytopes, we refer to \cite{zie07}.
	\begin{definition} \label{def_alpha_convexity}
		Let $\alpha\in (0,1)$ and $C \subset \R^n$ convex. A function $f \colon C \rightarrow [0,\infty)$ is called \textit{$\alpha$-convex} if $f^\alpha$ is convex. For $\alpha=1$, we drop the non-negativity requirement, i.e., a 1-convex function is simply a convex function $C \rightarrow \R$. Let $P \in \Pn$. For $\alpha \in (0,1]$, we denote by $\M^{\cvx}_\alpha(P)$ the space of all signed measures that are absolutely continuous with respect $\vol_{\dim P} \in \M(P)$ with an $\alpha$-convex density $f \colon \relint P \rightarrow \R$.
	\end{definition}
	We are now ready to state our main result, which gives a full characterization of the sets $W(P)$ and $\W^\pm(P)$ for $P \in \Pnn$. Here and in the following, $\Phi(P)$ is the set of proper faces of $P$ and $\Phi_{m}(P)$ is the set of faces of dimension $m$.%
	\begin{theorem} \label{thm_main_result}
		Let $P\in \Pnn$. The set of small perturbations of $P$ is given by the direct sum
		\begin{equation} \label{eq_thm_main_result}
			\W(P)=\sum_{F \in \Phi(P)} -\M^{\cvx}_{(n-\dim F)^{-1}}(F).
		\end{equation}
		This set is a convex cone.
	\end{theorem}
	
	Translated into the setting of the Reynolds transport theorem (see, e.g., \cite{rp22}), Theorem \ref{thm_main_result} asserts that the velocity function that describes the evolution of $\bd K_t$ has to be concave on the relative interior of each facet $F \subset P$. With regard to the lower-dimensional faces, Theorem \ref{thm_main_result} shows that the restriction of a small perturbation to a face $F \subset P$ with $\dim F < n-1$ is a negative measure, which has a density $f$ such that $-f$ is $(n-\dim F)^{-1}$-convex. This conclusion resembles the reverse Brunn-Minkowski inequality for coconvex sets \cite{kt14,sch18,fil17}.
	
	Theorem \ref{thm_main_result} leads to the following straightforward corollary.
	
	\begin{corollary} \label{cor_main_result}
		Let $P\in \Pnn$. The set of reversible small perturbations of $P$ is given by
		\begin{equation}
			\W^\pm(P)=\sum_{F \in \Phi_{n-1}(P)} [\M^{\cvx}_{1}(F)\cap -\M^{\cvx}_{1}(F)].
		\end{equation}
		This set is a vector space of dimension $\# \Phi_{n-1}(P) \cdot n$.
	\end{corollary}

	The rest of this text is organized as follows: In Section \ref{sect_weak_derivatives_signed_measures}, we discuss some functional analytic preliminaries and show that weak derivatives can be represented by signed measures. The four subsequent sections are devoted to the proof Theorem \ref{thm_main_result}. In Section \ref{sect_discrete_perturbations}, we construct polyhedral perturbations of a polytope, which are given by piecewise affine densities on the facets. In Section \ref{sect_wasserstein_distance}, we discuss the Wasserstein distance, which enables us to metrize weak convergence on compact sets. Using this tool, we show in Section \ref{sect_weak_limits_of discrete_perturbations} that the set of polyhedral perturbations is dense in the set that we claim to be $\W(P)$. To complete the proof of our main result, we show in Section \ref{sect_proof_of_main_result} that there are no other small perturbations of $P$. Finally, we discuss a necessary condition for polytopal maximizers of the isotropic constant in Section \ref{sect_application}.
	
	\section{Weak derivatives and signed measures} \label{sect_weak_derivatives_signed_measures}
	
	In this section, we study some basic properties of weakly differentiable families. Using the Riesz-Markov-Kakutani representation theorem, we show that weak derivatives correspond to signed measures on the boundary of $\bd K_0$.
	
	It is common to equip $\Kn$ with the \emph{Hausdorff metric}, which is defined by
	\begin{equation}
		\dH(K,L)\coloneqq \max \left\{\max_{x \in K} \dist(x,L),\,\max_{x \in L} \dist(x,K)\right\}.
	\end{equation}
	Another metric on $\Kn$ is defined in terms of the symmetric difference $K \triangle L \coloneqq K\setminus L \cup L \setminus K$. It is called the \textit{symmetric difference metric} and is given by
	\begin{equation}
		\dS(K,L)\coloneqq \vol(K \triangle L).
	\end{equation}
	It was shown in \cite{sw65} that $\dH$ and $\dS$ induce the same topology on $\Kn$.
	
	The following lemma shows that a weakly differentiable family $(K_t)_{t \in [0,1]}$ is continuous at $t=0$ with respect to both $\dH$ and $\dS$.
	
	\begin{lemma}\label{lemma_fast_convergence}
		Let $(K_t)_{t \in [0,1]}$ be weakly differentiable. 
		\begin{enumerate}[label=(\roman*)]
			\item \label{lemma_fast_convergence_1} We have \begin{equation}
				\limsup_{t \rightarrow 0} \frac{\vol(K_t \triangle K_0)}{t} < \infty.
			\end{equation}
		\end{enumerate}
		In particular, $K_t$ converges to $K_0$ in the Hausdorff metric as $t \rightarrow 0$, which has the following consequences:
		\begin{enumerate}[label=(\roman*)] \setcounter{enumi}{1}
			\item \label{lemma_fast_convergence_2} For every $\varepsilon>0$, there exists a $t_\varepsilon \in (0,1]$ such that $K_t \subset K_0+\varepsilon B_n$ for all $t \leq t_\varepsilon$, where $B_n$ denotes the centered Euclidean ball of radius 1.
			\item \label{lemma_fast_convergence_3} For every convex body $L \subset \interior K_0$, there exists a $t_L \in (0,1]$ such that $L \subset K_t$ for all $t \leq t_L$.
		\end{enumerate}
	\end{lemma}
	
	\begin{proof}
		We assume the negation of \ref{lemma_fast_convergence_1}, which is equivalent to the statement that there exists a decreasing sequence $(t_i)_{i \in \N}$ converging to zero with $\vol(K_{t_i} \triangle K_0) \geq 2^i \cdot t_i$ for all $i \in \N$. Without loss of generality, let $0 \in \interior K$. In the following, we will use the gauge function $\norm{\cdot}_{K_0}\colon \R^n \rightarrow \R$ of $K_0$, which is given by
		\begin{equation}
			\norm{x}_{K_0}\coloneqq \min\{t \geq 0 \mid x \in tK_0\}.
		\end{equation}
		We set
		\begin{equation}
			M_i \coloneqq \left\{x \in K_{t_i} \triangle K_0 \MID \norm{x}_{K_0} \notin [1-r_i,1+r_i]\right\} \quad \text{for } i \in \N.
		\end{equation}
		For every $i \in \N$, there exists an $r_i \in (0,1)$ with $\vol(M_i) \geq 2^{i-1}\cdot t_i$. Clearly, the values $r_i$ can be chosen such that $r_{i+1}<r_i$ for all $i\in\N$. The sequence $(r_i)_{i \in \N}$ induces a partition of $(0,1)$ into half-open subintervals
		\begin{equation}
			I_i \coloneqq \{s \in (0,1) \mid r_{i+1} \leq s < r_i \} \text{ for } i \geq 1; \quad \text{and } I_0 \coloneqq \{s \in (0,1) \mid r_1 \leq s\}.
		\end{equation}
		We define $f \colon [0,1] \rightarrow [0,\infty)$ by setting $f(r_i)\coloneqq \tfrac{1}{i}$, $f(0) \coloneqq 0$, $f(1) \coloneqq 1$ and interpola\-ting linearly on the intervals $I_i$. Since $f$ is continuous, the function $g \colon \R^n \rightarrow \R$ given by
		\begin{equation}
			g(x) \coloneqq \sgn(\norm{x}_{K_0}-1) \cdot f(\min\{|\norm{x}_{K_0}-1|,1\})
		\end{equation}
		is continuous as well. By construction, we have $g(x) \geq \frac{1}{i}$ on $M_i \cap K_0$ and $g(x) \leq -\frac{1}{i}$ on $M_{i} \setminus K_0$, leading to the estimate
		\begin{align}
			\frac{1}{t_i}\left(\int_{K_{t_i}} g(x) \dint x- \int_{K_0} g(x) \dint x\right) \geq \frac{1}{t_i \cdot i}\vol M_i\geq \frac{2^{i-1}}{i} \goestoinfty{i} \infty.
		\end{align}
		This shows that the function $\phi_g(t)$ from Definition \ref{def_weak_derivative} is not differentiable from the right at $t=0$.
		
		The claim of \ref{lemma_fast_convergence_2} is an immediate consequence of the definition of $\dH$, for \ref{lemma_fast_convergence_3} see \cite[Lem.~1.8.18]{sch13}.
	\end{proof}
	
	In order to state the Riesz-Markov-Kakutani representation theorem, we review some preliminaries from functional analysis. 
	
	Let $X \subset \R^n$. We denote by $C_b(X)$ the Banach space of bounded continuous functions $X \rightarrow \R$, equipped with the norm $\norm{f}_\infty \coloneqq \sup_{x \in X} |f(x)|$. Two subspaces of $C_b(X)$ are relevant for our purposes: 
	\begin{enumerate}[label=\roman*)]
		\item The subspace of $C_b(X)$ that contains all functions $f$ whose support
		\begin{equation}
			\supp f \coloneqq \cl \{x \in X \mid f(x)\neq 0\}
		\end{equation}
		is compact is denoted by $C_c(X)$. The normed space $C_c(x)$ is not complete.
		\item The closure of $C_c(x)$ is denoted by $C_0(X)$ and contains all functions $f$ that \textit{vanish at infinity} in the following sense: for every $\varepsilon >0$ there exists a compact set $K \subset X$ with $|f(x)| < \varepsilon$ for all $x \notin K$.
	\end{enumerate}
	The (topological) \textit{dual space} of $C_0(X)$ is denoted by $C_0(X)^*$ and consists of all continuous linear functionals $C_0(X) \rightarrow \R$. A linear functional $h \colon C_0(X) \rightarrow \R$ is continuous if and only if its operator norm
	\begin{equation}
		\norm{h}_* \coloneqq \sup \{|h(f)| \mid f \in C_0(X), \, \norm{f}_\infty \leq 1\}
	\end{equation}
	is finite.
	
	For $X \subset \R^n$, let $\mu \in \M(X)$ be a signed measure. Hahn's decomposition theorem \cite[Thm.~7.43]{kle20} asserts that there exists a set $X^+\subset X$ such that
	\begin{enumerate}[label=\roman*)]
		\item $\mu(A) \geq 0$ for all $A \subset X^+$; and
		\item $\mu(A) \leq 0$ for all $A \subset X^- \coloneqq X \setminus X^+$.
	\end{enumerate}
	The partition $(X^+,X^-)$ is called a \textit{Hahn decomposition} of $X$ with respect to $\mu$. Setting $\mu^+\coloneqq \mu(\cdotbox \cap X^+)$ and $\mu^-\coloneqq -\mu(\cdotbox \cap X^-)$, we have $\mu = \mu^+-\mu^-$. The pair $(\mu^+,\mu^-)$ is called the \textit{Jordan decomposition} of $\mu$; it does not depend on which particular Hahn decomposition is chosen (see \cite[Cor.~7.44]{kle20}). The positive measure $|\mu|\coloneqq \mu^++\mu^-$ is called the \textit{total variation} of $\mu$.
	\begin{definition}
		Let $\mu \in \M(X)$ be a signed measure with Jordan decomposition $(\mu^+,\mu^-)$. The expression
		\begin{align}
			\norm{\mu}_{\TV} \coloneqq |\mu|(X)= \mu^+(X)+\mu^-(X)
		\end{align}
		defines a norm on $\M(X)$ (see \cite[Cor.~7.45]{kle20}), which is called the \textit{total variation norm}.
	\end{definition}
	
	In the following, we will be mainly concerned with signed measures that are absolutely continuous with respect to the volume (Lebesgue measure) on some affine subspace of $\R^n$. For a Borel set $A \subset \R^n$ and $f \in L^1(A)$, we use the notation $[f]_{\M(A)}$ to denote the signed measure on $A$ that is absolutely continuous with respect to $\vol_{\dim \aff A}$ with density $f$. If $f$ is only defined on a subset of $A$, we use the trivial extension of $f$ to $A$ as a density. Moreover, if $\aff A = \R^n$, we simply write $[f]_{\M}$. For such signed measures, we have
	\begin{equation}
		\norm{[f]_{\M(A)}}_{\TV}=\norm{f}_{L^1(A)}.
	\end{equation}
	
	We now state a special case of the Riesz-Markov-Kakutani representation theorem; see for example \cite[Thm.~2.26]{els18} or \cite[Thm.~14.14]{ab99}.
	
	\begin{theorem}[Riesz-Markov-Kakutani] \label{thm_rmk}
		Let $X \subset \R^n$ be open or closed. For $\mu \in \M(X)$, let $I_\mu \colon C_0(X)\rightarrow \R$ be given by $I_\mu(f) = \int_X f \dint \mu$. The map $\mu \mapsto I_\mu$ is an isometric isomorphism between the normed spaces $(\M(X), \norm{\cdot}_{\TV})$ and $(C_0(X)^*, \norm{\cdot}_*)$.
	\end{theorem}
	
	Since $C_c(X)$ is dense in $C_0(X)$, a continuous linear functional $C_0(X) \rightarrow \R$ is uniquely determined by its values on $C_c(X)$. In particular, $C_c(X)$ separates points in $C_0(X)^*$.
	
	\begin{remark}
		Theorem \ref{thm_rmk} is usually stated for more general topological spaces $X$, namely locally compact Hausdorff spaces. %
		For arbitrary locally compact Hausdorff spaces, the map $\mu \mapsto I_{\mu}$ is not injective, unless we demand that the representing measure $\mu$ be \textit{regular}, i.e., we require that for every Borel set $A \subset X$ and $\varepsilon > 0$, there exist a compact set $K$ and an open set $U$ such that $K \subset A \subset U$ and $\mu(U\setminus K)< \varepsilon$ \cite[Def.~2.21]{els18}. In our setting, we can drop the requirement since by \cite[(13.7.9)]{die70} and \cite[(13.7.7)]{die70} every positive finite measure on $\R^n$ is regular, and hence, by virtue of the Jordan decomposition, every signed measure on $\R^n$ is regular.
	\end{remark}

	\begin{definition} \label{def_weak_convergence}
		Let $X \subset \R^n$ be open or closed and let $(\mu_i)_{i \in \N}$ be a sequence in $\M(X)$. We say that $(\mu_i)_{i \in \N}$ \textit{converges weakly} to $\mu \in \M(X)$ and write $\mu_i \warrow \mu$ or $\wlim_{i \rightarrow \infty}\mu_i=\mu$ if
		\begin{equation} \label{eq_def_weak_convergence}
			\int_X f \dint \mu_i \goestoinfty{i} \int_X f \dint \mu
		\end{equation}
		holds for all $f \in C_b(X)$.%
	\end{definition}
	
	It follows from Theorem \ref{thm_rmk} that weak limits are unique.
	
	\begin{remark} \label{rem_two_notions}
		Two additional notions of convergence are obtained if we require that condition \eqref{eq_def_weak_convergence} holds for all $f \in C_c(X)$ or for all $f \in C_0(X)$, respectively. For the purposes of this article, the differences between these notions are not important since we can, considering small perturbations of a given convex body $K \subset \R^n$, restrict our attention to signed measures in $\M(K+B_n)$. Some consequences concerning the differentiability of parametric families of measures encountered in the case of non-compact support are discussed in \cite[Ex.~2.12--Ex.~2.15]{wei09}. For a general comparison of topologies and convergence on spaces of measures, including the notions of convergence induced by the functions in $C_b(X)$ and $C_c(X)$, respectively, we refer to \cite[Sect.~13.20:~Prob.~1,2]{die70}. %
	\end{remark}

	\begin{proof}[Proof of Proposition \ref{prop_weak_derivatives_are_signed_measures}]
		We only have to show that every weak derivative can be represented by a signed measure $\mu \in \M(\bd K_0)$. In light of Theorem \ref{thm_rmk}, this amounts to showing that $K_0'$ can be written as
		\begin{equation}
			K_0' = h \circ (f \mapsto f|_{\bd K_0}),
		\end{equation}
		where $h$ is a continuous linear functional $C(\bd K_0) \rightarrow \R$. Let 
		\begin{equation}
			C \coloneqq \limsup_{t \rightarrow 0} \frac{\norm{[\One_{K_t}-\One_{K_0}]_{\M}}_{\TV}}{t}.
		\end{equation}
		By Lemma \ref{lemma_fast_convergence}(i), we have $C < \infty$.
		Let $f \in C(\R^n)$. We have
		\begin{equation}
			\left|\!\left.\ddt \phi_f(t) \right|_{t=0}\right|=\lim_{t \rightarrow 0^+} \left|\frac{1}{t}\left(\int_{K_t} f(x) \dint x - \int_{K_0} f(x) \dint x\right)\right| \leq C \cdot \norm{f|_{K_0+B_n}}_\infty,
		\end{equation}
		where the last step follows from Lemma \ref{lemma_fast_convergence}\ref{lemma_fast_convergence_2}. This shows that $K_0'$ can be written as
		\begin{equation}
			K_0' = h \circ (f \mapsto f|_{K_0+B_n}),
		\end{equation}
		where $h$ is a continuous linear functional $C(K_0+B_n) \rightarrow \R$. By Theorem \ref{thm_rmk}, there exists a $\mu \in \M(K_0+B_n)$ with $K_0'(f)=\int_{K_0+B_n} f \dint \mu$ for all $f \in C(\R^n)$. We consider an arbitrary function $f \in C_c(\R^n \setminus \bd K_0)$. Since $\supp f$ is compact, it satisfies $\dist(\supp f, \bd K_0)>0$. We obtain that the trivial extension $\overline{f}$ given by
		\begin{equation}
			\overline{f}(x)=\begin{cases}
				f(x) & \text{if } x\notin \bd K_0,\\
				0 & \text{if } x\in \bd K_0,
			\end{cases}
		\end{equation}
		is continuous, and, using Lemma \ref{lemma_fast_convergence}\ref{lemma_fast_convergence_2}--\ref{lemma_fast_convergence_3}, that $\int [\One_{K_t}-\One_{K_0}] \overline{f}(x) \dint x =0$ holds for all $t$ smaller than some $t_f >0$. Hence, we have
		\begin{equation}
			\int f \dint \mu|_{\R^n \setminus \bd K_0} = \int \overline{f} \dint \mu = \lim_{t \rightarrow 0^+} \frac{1}{t} \int [\One_{K_t}-\One_{K_0}] \overline{f}(x) \dint x = 0.
		\end{equation}
		Since $C_c(\R^n \setminus \bd K_0)$ separates points in $C_0(\R^n \setminus \bd K_0)^*$, this shows that $\mu|_{\R^n \setminus \bd K_0}$ is the zero measure. It follows that $\mu$ is supported on $\bd K_0$.
	\end{proof}

	\section{Polyhedral perturbations of polytopes} \label{sect_discrete_perturbations}
	
	This section begins the constructive part of the proof of Theorem \ref{thm_main_result}. The basic building blocks of our construction are \textit{polyhedral perturbations}, which are described by densities on the boundary of a polytope $P$ whose restrictions to $\relint F$ are piecewise affine and concave for each facet $F \subset P$. Using an elementary construction, we show that every polyhedral perturbation can be realized as a weak derivative.
	
	An \textit{affine functional} on an affine subspace $H \subset \R^n$ is an affine function $H \rightarrow \R$. %
	
	\begin{definition}
		Let $K \subset \R^n$ be convex. A concave function $f \colon K \rightarrow \R$ is called \textit{polyhedral} if it is of the form
		\begin{equation} \label{eq_def_pca}
			f(x) = \min_{h \in \mathcal{H}}h(x)
		\end{equation}
		for a finite family $\mathcal{H}$ of affine functionals $\aff(K)\rightarrow \R$. Given a polyhedral concave function $f \colon K \rightarrow \R$, there exists a unique family $\mathcal{H}$ that provides an \textit{irredundant} description of $f$ in the sense that for each $h \in \mathcal{H}$ there is an $x \in K$ with
		\begin{equation}
			h(x) < \min_{g \in \mathcal{H} \setminus \{h\}} g(x).
		\end{equation}
		Assuming that $\mathcal{H}$ is irredundant, the expression \eqref{eq_def_pca} can be read as describing a polyhedral function $\overline{f}\colon \aff(K) \rightarrow \R$. We call $\overline{f}$ the \textit{canonical extension} of $f$.
	\end{definition}
	
	The main result of this section reads as follows.
	
	\begin{proposition} \label{prop_discrete_perturbation}
		Let $P \in \Pnn$. Moreover, let $f\colon \bd P \rightarrow \R$ be a measurable function with the property that the restriction $f|_{ \relint F}$ is concave and polyhedral for each facet $F \subset P$. Then we have
		\begin{equation} \label{eq_prop_discrete_perturbation}
			\sum_{F \in \Phi_{n-1}(P)}\left[f|_{F}\right]_{\M(F)} \in \W(P).
		\end{equation}
		We call a small perturbation of the form \eqref{eq_prop_discrete_perturbation} a \emph{polyhedral perturbation} and denote the set of polyhedral perturbations of $P$ by $\W^{\triangleleft}(P)$.
	\end{proposition}
	
	The geometric idea behind the following construction is to choose a generic direction $v \in \sphere$, introduce a redundant constraint for every affine functional in the description of $f$, and shift the corresponding hyperplanes in such a way that the velocities in direction $v$ are proportional to the values of $f$ at the corresponding boundary points of $P$.
	
	\begin{construction} \label{construction_discrete_perturbation}
		Let $P$ and $f$ be as in Proposition \ref{prop_discrete_perturbation} and assume that $P$ is given by
		\begin{equation}
			P = \{x \in \R^n \mid \scpr{u_i,x} \leq b_i \text{ for } i \in [m]\}
		\end{equation}
		for $u_1,\dots,u_m \in \sphere$ and $b_1,\dots,b_m \in \R$. We assume that this description is irredundant and denote the facet with outer normal vector $u_i$ by $F_i$, $i \in [m]$. Let $v \in \sphere$ be \textit{generic} with respect to $P$ in the sense that $\scpr{u_i,v} \neq 0 \text{ for } i \in [m]$. We set
		\begin{equation}
			I^+ \coloneqq \{i \in [m] \mid \scpr{u_i,v}>0\}\quad \text{and} \quad I^- \coloneqq \{i \in [m] \mid \scpr{u_i,v}<0\}.
		\end{equation}
		Because $v$ is generic, we have $I^-=[m]\setminus I^+$.
		
		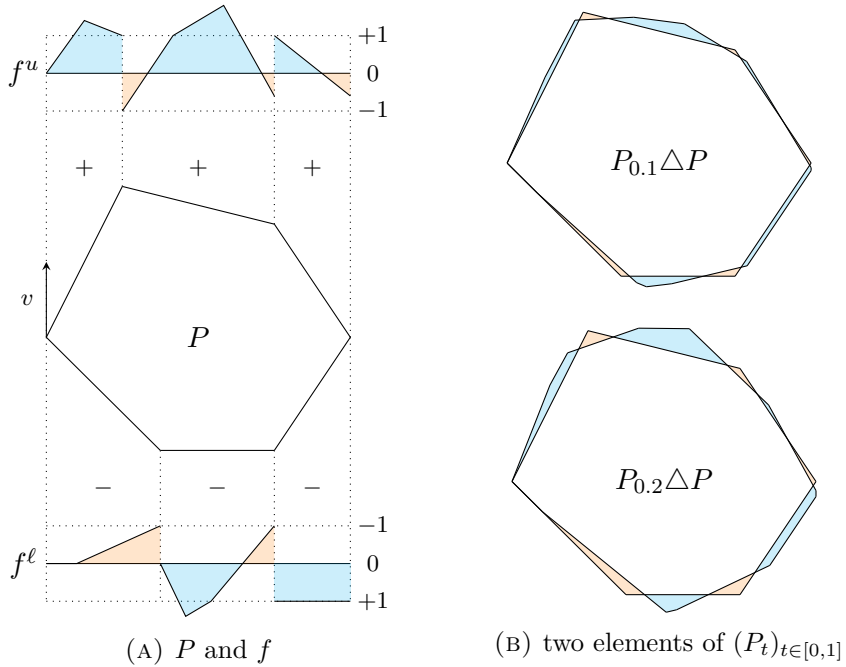
\begin{figure}[ht]
			\begin{subfigure}[c]{.4\linewidth}
				\centering
				\begin{tikzpicture}%
	[scale=1.000000,
	axis/.style={dotted, thin},
	axis2/.style={loosely dotted, thin},
	edge/.style={color=black},
	red/.style={color=red},
	orange/.style={fill=orange,fill opacity=0.200000},
	cyan/.style={fill=cyan,fill opacity=0.200000},
	vertex/.style={}]

	\draw[-stealth] (-2,0) -- (-2,1);
	
	\node at (-1.5,2.25) {\footnotesize$+$};
	\node at (0,2.25) {\footnotesize$+$};
	\node at (1.5,2.25) {\footnotesize$+$};
	
	\node at (-1.25,-2) {\footnotesize$-$};
	\node at (0.25,-2) {\footnotesize$-$};
	\node at (1.5,-2) {\footnotesize$-$};
	
	\draw[edge] (-2,0) -- (-1,2);
	\draw[edge,axis] (-2,0) -- (-2,4);
	\draw[edge,axis] (-2,0) -- (-2,-3.5);
	
	\draw[edge] (-1,2) -- (1,1.5);
	\draw[edge,axis] (-1,2) -- (-1,4);
	
	\draw[edge] (1,1.5) -- (2,0);
	\draw[edge,axis] (1,1.5) -- (1,4);
	
	\draw[edge] (2,0) -- (1,-1.5);
	\draw[edge,axis] (2,0) -- (2,4);
	\draw[edge,axis] (2,0) -- (2,-3.5);
	
	\draw[edge] (1,-1.5) -- (-0.5,-1.5);
	\draw[edge,axis] (1,-1.5) -- (1,-3.5);
	
	\draw[edge] (-0.5,-1.5) -- (-2,0);
	\draw[edge,axis] (-0.5,-1.5) -- (-0.5,-3.5);
	
	\node at (0,0) {$P$};
	
	\node at (-2.25,0.5) {\footnotesize $v$};

	\node at (-2.3,3.5) {$f^u$};
	\draw[edge,axis] (-2,4) -- (2,4);
	\node at (2.3,4) {\footnotesize $+1$};
	\draw[edge] (-2,3.5) -- (2,3.5);
	\node at (2.3,3.5) {\footnotesize $0$};
	\draw[edge,axis] (-2,3) -- (2,3);
	\node at (2.3,3) {\footnotesize $-1$};
	
	\node at (-2.3,-3) {$f^\ell$};
	\draw[edge,axis] (-2,-3.5) -- (2,-3.5);
	\node at (2.3,-3.5) {\footnotesize $+1$};
	\draw[edge] (-2,-3) -- (2,-3);
	\node at (2.3,-3) {\footnotesize $0$};
	\draw[edge,axis] (-2,-2.5) -- (2,-2.5);
	\node at (2.3,-2.5) {\footnotesize $-1$};
	
	\fill[cyan] (-2, 3.5+0) -- (-1.5, 3.5+0.7) -- (-1, 3.5+0.5) -- (-1, 3.5+0) -- cycle {};
	\draw[edge] (-2, 3.5+0) -- (-1.5, 3.5+0.7) -- (-1, 3.5+0.5);
	
	\fill[orange] (-1, 3.5-0.5) -- (-2/3, 3.5+0) -- (-1, 3.5+0) -- cycle {};
	\fill[cyan] (-2/3, 3.5+0) -- (-1/3, 3.5+0.5) -- (1/3, 3.5+0.9) -- (5/6, 3.5+0) -- cycle {};
	\fill[orange] (5/6, 3.5+0) -- (1, 3.5-0.3) -- (1, 3.5+0) -- cycle {};
	\draw[edge] (-1, 3.5-0.5) -- (-1/3, 3.5+0.5) -- (1/3, 3.5+0.9) -- (1, 3.5-0.3);
	
	\fill[cyan] (1, 3.5+0.5) -- (1+5/8, 3.5+0) -- (1, 3.5+0) -- cycle {};
	\fill[orange] (1+5/8, 3.5+0) -- (2, 3.5-0.3) -- (2, 3.5+0) -- cycle {};
	\draw[edge] (1, 3.5+0.5) -- (2, 3.5-0.3);
	
	\fill[orange] (-1.6, -3) -- (-0.5, -3) -- (-0.5, -3+.5) -- cycle {};
	\draw[edge] (-2, -3) -- (-1.6, -3) -- (-0.5, -3+0.5);
	
	\fill[cyan] (-0.5, -3) -- (-1/6, -3-.7) -- (1/6, -3-.5) -- (1/6+5/12, -3) -- cycle {};
	\fill[orange] (1/6+5/12, -3) -- (1,-3) -- (1,-3+.5) -- cycle {};
	\draw[edge] (-0.5, -3) -- (-1/6, -3-.7) -- (1/6, -3-.5) -- (1, -3+.5);
	
	\fill[cyan] (1, -3) -- (1, -3-.5) -- (2, -3-.5) -- (2, -3) -- cycle {};
	\draw[edge] (1, -3-.5) -- (2, -3-.5);
	
\end{tikzpicture}
				
				\caption{$P$ and $f$ \label{subfig_P_and_f}}
			\end{subfigure}
			\begin{subfigure}[c]{0.4\linewidth}
				\centering
				\begin{tikzpicture}%
	[scale=1.000000,
	back/.style={loosely dotted, thin},
	edge/.style={color=black},
	facet2/.style={fill=cyan,fill opacity=0.200000},
	facet1/.style={fill=orange,fill opacity=0.200000},
	facet/.style={fill=white,fill opacity=1.00000},
	vertex/.style={}]
	
	\coordinate (2.00000, 0.00000) at (2.00000, 0.00000);
	\coordinate (1.00000, -1.50000) at (1.00000, -1.50000);
	\coordinate (-0.50000, -1.50000) at (-0.50000, -1.50000);
	\coordinate (1.00000, 1.50000) at (1.00000, 1.50000);
	\coordinate (-2.00000, 0.00000) at (-2.00000, 0.00000);
	\coordinate (-1.00000, 2.00000) at (-1.00000, 2.00000);
	\fill[facet1] (-1.00000, 2.00000) -- (1.00000, 1.50000) -- (2.00000, 0.00000) -- (1.00000, -1.50000) -- (-0.50000, -1.50000) -- (-2.00000, 0.00000) -- cycle {};
	\fill[facet2] (-1.50000, 1.14000) -- (-1.10695, 1.89465) -- (-0.33333, 1.93333) -- (0.33333, 1.84667) -- (1.15238, 1.34705) -- (2.00000, -0.06000) -- (2.00000, -0.10000) -- (1.15873, -1.36190) -- (0.16667, -1.60000) -- (-0.16667, -1.64000) -- (-0.29554, -1.58587) -- (-1.60000, -0.40000) -- (-2.00000, 0.00000) -- cycle {};
	
	\fill[facet] (-2.00000, 0.00000) -- (-1.60000, -0.40000) -- (-0.39000, -1.50000) -- (0.58333, -1.50000) -- (1.07937, -1.38095) -- (1.98101, -0.02848) -- (1.62500, 0.56250) -- (1.06742, 1.39888) -- (0.83333, 1.54167) -- (-0.66667, 1.91667) -- (-1.05128, 1.89744) -- cycle {};
	
	\draw[edge] (-1.10695, 1.89465) -- (-0.33333, 1.93333);
	\draw[edge] (-1.10695, 1.89465) -- (-1.50000, 1.14000);
	\draw[edge] (2.00000, -0.10000) -- (1.15873, -1.36190);
	\draw[edge] (2.00000, -0.10000) -- (2.00000, -0.06000);
	\draw[edge] (1.15873, -1.36190) -- (0.16667, -1.60000);
	\draw[edge] (-1.60000, -0.40000) -- (-0.29554, -1.58587);
	\draw[edge] (-1.60000, -0.40000) -- (-2.00000, 0.00000);
	\draw[edge] (-0.33333, 1.93333) -- (0.33333, 1.84667);
	\draw[edge] (-0.29554, -1.58587) -- (-0.16667, -1.64000);
	\draw[edge] (0.16667, -1.60000) -- (-0.16667, -1.64000);
	\draw[edge] (0.33333, 1.84667) -- (1.15238, 1.34705);
	\draw[edge] (1.15238, 1.34705) -- (2.00000, -0.06000);
	\draw[edge] (-2.00000, 0.00000) -- (-1.50000, 1.14000);
	
	\draw[edge] (2.00000, 0.00000) -- (1.00000, -1.50000);
	\draw[edge] (2.00000, 0.00000) -- (1.00000, 1.50000);
	\draw[edge] (1.00000, -1.50000) -- (-0.50000, -1.50000);
	\draw[edge] (-0.50000, -1.50000) -- (-2.00000, 0.00000);
	\draw[edge] (1.00000, 1.50000) -- (-1.00000, 2.00000);
	\draw[edge] (-2.00000, 0.00000) -- (-1.00000, 2.00000);
	
	\node at (0,0) {$P_{0.1} \triangle P$};
	
\end{tikzpicture}
				\vspace*{0.5cm}
				
				\begin{tikzpicture}%
	[scale=1.000000,
	back/.style={loosely dotted, thin},
	edge/.style={color=black},
	facet2/.style={fill=cyan,fill opacity=0.200000},
	facet1/.style={fill=orange,fill opacity=0.200000},
	facet/.style={fill=white,fill opacity=1.00000},
	vertex/.style={}]
	
	\coordinate (2.00000, 0.00000) at (2.00000, 0.00000);
	\coordinate (1.00000, -1.50000) at (1.00000, -1.50000);
	\coordinate (-0.50000, -1.50000) at (-0.50000, -1.50000);
	\coordinate (1.00000, 1.50000) at (1.00000, 1.50000);
	\coordinate (-2.00000, 0.00000) at (-2.00000, 0.00000);
	\coordinate (-1.00000, 2.00000) at (-1.00000, 2.00000);
	\fill[facet1] (-1.00000, 2.00000) -- (1.00000, 1.50000) -- (2.00000, 0.00000) -- (1.00000, -1.50000) -- (-0.50000, -1.50000) -- (-2.00000, 0.00000) -- cycle {};
	\fill[facet2] (-1.50000, 1.28000) -- (-1.26846, 1.70604) -- (-0.33333, 2.03333) -- (0.33333, 2.02667) -- (1.37647, 1.01482) -- (2.00000, -0.12000) -- (2.00000, -0.20000) -- (1.39216, -1.11176) -- (0.16667, -1.70000) -- (0.02921, -1.73299) -- (-1.60000, -0.40000) -- (-2.00000, 0.00000) -- cycle {};
	
	\fill[facet] (-1.12121, 1.75758) -- (-0.66667, 1.91667) -- (0.83333, 1.54167) -- (1.22642, 1.16038) -- (1.62500, 0.56250) -- (1.96386, -0.05422) -- (1.19608, -1.20588) -- (0.58333, -1.50000) -- (-0.25556, -1.50000) -- (-1.60000, -0.40000) -- (-2.00000, 0.00000) -- cycle {};
	\draw[edge] (-1.26846, 1.70604) -- (-0.33333, 2.03333);
	\draw[edge] (-1.26846, 1.70604) -- (-1.50000, 1.28000);
	\draw[edge] (2.00000, -0.20000) -- (1.39216, -1.11176);
	\draw[edge] (2.00000, -0.20000) -- (2.00000, -0.12000);
	\draw[edge] (1.39216, -1.11176) -- (0.16667, -1.70000);
	\draw[edge] (-0.33333, 2.03333) -- (0.33333, 2.02667);
	\draw[edge] (-1.60000, -0.40000) -- (0.02921, -1.73299);
	\draw[edge] (-1.60000, -0.40000) -- (-2.00000, 0.00000);
	\draw[edge] (0.16667, -1.70000) -- (0.02921, -1.73299);
	\draw[edge] (0.33333, 2.02667) -- (1.37647, 1.01482);
	\draw[edge] (1.37647, 1.01482) -- (2.00000, -0.12000);
	\draw[edge] (-2.00000, 0.00000) -- (-1.50000, 1.28000);
	
	\draw[edge] (2.00000, 0.00000) -- (1.00000, -1.50000);
	\draw[edge] (2.00000, 0.00000) -- (1.00000, 1.50000);
	\draw[edge] (1.00000, -1.50000) -- (-0.50000, -1.50000);
	\draw[edge] (-0.50000, -1.50000) -- (-2.00000, 0.00000);
	\draw[edge] (1.00000, 1.50000) -- (-1.00000, 2.00000);
	\draw[edge] (-2.00000, 0.00000) -- (-1.00000, 2.00000);

	\node at (0,0) {$P_{0.2} \triangle P$};
\end{tikzpicture}
				
				\caption{two elements of $(P_t)_{t \in [0,1]}$}%
		\end{subfigure}
		\caption{An illustration of Construction \ref{construction_discrete_perturbation}. Subfigure (A) refers to the functions $f^u$ and $f^\ell$ introduced in Lemma \ref{lemma_discrete_perturbation} below. The sets $P_{0.1} \setminus P$ and $P_{0.2} \setminus P$ are shaded in cyan, whereas the sets $P \setminus P_{0.1}$ and $P \setminus P_{0.2}$ are shaded in orange.}
	\end{figure}

		Denoting the orthogonal projection onto $v^\bot \coloneqq \{x \in \R^n \mid \scpr{v,x}=0\}$ by $\pi_{v^\bot}$, we define
		\begin{equation}
			h_i\colon \R^n \rightarrow \R, \quad  h_i(x)\coloneqq\frac{-\scpr{u_i, \pi_{v^\bot}(x)}+b_i}{\scpr{u_i,v}} \quad \text{for } i\in [m].
		\end{equation}
		Then $P$ can be written as
		\begin{equation}
			P= \{ x \in \R^n \mid \scpr{v,x} \leq h_i(x) \text{ for } i \in I^+,\ \scpr{v,x} \geq h_i(x) \text{ for } i \in I^-\}. 
		\end{equation}
		Let $i \in I$. Since $v$ is generic, the restriction $\pi_{v^\bot}|_{\aff F_i}$ is a bijection $\aff F_i \rightarrow v^\bot$. We denote this bijection by $\pi_i$ and define a function $f_i \colon v^\bot \rightarrow \R$ via
		\begin{equation} \label{eq_def_f_i}
			f_i(x) \coloneqq \frac{\overline{f|_{\relint F_i}}(\pi_i^{-1}(x))}{|\scpr{u_i,v}|},
		\end{equation}
		where $\overline{f|_{\relint F_i}}\colon \aff F_i \rightarrow \R$ is the canonical extension of $f|_{\relint F_i}$. 
		
		Now we define two functions $u,\ell \colon v^\bot \times [0,1]\rightarrow \R$ by
		\begin{equation} \label{eq_def_u_ell}
			u(x,t)\coloneqq \min_{i \in I^+} \left[h_i(x)+ tf_i(x)\right] \eqand
			\ell(x,t)\coloneqq\max_{i \in I^-} \left[h_i(x)- t f_i(x)\right].
		\end{equation}
		By construction, $u(\cdotbox,t)$ and $-\ell(\cdotbox,t)$ are concave and polyhedral for all $t\in [0,1]$. Therefore,
		\begin{equation} \label{eq_def_tilde_P}
			\tilde{P}_t \coloneqq \left\{ x \in \pi_{v^\bot}^{-1}(P) \MID \ell(\pi_{v^\bot}(x),t) \leq \scpr{v,x} \leq u(\pi_{v^\bot}(x),t)\right\}
		\end{equation}
		is either a polytope or the empty set for $t \in [0,1]$. We obtain a family of full-dimensional polytopes $(P_t)_{t \in [0,1]}$ that satisfies $P_0=P$ by setting
		\begin{equation}
			\renewcommand\qedsymbol{$\lozenge$}
			\pushQED{\qed} 
			P_t \coloneqq \begin{cases}
				\tilde{P}_t & \text{if }\dim \tilde{P}_t = n,\\
				P  & \text{otherwise}.
			\end{cases} \qedhere
			\popQED
		\end{equation}
	\end{construction}

	\begin{remark} \label{rem_construction_discrete_perturbation_volume_bound}
		For later reference we note that Construction \ref{construction_discrete_perturbation}  has the following property: For any $t \in [0,1]$, the volume of $P_t \setminus P$ is given by
		\begin{align}
			 \int_{\pi_{v^\bot}(P)} \max\left\{u(x,t)-\min_{i \in I^+} h_i(x),0\right\} \dint x + \int_{\pi_{v^\bot}(P)} \max\left\{\max_{i \in I^-} h_i(x)-\ell(x,t),0\right\}\dint x
		\end{align}
		and hence, using \eqref{eq_def_u_ell}, satisfies the bound
		\begin{align}
			\vol(P_t \setminus P) \leq \sum_{i \in [m]}\int_{\pi_{v^\bot}(F_i)} \max \{tf_i(x),0\} \dint x= t\cdot \left(\sum_{F \in \Phi_{n-1}(P)}\left[f|_{F}\right]_{\M(F)}\right)^+.
		\end{align}
	\end{remark}

	\begin{example}
		Let $P \subset \R^3$ be a simplicial polytope and let $v \in \mathbb{S}^2$ be generic with respect to $P$. Moreover, let $F \subset P$ be a facet. For the sake of simplicity, we assume that $\pi_{v^\bot}(F) \subset \relint(\pi_{v^\bot}(P))$. Performing Construction \ref{construction_discrete_perturbation} with the following functions $f$, we obtain various basic perturbations of $P$.
		\begin{enumerate}[label=(\roman*)]
			\item \label{example_shift} \textit{Shifting $F$:} If $f = \One_{F}$, then the family $(P_t)_{t \in [0,1]}$ corresponds to a continuous outward parallel shift of the facet $F$ with constant velocity.
			\item \label{example_hinge} \textit{Hinging $F$:} Let $E \subset F$ be an edge and set $f \coloneqq \One_F \cdot \dist(\cdotbox, E)$. Then the family $(P_t)_{t \in [0,1]}$ corresponds to hinging the facet $F$ outward around the edge $E$, where the dihedral angle between $F$ and the hinged facet $\tilde{F}_t$ at time $t \in [0,1]$ is equal to $\arctan(t)$.
			\item \label{example_stack} \textit{Stacking a pyramid onto $F$:} If $f \coloneqq \One_F \cdot \dist(\cdotbox, \relbd F)$, then the family $(P_t)_{t \in [0,1]}$ corresponds to the operation of stacking a pyramid $Q$ onto the facet $F$, where the apex of $Q$ moves outwards with constant velocity. 
		\end{enumerate}
		Recalling Corollary \ref{cor_main_result}, we note that \ref{example_shift} and \ref{example_hinge} describe reversible small perturbations, whereas \ref{example_stack} is not reversible.
	\end{example}
	
	Proposition \ref{prop_discrete_perturbation} follows from the following claim.
	
	\begin{claim} \label{claim_discrete_perturbations}
		The family $(P_t)_{t \in [0,1]}$ from Construction \ref{construction_discrete_perturbation} is weakly differentiable with
		\begin{equation}
			P_0' = \sum_{F \in \Phi_{n-1}(P)}\left[f|_{F}\right]_{\M(F)}.
		\end{equation}
	\end{claim}
	
	Before we proceed to show Claim \ref{claim_discrete_perturbations}, we prove some auxiliary statements.
	
	\begin{lemma} \label{lemma_discrete_perturbation}
		Let $P$ and $f$ be as in Proposition \ref{prop_discrete_perturbation} and $f_i$, $i\in [m]$, be as in Construction \ref{construction_discrete_perturbation}. We define two functions $v^\bot \rightarrow \R$ by
		\begin{equation}
			f^u \coloneqq \sum_{i \in I^+} f_i \cdot \One_{\pi_{v^\bot}(\relint F_i)} \eqand f^\ell \coloneqq \sum_{i \in I^-} f_i \cdot \One_{\pi_{v^\bot}(\relint F_i)}.
		\end{equation}
		The functions $f^u$ and $f^\ell$ and the functions $u$, $\ell$ from Construction \ref{construction_discrete_perturbation} have the following properties:
		\begin{enumerate}[label=(\roman*)]
			\item \label{lemma_discrete_perturbation:epsilon_x} For almost all $x \in \pi_{v^\bot}(P)$ there is an $\varepsilon_x>0$ such that \begin{equation}
				u(x,t)=u(x,0)+t f^u(x) \eqand \ell(x,t)=\ell(x,0)-t f^\ell(x)
			\end{equation}
			for $t\in [0,\varepsilon_x]$.
			\item \label{lemma_discrete_perturbation:speed_limit} Given $P$ and $f$, there is a $C>0$ such that \begin{equation}
				\max\{|u(x,t)-u(x,0)|,|\ell(x,t)-\ell(x,0)|\}\leq Ct
			\end{equation}
			for $t \in [0,1]$ and $x \in \pi_{v^\bot}(P)$.
		\end{enumerate}
	\end{lemma}
	
	\begin{proof}
		For \ref{lemma_discrete_perturbation:epsilon_x}, we first observe that $N^+\coloneqq\pi_{v^\bot}(P) \setminus \left(\bigcup_{i \in I^+} \pi_{v^\bot}(\relint F_i)\right)$ and $N^-\coloneqq\pi_{v^\bot}(P) \setminus \left(\bigcup_{i \in I^-} \pi_{v^\bot}(\relint F_i)\right)$ satisfy $\volt{n-1}(N^+)=\volt{n-1}(N^-)=0$. Let $x \in \pi_{v^\bot}(P)\setminus (N^+ \cup N^-)$. There is a unique $i \in I^+$ such that $x \in \pi_{v^\bot}(\relint F_i)$. Let $f_{\min}\coloneqq \min_{j \in I^+}f_j(x)$ and $h_{\min} \coloneqq \min_{j \in I^+\setminus \{i\}}h_j(x)$. Since $x \in \pi_{v^\bot}(\relint F_i)$, we have $h_i(x) < h_{\min}$. If $f_{\min}=f_i(x)$, then $u(x,t)=h_i(x)+t f_i(x)$ holds for all $t \geq 0$ and we can choose $\varepsilon^+_x>0$ arbitrarily. Otherwise, the choice $\varepsilon^+_x \coloneqq \frac{h_{\min}-h_i(x)}{f_i(x)-f_{\min}}$ implies
		\begin{equation}
			h_i(x)+t f_i(x) \leq h_{\min}+t f_{\min} \leq \min_{j \in I^+\setminus \{i\}} h_j(x)+f_j(x) \quad \text{for } t \in [0,\varepsilon^+_x].
		\end{equation}
		In both cases, we have
		\begin{equation}
			u(x,t)=h_i(x)+t f_i(x)=u(x,0)+t f^u(x)  \quad \text{for } t \in [0,\varepsilon^+_x].
		\end{equation}
		By a similar argument, we find an $\varepsilon^-_x>0$ such that $\ell(x,t)=\ell(x,0)-t f^\ell(x)$ for $t \in [0,\varepsilon^-_x]$. Setting $\varepsilon_x\coloneqq \min \{\varepsilon^+_x,\varepsilon^-_x\}$, we obtain the claim.
		
		For \ref{lemma_discrete_perturbation:speed_limit}, we set $C \coloneqq \max_{i \in I}\max_{x \in \pi_{v^\bot(P)}}|f_i(x)|$. Since $I$ is finite and $\pi_{v^\bot}(P)$ is compact, we have $C<\infty$. Let $ x \in \pi_{v^\bot}(P)$ and $t \in [0,1]$. We have
		\begin{equation}
			u(x,t) \geq \min_{i \in I^+}h_i(x)+t \cdot\min_{i \in I^+}f_i(x)\geq u(x,0) - Ct
		\end{equation}
		and
		\begin{equation}
			u(x,t) \leq \min_{i \in I^+} h_i(x)+ t \cdot \max_{i \in I^+}f_i(x) \leq u(x,0)+Ct,
		\end{equation}
		implying $|u(x,t)-u(x,0)|\leq Ct$. The bound on $|\ell(x,t)-\ell(x,0)|$ is shown similarly.
	\end{proof}
	
	We are now ready to verify that the family of polytopes $(P_t)_{t \in [0,1]}$ from Construction \ref{construction_discrete_perturbation} realizes the small perturbation given by the density $f$.
	
	\begin{proof}[Proof of Claim \ref{claim_discrete_perturbations}]
		Let $g \in C(\R^n)$ be arbitrary. Setting
		\begin{equation}
			G(t)\coloneqq\frac{1}{t} \left[\int_{P_t}g \dint \volt{n} - \int_{P}g \dint \volt{n}\right], \quad t\in [0,1],
		\end{equation}
		and 
		\begin{equation}
			\mu \coloneqq \sum_{F \in \Phi_{n-1}(P)}\left[f|_{ \relint F}\right]_{\M(\relint F)},
		\end{equation}
		our goal is to show that $\lim_{t \rightarrow 0^+} G(t) = \int g \dint \mu$.  Without loss of generality, we can assume that $v$ is the $n$-th standard basis vector $e_n$ and that $g$ is bounded (as the claim only depends on the values of $g$ on, say, $P+B_n$). Identifying $\R^n$ with $\R^{n-1} \times \R$, we write $g(x,y)\coloneqq g((x_1,\dots,x_{n-1},y)^\transp)$. Recalling \eqref{eq_def_f_i}, we observe that
		\begin{equation}
			\int g \dint \mu = \int_{\pi_{v^\bot}(P)} \left(g(x,u(x,0))\cdot f^u(x)+g(x,\ell(x,0))\cdot f^\ell(x)\right) \mathrm{d} x.
		\end{equation}
		Since $\pi_{v^\bot}(P_t) \subset \pi_{v^\bot}(P)$, we can decompose $G(t)$ into two parts, $G(t)=G_1(t)+G_2(t)$, where
		\begin{align}
			G_1(t)%
			\coloneqq \int_{\pi_{v^\bot}(P_t)\cap \pi_{v^\bot}(P)} \frac{1}{t}\left(\int_{\ell(x,t)}^{u(x,t)} g(x,y) \dint y -\int_{\ell(x,0)}^{u(x,0)} g(x,y) \dint y\right) \mathrm{d} x
		\end{align}
		and
		\begin{align}
			G_2(t)\coloneqq\int_{\pi_{v^\bot}(P) \setminus \pi_{v^\bot}(P_t)} \frac{1}{t}\int_{\ell(x,0)}^{u(x,0)} g(x,y) \dint y \dint x.
		\end{align}
		We first consider $G_2(t)$. On the set $\pi_{v^\bot}(P) \setminus \pi_{v^\bot}(P_t)$ we have $u(x,t)<\ell(x,t)$, which implies, using Lemma \ref{lemma_discrete_perturbation}\ref{lemma_discrete_perturbation:speed_limit},
		\begin{align}
			u(x,0)-\ell(x,0) &< \ell(x,t)-\ell(x,0)-u(x,t)+u(x,0)\leq 2Ct \quad\text{for } x\in \pi_{v^\bot}(P) \setminus \pi_{v^\bot}(P_t).
		\end{align}
		It follows that, for all $(x,t) \in (\pi_{v^\bot}(P) \setminus \pi_{v^\bot}(P_t)) \times [0,1]$,
		\begin{equation}
			\frac{1}{t}\int_{\ell(x,0)}^{u(x,0)} g(x,y) \dint y \leq 2 C  \cdot \norm{g}_\infty.
		\end{equation}
		
		Therefore, $x \mapsto 2 C  \cdot \norm{g}_\infty\cdot \One_{\pi_{v^\bot}(P)}(x)$ is a dominating integrable function and we can apply the Lebesgue dominated convergence theorem to obtain
		\begin{align}
			\lim_{t \rightarrow 0}G_2(t)&=\lim_{t \rightarrow 0} \int \underbrace{\One_{\pi_{v^\bot}(P) \setminus \pi_{v^\bot}(P_t)}(x)}_{\rightarrow 0 \text{ pointwise}} \frac{1}{t}\int_{\ell(x,0)}^{u(x,0)} g(x,y) \dint y \dint x=0.
		\end{align}
		Turning to $G_1(t)$, we observe that it can be rewritten as	
		\begin{align}
			G_1(t)&= \int_{\pi_{v^\bot}(P_t)\cap \pi_{v^\bot}(P)} \frac{1}{t}\left(\int_{u(x,0)}^{u(x,t)} g(x,y) \dint y -\int_{\ell(x,0)}^{\ell(x,t)} g(x,y) \dint y\right) \mathrm{d} x.
		\end{align}
		Again using Lemma \ref{lemma_discrete_perturbation}\ref{lemma_discrete_perturbation:speed_limit}, we have 
		\begin{equation}
			\left|\frac{1}{t}\int_{u(x,0)}^{u(x,t)} g(x,y) \dint y\right|\leq \frac{1}{t} \cdot |u(x,t)-u(x,0)| \cdot \norm{g}_\infty\leq C \cdot \norm{g}_\infty
		\end{equation}
		and a similar bound for $\left|\frac{1}{t}\int_{\ell(x,0)}^{\ell(x,t)} g(x,y) \dint y\right|$. Therefore, $x \mapsto 2 C  \cdot \norm{g}_\infty\cdot \One_{\pi_{v^\bot}(P)}(x)$ is again a dominating integrable function. Applying Lemma \ref{lemma_discrete_perturbation}\ref{lemma_discrete_perturbation:epsilon_x}, the Lebesgue dominated convergence theorem and the fundamental theorem of calculus, we obtain
		\begin{align}
			\lim_{t \rightarrow 0}G_1(t)&=\lim_{t \rightarrow 0} \int_{\pi_{v^\bot}(P)} \frac{1}{t}\left(\int_{u(x,0)}^{u(x,0)+t f^u(x)} g(x,y) \dint y -\int_{\ell(x,0)}^{\ell(x,0)-t f^\ell(x)} g(x,y) \dint y\right) \mathrm{d} x\\
			&= \int_{\pi_{v^\bot}(P)} \left(g(x,u(x,0))\cdot f^u(x)+g(x,\ell(x,0))\cdot f^\ell(x)\right) \dint x= \int g \dint \mu.
			\qedhere
		\end{align}
	\end{proof}

	\section{The Wasserstein distance} \label{sect_wasserstein_distance}
	
	We say that a sequence $(\mu_i)_{i\in\N}$ in $\M(\R^n)$ is \textit{bounded} if $\sup_{n \in \N}\norm{\mu_n}_{\TV}<\infty$. The goal of this section is to metrize weak convergence of bounded sequences in $\M(X)$ for a compact set $X \subset \R^n$. As a byproduct of the metrization, we will obtain the following result.
	
	\begin{proposition} \label{prop_set_of_perturbations_weakly_closed}
		Let $P \in \Pnn$. The sequential closure of the set $\W^\triangleleft(P) \subset \M(\bd P)$ in the topology of weak convergence is contained in $\W(P)$, i.e., if $(\mu_i)_{i \in \N}$ is a sequence of polyhedral perturbations of $P$ and $\mu_i \warrow \mu$, then $\mu \in \W(P)$.
	\end{proposition}

	\begin{remark} \label{rem_bounded_subseq}
		Let $X \subset \R^n$ be compact. Then the $\norm{\cdot}_{\TV}$-unit ball of $\M(X)$, which we denote by $B_{\M(X)}$, is compact and metrizable (and hence also sequentially compact) with respect to the vague topology, which is induced by the family of all functions $g \in C_c(\R^n)$ via $\mu \mapsto \int g \dint \mu$. By Remark \ref{rem_two_notions} and the compactness of $X$, the vague topology coincides with the topology induced by all functions $g \in C_b(\R^n)$ via $\mu \mapsto \int g \dint \mu$. We consider a bounded sequence $(\mu_i)_{i \in \N}$ in $\M(X)$. Then there exists a signed measure $\mu \in \M(X)$ and a subsequence $(\mu_{k(i)})_{i \in \N}$ with $\mu_{k(i)} \warrow \mu$. If $\mu_i \in \M^+(X)$ for all $i \in \N$, then $\mu \in \M^+(X)$. This is a direct consequence of \cite[(13.4.1)]{die70} and \cite[(13.4.2)(ii)]{die70} adapted to real measures (see \cite[Sect.~13.2]{die70}), using that $B_{M(X)}$ is for compact $X$ closed in the vague topology on the space of all real measures on $\R^n$. That the real measures in the sense of \cite[Sect.~13.2]{die70} include the signed measures in the sense of this article follows from \cite[(13.7.9)]{die70} and \cite[(13.3.6)]{die70}. 
	\end{remark}
	
	Instead of using the mentioned metrizability result ``off the shelf'', this section uses tools from optimal transport theory to develop a metrization of weak convergence %
	on the unit ball of $\M(X)$. The chosen metrization is based on the \textit{Wasserstein norm} from \cite{prt21}. In addition to proving Proposition \ref{prop_set_of_perturbations_weakly_closed}, this metrization provides a useful tool for analyzing weak convergence of concrete sequences in $\M(\R^n)$.	
	
	We first recall some basic notions from optimal transport theory. For a comprehensive introduction to this subject, we refer to \cite{vil09}. Let $X_1,X_2\subset \R^n$ and $\mu \in \M^+(X_1 \times X_2)$. The measure $\mu$ is said to have \textit{marginals} $\mu_1 \in \M^+(X)$ and $\mu_2 \in \M^+(X_2)$ if $(\pi_1)_* \mu= \mu_1$ and $(\pi_2)_* \mu = \mu_2$, where $\pi_i$ denotes the projection onto the $i$-th coordinate and $(\pi_i)_* \mu$ is the \textit{pushforward} of $\mu$, given by
	\begin{equation} \label{eq_def_pushforward}
		(\pi_i)_* \mu (A) = \mu(\pi_i^{-1}(A)) \quad \text{for all Borel sets } A \subset X_i.
	\end{equation}
	More generally, we define the pushforward $f_*\mu \in \M(Y)$ for any measurable map $f\colon X\rightarrow Y$ and a signed measure $\mu \in \M(X)$ as in \eqref{eq_def_pushforward}. 
	
	Let $X\subset \R^n$ and let $\mu_1,\mu_2 \in \M^+(X)$ with $\mu_1(X)=\mu_2(X)$. The \textit{Wasserstein distance} (of order 1) between $\mu_1$ and $\mu_2$ is given by	
	\begin{equation} \label{def_wasserstein_distance}
		W(\mu_1,\mu_2) \coloneqq \inf_\tau \int_{X \times X} \norm{x-y}_2 \dint \tau(x,y),
	\end{equation}
	where the infimum is taken over all $\tau \in \M^+(X \times X)$ with marginals $\mu_1$ and $\mu_2$. A measure $\tau$ as in \eqref{def_wasserstein_distance} is called a \textit{transference plan}. Every measurable function $f \colon X \rightarrow X$ with $f_* \mu_1 = \mu_2$ induces a transference plan $(\id,f)_*\mu_1$, where $(\id,f)$ denotes the map $x \mapsto (x,f(x))$.
	
	The following definition is a special case of \cite[Def.~11]{prt21}.
	
	\begin{definition} \label{def_generalized_wasserstein_distance}
		Let $\mu, \nu \in \M^+(\R^n)$. The \textit{generalized Wasserstein distance} of $\mu$ and $\nu$ is defined by
		\begin{equation} \label{eq_def_generalized_wasserstein_distance}
			\overline{W}(\mu,\nu)= \inf_{\tilde{\mu},\tilde{\nu}}\big( \norm{\mu-\tilde{\mu}}_{\TV}+\norm{\nu-\tilde{\nu}}_{\TV}+W(\tilde{\mu},\tilde{\nu})\big),
		\end{equation}
		where the infimum is taken over all $\tilde{\mu},\tilde{\nu}\in \M^+(\R^n)$ with $\tilde{\mu}(\R^n)=\tilde{\nu}(\R^n)$.
	\end{definition}
	
	As shown in \cite[Prop.~1]{pr14}, the infimum in \eqref{eq_def_generalized_wasserstein_distance} is attained and does not change if we impose the additional condition that $\tilde{\mu}(A) \leq \mu(A)$ and $\tilde{\nu}(A) \leq \nu(A)$ holds for all Borel sets $A \subset \R^n$.
	
	We say that a sequence $(\mu_i)_{i\in\N} $ in $M(\R^n)$ is \textit{tight} if for every $\varepsilon>0$ there exists a compact set $K \subset \R^n$ with $\sup_{i \in \N}|\mu_i|(\R^n \setminus K)<\varepsilon$. The significance of the generalized Wasserstein distance lies in the fact that it metrizes weak convergence of tight sequences of positive measures. The following result is taken from \cite[Thm.~3]{pr14}.
	
	\begin{proposition} \label{prop_wasserstein_distance_metrizes_positive_weak_convergence}
		Let $(\mu_i)_{i \in \N}$ be a sequence in $\M^+(\R^n)$ and $\mu \in \M^+(\R^n)$. We have
		\begin{equation}
			\overline{W}(\mu_i,\mu) \underset{i \rightarrow \infty}{\longrightarrow} 0 \quad \text{if and only if} \quad \mu_i \warrow \mu \text{ and } (\mu_i)_{i \in \N} \text{ is tight.}
		\end{equation}
	\end{proposition}
	
	We now turn to signed measures.
	
	\begin{definition}
		Let $\mu\in \M(X)$ with Jordan decomposition $(\mu^+,\mu^-)$. Following \cite[Def.~17]{prt21}, we define the \textit{Wasserstein norm} by
		\begin{equation}
			\norm{\mu}_{\Wt} \coloneqq \overline{W}(\mu^+,\mu^-).
		\end{equation}
	\end{definition}
	
	It is shown in \cite[Prop.~21]{prt21} that $\norm{\cdot}_{\Wt}$ is indeed a norm, and it is easy to see that $\norm{\mu}_{\Wt} \leq \norm{\mu}_{\TV}$ holds for all $\mu \in \M(X)$. The following sequential compactness result will be used extensively below, especially in Section \ref{sect_proof_of_main_result}.

	The utility of $\norm{\cdot}_{\Wt}$ for our purposes lies in the following fact.
	
	\begin{proposition} \label{prop_wasserstein_metrizes_weak_convergence}
		Let $X\subset \R^n$ be compact. The Wasserstein norm metrizes weak convergence of bounded sequences in $\M(X)$, i.e., for a bounded sequence $(\mu_i)_{i\in\N}$ in $\M(X)$ we have
		\begin{equation}
			\norm{\mu_i}_{\Wt} \underset{i \rightarrow \infty}{\longrightarrow} 0 \quad \text{if and only if} \quad \mu_i \warrow 0.
		\end{equation}
	\end{proposition}
	
	\begin{proof}
		We begin with the forward implication. Let $f \in C(X)$ and $\mu \in \M(X)$ with $\norm{\mu}_{\Wt} < \varepsilon$ for some $\varepsilon>0$. Then there exists a transference plan $\tau \in \M^+(X \times X)$ with $\tau(X\times X) \leq \mu^+(X)$ and
		\begin{equation}
			\norm{\mu}_{\Wt}= \norm{\mu^+-(\pi_1)_* \tau}_{\TV}+\norm{\mu^--(\pi_2)_* \tau}_{\TV}+\int_{X \times X} \norm{x-y}_2 \dint \tau(x,y).
		\end{equation} 
		Since $\norm{\mu}_{\Wt}< \varepsilon$, we have $\int_{X \times X} \norm{x-y}_2 \dint \tau(x,y)< \varepsilon$ and hence
		\begin{equation} \label{eq_prop_wasserstein_metrizes_weak_convergence_eq1}
			\tau\left(\{(x,y) \in X \times X \mid \norm{x-y}_2 \geq \sqrt{\varepsilon}\}\right) \leq \sqrt{\varepsilon}.
		\end{equation}
		Writing $\mu$ as $\mu^+-(\pi_1)_* \tau -\mu^-+(\pi_2)_* \tau +(\pi_1)_* \tau-(\pi_2)_* \tau$, we get
		\begin{align}
			\left|\int f \dint \mu\right| &\leq \norm{f}_\infty \left(\norm{\mu^+-(\pi_1)_* \tau}_{\TV}+\norm{\mu^--(\pi_2)_* \tau}_{\TV}\right)+\left|\int f \dint [(\pi_1)_* \tau-(\pi_2)_* \tau]\right|\\
			&\leq  \norm{f}_\infty \cdot \varepsilon+\left|\int f(x) - f(y) \dint \tau(x,y)\right|.
		\end{align}
		To bound the integral on the right-hand side, we first observe that
		\begin{align}
			\left|\int_{\{(x,y) \in X \times X \mid \norm{x-y}_2 < \sqrt{\varepsilon}\}} f(x) - f(y) \dint \tau(x,y)\right|\leq \sup_{\substack{x,y \in X\\ \norm{x-y}_2 < \sqrt{\varepsilon}}}|f(x)-f(y)| \cdot \tau(X\times X).
		\end{align}
		Since $X$ is compact, $f$ is uniformly continuous. Therefore, the function
		\begin{equation}
			\omega\colon [0,\diam X] \rightarrow [0,\infty), \quad t \mapsto \sup \big\{|f(x)-f(y)| \mid x,y \in X, \, \norm{x-y}_2 \leq t \big\}
		\end{equation}
		satisfies $\lim_{t \rightarrow 0^+}\omega(t)=\omega(0)=0$. Recalling \eqref{eq_prop_wasserstein_metrizes_weak_convergence_eq1}, we have
		\begin{align}
			\left|\int_{\{(x,y) \in X \times X \mid \norm{x-y}_2 \geq \sqrt{\varepsilon}\}} f(x) - f(y) \dint \tau(x,y)\right| \leq \sup_{x,y \in X} |f(x)-f(y)| \cdot \sqrt{\varepsilon}.
		\end{align}
		Putting the estimates together, we have
		\begin{equation}
			\left|\int f \dint \mu\right| \leq \norm{f}_\infty \cdot \varepsilon + \omega(\sqrt{\varepsilon}) \cdot \tau(X \times X) + 2 \cdot \norm{f}_\infty \cdot \sqrt{\varepsilon} \goestosth{\varepsilon}{0} 0,
		\end{equation}
		which shows the forward implication.
		
		For the other implication, we assume that $\mu_i \warrow 0$. By Remark \ref{rem_bounded_subseq}, there exists a subsequence of $\mu_{k(i)}$ whose Jordan decomposition satisfies ${\mu_{k(i)}}^+ \warrow \mu_+$ and ${\mu_{k(i)}}^- \warrow \mu_-$ for certain $\mu_+,\mu_- \in \M^+(X)$. Because $\mu_i \warrow 0$, we have
		\begin{equation}
			\int f \dint \mu_+ - \int f \dint \mu_- = \lim_{i \rightarrow \infty} \int f \dint \mu_{k(i)}=0
		\end{equation}
		for all $f \in C(X)$, which implies $\mu_+=\mu_-$. As $X$ is compact, both sequences $({\mu_{k(i)}}^+)_{i \in \N}$ and $({\mu_{k(i)}}^-)_{i \in \N}$ are tight, so Proposition \ref{prop_wasserstein_distance_metrizes_positive_weak_convergence} implies that \begin{equation}
			\overline{W}({\mu_{k(i)}}^+,\mu_+)\underset{i \rightarrow \infty}{\longrightarrow} 0 \quad \text{and} \quad \overline{W}({\mu_{k(i)}}^-,\mu_-)\underset{i \rightarrow \infty}{\longrightarrow} 0.
		\end{equation}
		The triangle inequality yields
		\begin{equation}
			\norm{\mu_{k(i)}}_{\Wt}=\overline{W}({\mu_{k(i)}}^+,{\mu_{k(i)}}^-)\leq\overline{W}({\mu_{k(i)}}^+,\mu_+)+\overline{W}(\mu_+,{\mu_{k(i)}}^-)\underset{i \rightarrow \infty}{\longrightarrow} 0. \qedhere
		\end{equation}
	\end{proof}
	
	\begin{remark} \label{rem_wasserstein_metrizes_weak_convergence}
		For $x \in [-1,1]$, let $\delta_x$ be the Dirac measure supported on $\{x\}$. The sequence $(\sqrt{i}(\delta_{1/i}-\delta_{-1/i}))_{i \in \N}$ shows that the forward implication in Proposition \ref{prop_wasserstein_metrizes_weak_convergence} fails without the assumption that $(\mu_i)_{i \in \N}$ is bounded. On the other hand, if we assume that the sequence $(\mu_i)_{i \in \N}$ is weakly convergent, then its boundedness follows from the uniform boundedness principle \cite[Thm.~14.1]{con07}.%
	\end{remark}

	We conclude this section with the proof of the result from its beginning.
	
	\begin{proof}[Proof of Proposition \ref{prop_set_of_perturbations_weakly_closed}]
		Let $v \in \sphere$ be generic with respect to $P$. For each $i \in \N$, let $(P_t^i)_{t \in [0,1]}$ be the family from Construction \ref{construction_discrete_perturbation} that realizes $\mu_i \in \W^\triangleleft(P)$. We fix $i \in \N$. By Remark \ref{rem_construction_discrete_perturbation_volume_bound}, we have $\frac{1}{t} \vol(P^i_t \setminus P) \leq \mu_i^+(\R^n)$ for all $t \in (0,1]$. Since $(P_t^i)_{t \in [0,1]}$ is weakly differentiable, we also have
		\begin{equation}
			\frac{1}{t} \left[\vol(P \setminus P^i_t)-\vol(P^i_t \setminus P)\right] \goestosth{t}{0} \mu_i^-(\R^n)-\mu_i^+(\R^n).
		\end{equation}
		Combining these two facts, we obtain that
		\begin{equation}
			\limsup_{t \rightarrow 0}\frac{1}{t} \left[\vol(P^i_t \setminus P)+\vol(P \setminus P^i_t)\right] \leq \mu_i^+(\R^n)+\mu_i^-(\R^n) = \norm{\mu_i}_{\TV}.
		\end{equation}
		Hence, by Lemma \ref{lemma_fast_convergence}\ref{lemma_fast_convergence_1} and Proposition \ref{prop_wasserstein_metrizes_weak_convergence}, there exists an $r_i>0$ such that 
		\begin{equation}
			\norm{\mu_i-\tfrac{1}{t}\cdot[\One_{P_t^i}-\One_P]_{\M}}_{\Wt}<\frac{1}{i} \quad \text{and} \quad \norm{\mu_i-\tfrac{1}{t}\cdot[\One_{P_t^i}-\One_P]_{\M}}_{\TV} \leq 2\norm{\mu_i}_{\TV} \quad \text{for all }t \in [0,r_i].
		\end{equation}
		Clearly, the values $r_i$ can be chosen such that $r_{i+1} \leq \tfrac{1}{i}$ for all $i \in \N$. For $t \in (0,1]$, we set $i(t) \coloneqq \max \{i \in \N \mid r_i \geq t\}$ and $P_t \coloneqq P_t^{i(t)}$. Finally, we set $P_0 \coloneqq P$. By Remark \ref{rem_wasserstein_metrizes_weak_convergence}, we have $\sup_{i \in \N}\norm{\mu_i}_{\TV} \leq C$ for some $C>0$, which implies that
		\begin{equation}
			\limsup_{t \rightarrow 0}\NORM{\mu_{i(t)}-\tfrac{1}{t}\cdot[\One_{P_t^{i(t)}}-\One_P]_{\M}}_{\TV} \leq 2\norm{\mu_{i(t)}}_{\TV} \leq 2C.
		\end{equation}
		Using Proposition \ref{prop_wasserstein_metrizes_weak_convergence} in combination with the triangle inequality for $\norm{\cdot}_{\W}$, it follows that $(P_t)_{t \in [0,1]}$ realizes the small perturbation $\mu$.
	\end{proof}

	\section{Weak limits of polyhedral perturbations} \label{sect_weak_limits_of discrete_perturbations}
	
	This section completes the constructive part of the proof of Theorem \ref{thm_main_result}. Using an approximation by polyhedral perturbations, we will show the first inclusion of Theorem \ref{thm_main_result}, which reads as follows.
	
	\begin{proposition} \label{prop_main_result_first_inclusion}
		Let $P\in \Pnn$. The set of small perturbations of $P$ satisfies
		\begin{equation}
			\W(P) \supset \sum_{F \in \Phi(P)} -\M^{\cvx}_{(n-\dim F)^{-1}}(F).
		\end{equation}
	\end{proposition}
	
	The following theorem is a classical result in convex analysis \cite[Thm.~10.8]{roc72}.
	
	\begin{theorem} \label{thm_rockafellar}
		Let $C \subset \R^n$ be relatively open and convex and let $(f_i)_{i \in \N}$ be a sequence of convex functions $C \rightarrow \R$. Let $X \subset C$ be dense in $C$, i.e., $C \subset \cl X$. If $(f_i)_{i \in \N}$ converges pointwise on $X$, then $(f_i)_{i \in \N}$ converges pointwise on $C$ to a convex function $f \colon C \rightarrow \R$, and the convergence is uniform on each closed bounded subset of $C$.
	\end{theorem}

	A convex function $f$ defined on a convex set $K\subset \R^n$ is called \textit{polyhedral} if $-f$ is polyhedral. The first step of our approximation argument is contained in the following lemma.
	
	\begin{lemma}\label{lemma_pac_approximation}
		Let $P \in \Pn$. For $\alpha \in (0,1]$, let $f \colon \relint P \rightarrow \R$ be an integrable $\alpha$-convex function. For every $\varepsilon>0$, there exists a polyhedral convex function $g\colon \relint P \rightarrow \R$ with
		\begin{equation}
			\norm{[f-g^{\frac{1}{\alpha}}]_{\M(P)}}_{\TV}=\int_{\interior P} |f(x)-g^{\frac{1}{\alpha}}(x)| \dint x < \varepsilon.
		\end{equation}
		In particular, we have $\norm{[f-g^{\frac{1}{\alpha}}]_{\M(P)}}_{\Wt}<\varepsilon$.
	\end{lemma}
	
	\begin{proof}
		Let $\varepsilon>0$. Without loss of generality, we assume that $P$ is full-dimensional. By assumption, $f^\alpha$ is convex. Let $(x_i)_{i \in \N}$ be a sequence whose image is dense in $\interior P$. For every $i \in \N$, let $h_i\colon \interior P \rightarrow \R$ be an affine functional with $h_i(x_i)=f^\alpha(x_i)$ and $h_i \leq f^\alpha$. We define $g_i$ by
		\begin{equation}
			g_i(x) \coloneqq \begin{cases}
				\max_{j \in [i]}h_j(x) & \text{if } \alpha =1,\\
				\max\left\{\max_{j \in [i]}h_j(x),0\right\} & \text{if } \alpha <1.
			\end{cases}
		\end{equation}
		The function $g_1^{1/\alpha}\colon \interior P \rightarrow \R$ is bounded and hence integrable. Therefore, there exists a compact set $Q \subset \interior P$ with 
		\begin{equation}
			\int_{(\interior P) \setminus Q} |f-g_i^{1/\alpha}(x)| \dint x \leq \int_{(\interior P) \setminus Q} |f-g_1^{1/\alpha}(x)| \dint x < \frac{\varepsilon}{2} \quad \text{for all } i\in \N.
		\end{equation}
		Since $(g_i)_{i \in \N}$ converges pointwise to $f^\alpha$ on the dense set $\{x_i\}_{i \in \N} \subset P$, Theorem \ref{thm_rockafellar} asserts that $(g_i)_{i \in \N}$ converges uniformly to $f^\alpha$ on the compact set $Q$. Since the sequence $(g_i)_{i \in \N}$ is monotonically increasing and bounded from above by the continuous function $f^\alpha$, there exists an interval $[a,b]$ (if $\alpha<1$, we can choose $a \geq 0$) with
		\begin{equation}
			f^\alpha(Q) \cup g_i(Q) \subset [a,b]  \quad \text{for all } i\in \N.
		\end{equation}
		Because the function $x \mapsto x^{\frac{1}{\alpha}}$ is Lipschitz continuous on $[a,b]$, the sequence $(g_i^{1/\alpha})_{i \in \N}$ converges uniformly to $f$ on $Q$. This implies the existence of an $N\in \N$ with
		\begin{equation}
			\int_{Q} |f(x)-g_i^{1/\alpha}(x)| \dint x < \frac{\varepsilon}{2} \quad \text{for all } i \geq N. \qedhere
		\end{equation}
	\end{proof}
	
	The following lemma allows us to approximate a measure with an ${(n-\dim F)^{-1}}$-convex density on a lower-dimensional face $F$ of $P$ by an iterative argument over the face lattice.
	
	\begin{lemma} \label{lemma_approximation_on_next_face}
		Let $P \in \Pn$ and $F \subset P$ a facet. For $m \geq 2$, let $f\colon \relint F \rightarrow \R$ be a $\tfrac{1}{m}$-convex integrable function. For every $\varepsilon>0$, there exists a $\tfrac{1}{m-1}$-convex integrable function $g\colon\relint P \rightarrow \R$ with $\norm{[f]_{\M(F)}-[g]_{\M(P)}}_{\Wt}<\varepsilon$ and $\norm{[f]_{\M(F)}-[g]_{\M(P)}}_{\TV}< 2\cdot\norm{[f]_{\M(F)}}_{\TV}+\varepsilon$.
	\end{lemma}
	
	\begin{proof}
		Let $\varepsilon>0$. By Lemma \ref{lemma_pac_approximation}, there exists a polyhedral convex function $h$ on $\relint F$ such that $\norm{[f-h^m]_{\M(F)}}_{\TV}<\varepsilon$. We denote the orthogonal projection onto $\aff F$ by $\pi_{\aff F}$ and the canonical extension of $h$ by $\overline{h}$. For $\delta>0$, we define $h_\delta\colon \relint P \rightarrow \R$ by
		\begin{equation}
			h_\delta(x) \coloneqq \max\left\{\overline{h}(\pi_{\aff F}(x))-\tfrac{1}{\delta}\cdot\dist(x,\aff F),0\right\}.
		\end{equation}
		Since $x \mapsto h(\pi_{\aff F}(x))$, $x\mapsto -\dist(x,\aff F)$ and $x \mapsto 0$ are convex functions on $\relint P$, $h_\delta$ is convex as well. Therefore, the function
		\begin{equation}
			g_\delta(x) \coloneqq \frac{m}{\delta} \cdot h_\delta(x)^{m-1}
		\end{equation}
		is $\tfrac{1}{m-1}$-convex. 
		
		Let $\mu_\delta  \coloneqq [g_\delta]_{\M(P)}$ and $\nu \coloneqq [h^m]_{\M(F)}$. It remains to show that $\delta>0$ can be chosen such that $\overline{W}(\mu_\delta,\nu)$ becomes arbitrarily small. We first choose a polytope $G \subset \relint F$ such that $\nu(F \setminus G)< \varepsilon$. Let $u_F \in \sphere$ be the outer normal vector of $F$. Since $G \subset \relint F$, there exists a $\delta_1>0$ such that $\conv\{G,G-\delta_1\cdot u_F\} \subset P$. We decompose $P$ into three parts:
		\begin{equation}
			P_1 \coloneqq P \cap \left(\pi_{\aff F}^{-1}(G)\right), \quad P_2 \coloneqq P \cap \left(\pi_{\aff F}^{-1}(F \setminus G)\right), \quad P_3 \coloneqq P \setminus\left(\pi_{\aff F}^{-1}(F)\right).
		\end{equation}
		We first consider $P_1$. If $\delta \leq \delta_1$, then Fubini's theorem implies that a density of $(\pi_{\aff F})_*(\mu_{\delta}|_{P_1})$ is given by
		\begin{equation}
			G\rightarrow [0, \infty), \quad y\mapsto \int_{P \cap [y,y-\delta \cdot h(y) \cdot u_F]} g_\delta(x) \dint x=\frac{m}{\delta} \cdot\int_0^{\delta h(y)} (h(y)-\tfrac{z}{\delta})^{m-1}\dint z=h(y)^m,
		\end{equation}
		which shows that $(\pi_{\aff F})_*(\mu_\delta|_{P_1})=\nu|_G$. Using the transference plan $\tau$ induced by the function $\pi_{\aff F}$, i.e., setting $\tau \coloneqq (\id, \pi_{\aff F})_*(\mu_\delta|_{P_1})$, we obtain that the Wasserstein distance of $\mu_\delta|_{P_1}$ and $\nu|_G$ is bounded by
		\begin{equation}
			W(\mu_\delta|_{P_1},\nu|_G) \leq \int_{G \times P_1} \norm{x-\pi_{\aff F}(x)}_2 \dint \tau(x,y) \leq \nu(G) \cdot \delta \cdot \max_{y \in G} h(y).
		\end{equation}
		In principle, the densities of $(\pi_{\aff F})_*(\mu_\delta|_{P_2})$ and $(\pi_{\aff F})_*(\mu_\delta|_{P_3})$ can be computed similarly, but for $y \in \pi_{\aff F}(P)\setminus G$ the segment $[y,y-\delta \cdot \overline{h}(y) \cdot u_F]$ is not necessarily contained in $P$. Therefore, we have the upper bounds $\mu_\delta(P_2) \leq \nu(F \setminus G)$ and
		\begin{equation}
			\mu_\delta(P_3) \leq \int_{\pi_{\aff F}(P_3)} \One_{\{P \cap [y,y-\delta \cdot \overline{h}(y) \cdot u_F] \neq \emptyset\}} \cdot \overline{h}(y)^m \dint y.
		\end{equation}
		Since $\overline{h}$ is bounded on $\pi_{\aff F}(P)$, there is a $\delta_2>0$ such that $\delta \leq \delta_2$ implies $\mu_\delta(P_3)<\varepsilon$.
		
		We now combine the bounds for $P_1$, $P_2$ and $P_3$. If 
		\begin{equation}
			\delta \leq \min\left\{\delta_1,\delta_2, \varepsilon \cdot \left[\nu(G)\cdot\textstyle\max_{y \in G} h(y)\right]^{-1}\right\},
		\end{equation}
		then by the estimates above we have
		\begin{equation}
			\overline{W}([g_\delta]_{\M(P)},[h^m]_{\M(F)})\leq \mu_\delta(P_2 \cup P_3)+\nu(F \setminus G)+W(\mu_\delta|_{P_1},\nu|_G)< 4\varepsilon.
		\end{equation}
		By the triangle inequality, it follows that
		\begin{equation}
			\norm{[f]_{\M(F)}-[g_\delta]_{\M(P)}}_{\Wt} \leq \norm{[f-h^m]_{\M(F)}}_{\TV}+\norm{[h^m]_{\M(F)}-[g_\delta]_{\M(P)}}_{\Wt} < 5 \varepsilon.
		\end{equation}
		With regard to the total variation norm, we observe that
		\begin{equation}
			\norm{[h^m]_{\M(F)}-[g_\delta]_{\M(P)}}_{\TV} \leq 2\cdot\nu(F)+\mu_\delta(P_3) <  2\cdot\norm{[f]_{\M(F)}}_{\TV}+\varepsilon,
		\end{equation}
		which implies
		\begin{align}
			\norm{[f]_{\M(F)}-[g_\delta]_{\M(P)}}_{\TV} &\leq \norm{[f-h^m]_{\M(F)}}_{\TV}+\norm{[h^m]_{\M(F)}-[g_\delta]_{\M(P)}}_{\TV} \\&< 2\cdot\left(\norm{[f]_{\M(F)}}_{\TV}+\varepsilon\right). \qedhere
		\end{align}
	\end{proof}

	We are now ready to prove the main result of this section.
	
	\begin{proof}[Proof of Proposition \ref{prop_main_result_first_inclusion}]
		Let $\mu \in \sum_{F \in \Phi(P)} -\M^{\cvx}_{(n-\dim F)^{-1}}(F)$. We consider a face $F \in \Phi(P)$ and define $f_F \colon \relint F \allowbreak\rightarrow \R$ as the unique $(n-\dim F)^{-1}$-convex function such that $\mu|_{\relint F} = -[f_F]_{\M(\relint F)}$. Let $F_{n-1}$ be a facet of $P$ that contains $F$ and let $F \subset F_{\dim F+1} \subset \dots \subset F_{n-1}$ be a chain of faces that contains an element of every intermediate dimension $d \in \{\dim F, \dots, n-1\}$. Let $i\in \N$. An iterated application of Lemma \ref{lemma_approximation_on_next_face}, using the triangle inequality for $\norm{\cdot}_{\Wt}$, and, finally, an application of Lemma \ref{lemma_pac_approximation} yields a polyhedral convex function $g_F \colon \relint F_{n-1} \rightarrow \R$ with $\norm{[f_F]_{\M(F)}-[g_F]_{\M(F_{n-1})}}_{\Wt}< \frac{1}{i }$ and $\norm{[f_F]_{\M(F)}-[g_F]_{\M(F_{n-1})}}_{\TV}< 3\norm{[f_F]_{\M(F)}}_{\TV}$. We repeat this construction for every $F \in \Phi(P)$ and set
		\begin{equation}
			\mu_i \coloneqq \sum_{F \in \Phi(P)} -[g_F]_{\M(\supp g_F)} \in \W^{\triangleleft}(P).
		\end{equation} 
		Performing this procedure for every $i \in \N$, we obtain a sequence $(\mu_i)_{i \in \N}$ in $\W^\triangleleft(P)$ which satisfies
		\begin{equation}
			\norm{\mu_i-\mu}_{\TV} \leq \sum_{F \in \Phi(P)} 3\norm{[f_F]_{\M(F)}}_{\TV} = 3\norm{\mu}_{\TV}
		\end{equation}
		as well as
		\begin{equation}
			\norm{\mu_i-\mu}_{\Wt} \leq \sum_{F \in \Phi(P)} \frac{1}{i} = \frac{\# \Phi(P)}{i} \goestoinfty{i} 0.
		\end{equation}
		By Proposition \ref{prop_wasserstein_metrizes_weak_convergence} and Proposition \ref{prop_set_of_perturbations_weakly_closed}, it follows that $\mu \in \W(P)$.
	\end{proof}

	\section{Weak limits of contrast sequences} \label{sect_proof_of_main_result}
	
	In this section, we complete the proof of Theorem \ref{thm_main_result}. It remains to show the second inclusion, i.e., that there are no other small perturbations than the ones that we constructed in the preceding sections.%
	
	\begin{definition} \label{def_P_coconvex_sequence}
		Let $P \subset \R^n$ be a full-dimensional polytope. A sequence of signed measures $(\mu_i)_{i \in \N}$ in $\M(\R^n)$ is called a \textit{$P$-contrast sequence} if for every $i\in \N$ there exist $\beta_i >0$ and $K_i \in \Knn \cup \{\emptyset\}$ such that
		\begin{equation}
			\mu_i = \beta_i[\One_{K_i}-\One_P]_{\M}.
		\end{equation}
		Clearly, the sequence $[(\beta_i,K_i)]_{i \in \N}$ is uniquely determined by $(\mu_i)_{i \in \N}$; we call it the \textit{component sequence} of $(\mu_i)_{i \in \N}$.
		A $P$-contrast sequence $(\mu_i)_{i \in \N}$ is called \textit{admissible} with \textit{limit} $\mu \in \M(\R^n)$ if
		\begin{enumerate}[label=(\roman*)]
			\item $(\mu_i)_{i \in \N}$ converges weakly to $\mu$,
			\item $(\beta_i)_{i \in \N}$ diverges to infinity, and
			\item \label{def_P_coconvex_sequence_property_3} there exists a compact set $C \subset \R^n$ such that $\mu_i$ is supported on $C$ for all $i\in \N$.
		\end{enumerate}
		We note in passing that condition \ref{def_P_coconvex_sequence_property_3} is logically redundant and merely included for convenience. 	
	\end{definition}
	
	Evidently, Definition \ref{def_P_coconvex_sequence} is closely related to weak differentiability: If $(K_t)_{t \in [0,1]}$ is weakly differentiable, then $[(i,K_{1/i})]_{i \in \N}$ is the component sequence of a $P$-contrast sequence. Translated into the terminology of Definition \ref{def_P_coconvex_sequence}, the problem we face is to characterize limits of $P$-contrast sequences. %
	As we will discuss below, parts of our proof are inspired by the ele\-gant proof of the reverse Brunn-Minkowski inequality for coconvex bodies by Fillastre \cite{fil17}.
	
	We begin with some basic properties of contrast sequences.
	
	\begin{lemma} \label{lemma_basic_properties_of_P_coconvex_sequences}
		Let $P \in \Pnn$ and $(\mu_i)_{i \in \N}$ an admissible $P$-contrast sequence with limit $\mu$.
		\begin{enumerate}[label=(\roman*)]
			\item \label{lemma_basic_properties_of_P_coconvex_sequences_1} The sequence $(\mu_i)_{i \in \N}$ is bounded.
			\item \label{lemma_basic_properties_of_P_coconvex_sequences_2} The sequence $(K_i)_{i \in \N}$ converges to $P$ with respect to $\dH$.
			\item \label{lemma_basic_properties_of_P_coconvex_sequences_3} The limit $\mu$ is supported on $\bd P$.
		\end{enumerate}
	\end{lemma}

	\begin{proof}
		As discussed in Remark \ref{rem_wasserstein_metrizes_weak_convergence}, \ref{lemma_basic_properties_of_P_coconvex_sequences_1} follows from the uniform boundedness principle, using that $\mu_i$ is supported on a fixed compact set $C$ for all $i\in\N$. Since $\beta_i$ diverges to infinity, the boundedness of $(\mu_i)_{i \in \N}$ implies that $\dS(K_i, P)$ converges to zero. Using the fact that $\dH$ and $\dS$ induce the same topology on $\Kn$, this shows \ref{lemma_basic_properties_of_P_coconvex_sequences_2}. Finally, \ref{lemma_basic_properties_of_P_coconvex_sequences_3} follows from \ref{lemma_basic_properties_of_P_coconvex_sequences_2} by a similar argument as in the proof of Proposition \ref{prop_weak_derivatives_are_signed_measures}.
	\end{proof}

	The following lemma is a partial version of the continuous mapping theorem for signed measures.
	\begin{lemma} \label{lemma_continuous_mapping}
		Let $(\mu_i)_{i \in \N}$ be a weakly convergent sequence in $\M(\R^n)$ and $\pi\colon \R^n \rightarrow \R^n$ a continuous map. We have
		\begin{equation}
			\pi_* \wlim_{i \rightarrow \infty} \mu_i = \wlim_{i \rightarrow \infty} \pi_*\mu_i.
		\end{equation}
		In particular, if $\pi_* \wlim_{i \rightarrow \infty} \mu_i$ is supported on $\{x \in \R^n \mid \pi(x)=x\}$, then
		\begin{equation}
			\wlim_{i \rightarrow \infty} \mu_i = \wlim_{i \rightarrow \infty} \pi_*\mu_i.
		\end{equation}
	\end{lemma}
	\begin{proof}
		Let $f \in C(\R^n)$. Because $\mu_i \warrow \mu$, we have
		\begin{equation}
			\int f \dint \pi_*\mu_i = \int f \circ \pi \dint \mu_i \goestoinfty{i} \int f \circ \pi \dint \mu = \int f \dint \pi_*\mu. \qedhere
		\end{equation}
	\end{proof}
	
	Our basic strategy for characterizing limits of $P$-contrast sequences is to consider the signed measures locally, i.e., on appropriately chosen compact subsets of $\R^n$.
	
	\begin{lemma} \label{lemma_local_restriction_1}
		Let $(\mu_i)_{i \in \N}$ be a sequence in $\M(\R^n)$ that converges weakly to $\mu \in \M(\R^n)$. Let $Q \subset \R^n$ be compact. By Remark \ref{rem_bounded_subseq}, there is a subsequence $(\mu_{k(i)})_{i \in \N}$ such that $(\mu_{k(i)}|_Q)_{i \in \N}$ converges weakly to a signed measure $\nu \in \M(Q)$. We have $\nu|_{\interior Q}= \mu|_{\interior Q}$.
	\end{lemma}

	\begin{proof}
		Let $f \in C_c(\interior Q)$. We note that the trivial extension $\overline{f}\colon \R^n \rightarrow \R$ is an element of $C(\R^n)$. %
		Therefore, we have
		\begin{equation}
			\int f \dint\nu|_{\interior Q}= \lim_{i \rightarrow \infty} \int_{Q} \overline{f} \dint\mu_{k(i)}= \int \overline{f} \dint\mu = \int \overline{f} \dint\mu|_{\interior Q}.
		\end{equation}
		Since $C_c(\interior Q)$ separates points in $C_0(\interior Q)^*$, it follows that $\nu|_{\interior Q}= \mu|_{\interior Q}$.
	\end{proof}

	After restricting the $P$-contrast sequence under consideration to an appropriately chosen compact subset of $\R^n$, we want to show that the weak limit of the restricted sequence is absolutely continuous with a concave density. For this, we need the following lemma.
	
	\begin{lemma} \label{lemma_uniform_concave_limit}
		Let $P,Q \in \Pn$ with $Q \subset \relint P$ and $m \coloneqq \dim P = \dim Q$. Let $(f_i)_{i \in \N}$ be a sequence of concave functions $P \rightarrow \R$ that is bounded in $L^1(P)$, i.e., we have $C \coloneqq\sup_{i \in \N} \int_P |f_i(x)| \dint x < \infty$. 
		\begin{enumerate}[label=(\roman*)]
			\item \label{lemma_uniform_concave_limit_boundedness} There exist constants $C_\ell, C_u \in \R$ with $f_i(x) \in [C_\ell,C_u]$ for all $x \in Q$ and $i \in \N$.
			\item \label{lemma_uniform_concave_limit_sequence} There exists a subsequence $(f_{k(i)})_{i \in \N}$ that converges uniformly on $Q$ to a concave function $f \colon Q \rightarrow \R$.
		\end{enumerate}
	\end{lemma}
	
	\begin{proof}
		We start with (i). Setting $T \coloneqq \frac{1}{2}(P+Q)$, we have $Q \subset \relint T \subset T \subset \relint P$. We first show that $\{f_i|_{T}\}_{i \in \N}$ is uniformly bounded from below.
		For an arbitrary $z \in P$, let $S \coloneqq \sphere \cap (\aff P -z)$. %
		We define a function
		\begin{equation}
			h\colon T \times S \rightarrow \R, \quad (y,u) \mapsto \volt{m}\left(\{x \in P \mid \scpr{u,x} \leq \scpr{u,y}\}\right).
		\end{equation}
		Because $T \subset \relint P$, the set $\{x \in P \mid \scpr{u,x} \leq \scpr{u,y}\}$ is an $m$-dimensional polytope for all $(y,u) \in T \times S$. Therefore, we have $h(y,u)>0$ for all $(y,u) \in T \times S$. Since $T \times S$ is compact and $h$ is continuous, it follows that
		\begin{equation}
			M \coloneqq \min_{(y,u) \in T \times S}  h(y,u)>0.
		\end{equation}
		Let $i \in \N$ and $x \in T$. If $f_i(x)\geq 0$, there is nothing to show, so we assume that $f_i(x)< 0$. Because $f_i$ is concave, the superlevel sets of $f_i$ are convex, and hence there exists a $u \in S$ such that $f_i(y)\leq f_i(x)$ for all $y \in \{x \in P \mid \scpr{u,x} \leq \scpr{u,y}\}$. This leads to the estimate
		\begin{equation}
			-f_i(x) \cdot M \leq \int_P |f_i(x)| \dint x \leq C,
		\end{equation}
		which implies $f_i(x) \geq -\frac{C}{M} \eqqcolon C_\ell$ for all $x \in T \supset Q$ and $i \in \N$.
		
		To complete the proof of (i), we now show that $\{f_i|_{Q}\}_{i \in \N}$ is uniformly bounded from above. Let $i\in \N$ and $x \in Q$. We distinguish two cases:
		\begin{enumerate}[label=\alph*)]
			\item If $f_i(y) \geq 0$ for all $y \in Q$, then by the concavity of $f_i$, the $(m+1)$-dimensional pyramid
			\begin{equation}
				K\coloneqq \conv\big((Q \times \{0\}) \cup \{(x,f_i(x))\}\big) \subset Q \times \R
			\end{equation}
			is contained in the subgraph of $f_i$. Therefore, we have
			\begin{equation}
				\tfrac{1}{m+1} \cdot f_i(x) \cdot \volt{m}(Q) =\volt{n}(K)\leq \int_P |f_i(x)| \dint x\leq C,
			\end{equation}
			which yields a uniform upper bound on $f_i(x)$.
			\item We assume there exists a $v \in Q$ with $f_i(v)<0$. We define a map 
			\begin{equation}
				b \colon \{(x,y)  \in Q \times Q \mid x \neq y\} \rightarrow \relbd T
			\end{equation} 
			by stipulating that $b(x,y)$ is the unique intersection point of $\relbd T$ with the ray that emanates from $x$ and goes through $y$. Because $Q \subset\relint T$, we have
			\begin{equation}
				\beta \coloneqq \min_{x \in Q} \dist(x, \relbd T) > 0.
			\end{equation}
			If $f_i(x) \leq f_i(v)$, then $f_i(v) < 0$ yields the desired upper bound, so we assume that $f_i(x)> f_i(v)$. Then, by the concavity of $f_i$, we have
			\begin{equation}
				\frac{f_i(x)-f_i(v)}{\norm{x-v}_2} \leq \frac{f_i(v)-f_i(b(x,v))}{\norm{v-b(x,v)}_2} \leq \frac{-C_\ell}{\beta}
			\end{equation}
			and hence
			\begin{equation}
				f_i(x) \leq \diam(Q)\cdot \frac{-C_\ell}{\beta}.
			\end{equation}
		\end{enumerate}
		Combining both cases, we obtain that $C_u \coloneqq \max\left\{\frac{C(m+1)}{\vol_{m}(Q)},\frac{-C_\ell\diam(Q)}{\beta}\right\}$ is a uniform upper bound on $f_i|_Q$.
		
		To show (ii), we construct the subsequence $(f_{k(i)})_{i \in \N}$. Let $X = \{x_j\}_{j \in \N} \subset P$ be a countable set with $P \subset\cl X$. Because $\{f_i(x_j)\}_{j \in \N} \subset [C_\ell,C_u]$ for every $i \in \N$, the Bolzano-Weierstraß theorem implies that there exists a nested sequence of index sets $k_1(\N) \supset k_2(\N) \supset \dots$ such that $(f(x_j)_{k_j(i)})_{i \in \N}$ converges for every $j \in \N$. Setting $k(i) \coloneqq k_i(i)$, we obtain a subsequence $(f_{k(i)})_{i \in \N}$ that converges pointwise on $X$. Now Theorem \ref{thm_rockafellar} implies that $f$ converges uniformly on $Q$ to a concave function $f \colon Q \rightarrow \R$.
	\end{proof}

	As the first step towards a characterization of limits of $P$-contrast sequences, we determine their behavior on the facets of $P$.
	
	\begin{lemma} \label{lemma_facet_concave_density}
		Let $P \subset \R^n$ be a full-dimensional polytope and  $(\mu_i)_{i \in \N}$ an admissible $P$-contrast sequence with limit $\mu$. Let $F \subset P$ be a facet and $G \subset \relint F$ a polytope with $\dim G= \dim F$. Then $\mu|_{\relint G}$ is absolutely continuous with respect to $\vol_{n-1} \in \M(F)$ with a concave density $f\colon \relint G \rightarrow \R$.
	\end{lemma}
	
	\begin{proof}
		Let $u_F \in \sphere$ be the unit outer normal vector of $F$. Setting $\tilde{G} \coloneqq \frac{1}{2}(F+G)$, we have $G \subset \relint \tilde{G} \subset \tilde{G} \subset \relint F$.
		Let $\alpha>0$ be such that
		\begin{equation}
			\tilde{G}-\alpha u_F \subset \interior P.
		\end{equation}
		Let $\big[(\beta_i,K_i)\big]_{i \in \N}$ be the component sequence of $(\mu_i)_{i \in \N}$. By Lemma \ref{lemma_basic_properties_of_P_coconvex_sequences}\ref{lemma_basic_properties_of_P_coconvex_sequences_2}, we can assume without loss of generality that $\tilde{G}-\alpha u_F \subset K_i$ for all $i \in \N$. %
		We define two polytopes
		\begin{equation}
			Q \coloneqq G+[-\alpha u_F,\alpha u_F] \quad \text{and} \quad \tilde{Q} \coloneqq \tilde{G}+[-\alpha u_F,\alpha u_F].
		\end{equation}
		By Remark \ref{rem_bounded_subseq}, there exists a subsequence $(\mu_{k(i)})_{i \in \N}$ with $(\mu_{k(i)}|_{Q})_{i \in \N} \warrow \nu$ for some $\nu \in \M(Q)$. Arguing as in the proof of Proposition \ref{prop_weak_derivatives_are_signed_measures}, we obtain that $\nu$ is supported on $G$. Using Lemma \ref{lemma_continuous_mapping} and Lemma \ref{lemma_local_restriction_1}, we have
		\begin{equation}
			\mu|_{\interior Q} = \nu|_{\interior Q} = ((\pi_{\aff  F})_*\nu)|_{\interior Q} = \left.\left(\wlim_{i \rightarrow \infty} (\pi_{\aff  F})_*(\mu_{k(i)}|_{Q})\right)\right|_{\interior Q}.
		\end{equation}
		For every $i \in \N$, we define $g_i\colon \tilde{G} \rightarrow \R$ by
		\begin{equation}
			g_i(x)\coloneqq \beta_i\cdot \max\{t \in \R \mid x+t\cdot u_F \in K_i\}.
		\end{equation}
		Because we assumed that $\tilde{G}- \alpha u_F \subset K_i$ for all $i \in \N$, every $g_i$ is well-defined and bounded from below by $\alpha \beta_i$. Moreover, the convexity of $K_i$ implies that $g_i$ is concave. It follows from Fubini's theorem that $g_i|_{G}$ is a density of $(\pi_{\aff F})_*(\mu_i|_{Q})$. By Lemma \ref{lemma_basic_properties_of_P_coconvex_sequences}\ref{lemma_basic_properties_of_P_coconvex_sequences_1}, there exists a constant $C>0$ such that
		\begin{equation}
			\int_{\tilde{G}} |g_i(x)| \dint x \leq \norm{\mu_i}_{\TV} \leq C \quad \text{for all }i \in \N.
		\end{equation}
		Hence, Lemma \ref{lemma_uniform_concave_limit}\ref{lemma_uniform_concave_limit_sequence} implies that we can find a subsequence $\{g_{k'(i)}|_{G}\}_{i \geq N}$ of $\{g_{k(i)}|_{G}\}_{i \geq N}$ that converges uniformly on $G$ to a concave function $g^*\colon G \rightarrow \R$. Since uniform convergence implies $L^1$-convergence, it follows that $[g_{k'(i)}|_{G}]_{\M(G)} \warrow [g^*]_{\M(G)}$. By the uniqueness of weak limits, it follows that $[g^*]_{\M(G)}=[g^*]_{\M(G)}|_{\interior Q} = \mu|_{\interior Q}$.
	\end{proof}

	We now extend the assertion of Lemma \ref{lemma_facet_concave_density} to whole facets of $P$.
	
	\begin{proposition} \label{prop_facet_concave_density}
		Let $P \subset \R^n$ be a full-dimensional polytope, $F \subset P$ a facet and  $(\mu_i)_{i \in \N}$ an admissible $P$-contrast sequence with limit $\mu$. Then $\mu|_{\relint F}$ is absolutely continuous with respect to $\vol_{n-1} \in \M(F)$ with a concave density $f\colon \relint F \rightarrow \R$.
	\end{proposition}
	
	\begin{proof}
		Let $F_1 \subset F_2 \subset \dots \subset \relint F$ be a sequence of polytopes with $\bigcup_{i \in \N}F_i=\relint F$. %
		For $i \in \N$, let $f_i \colon \relint F_i \rightarrow \R$ be the concave density of $\mu|_{\relint F_i}$ whose existence was shown in Lemma \ref{lemma_facet_concave_density}. Let $j,k \in \N$ with $j<k$. Then $f_k$ and $f_j$ agree on $\relint F_j$, since both densities are continuous and describe the same signed measure on $\relint F_j$. Therefore, the function $f \colon \relint F\rightarrow \R$ given by $f(x) \coloneqq f_{\min \{i \in \N \mid x\in \relint F_i\}}(x)$ is concave. Let $g \in C_c(\relint F)$. Since $\supp g$ is compact, there is an $i\in \N$ with $\supp g \subset \relint F_i$. By construction, $f|_{\relint F_i}$ is a density of $\mu|_{\relint F_i}$, so we have $\int g \dint \mu = \int gf \dint\volt{n-1}(F)$. Since $C_c(\interior Q)$ separates points in $C_0(\interior Q)^*$, it follows that $\mu|_{\relint F}=[f]_{\M(\relint F)}$.
	\end{proof}

	In the second step, we determine the behavior of limits of $P$-contrast sequences on lower-dimensional faces of $P$. The following characterization of weak convergence is part of the well-known Portemanteau theorem \cite[Thm.~13.16]{kle20}.
	
	\begin{theorem}[Portemanteau] \label{thm_portemanteau}
		Let $X$ be a metric space and $\mu, \mu_1, \mu_2, \ldots \in \M^+(X)$ with $\sup_{i \in \N} \mu(X)< \infty$. Then
		$\mu_i \warrow \mu$ if and only if
		\begin{equation}
			\limsup_{i \rightarrow \infty} \mu_n(X)\leq \mu(X) \quad \text{and} \quad \liminf_{i \rightarrow \infty} \mu_i(G) \geq \mu(G) \quad\text{for all open sets } G \subset X.
		\end{equation}
	\end{theorem}

	We use Theorem \ref{thm_portemanteau} to show that the restriction of a limit of a $P$-contrast sequence to a face $F \subset P$ with $\dim F < n-1$ is a negative measure.
	
	\begin{lemma} \label{lemma_m_face_negative_measure}
		Let $P \subset \R^n$ be a full-dimensional polytope and  $(\mu_i)_{i \in \N}$ an admissible $P$-contrast sequence with limit $\mu$. Let
		\begin{equation}
			S\coloneqq \bd P \setminus \left(\bigcup_{\substack{F \in \Phi_{n-1}(P)}}\relint F\right),
		\end{equation}
		i.e., $S$ is the union of all proper faces of $F$ that are not facets.
		Then there exists a restricted subsequence $(\mu_{k(i)}|_P)_{i \in \N}$ in $\M^-(P)\coloneqq -\M^+(P)$ that converges weakly to a signed measure $\nu$ with $\nu|_S=\mu|_S$. In particular, we have $\mu(A) \leq 0$ for all Borel sets $A \subset S$.
	\end{lemma}
	\begin{proof}
		Let $P= \{x \in \R^n \mid \scpr{a_i,x} \leq b_i, i \in [m]\}$ for $a_i \in \R^n$ and $b_i \in \R$. For $j \in [m]$, let
		\begin{equation}
			Q_j \coloneqq \{x \in \R^n \mid \scpr{a_j,x} \geq b_j, \, \scpr{a_i,x} \leq b_i \text{ for } 1 \leq i < j \}.
		\end{equation} 
		We decompose $\mu_i$ into the sum
		\begin{equation}
			\mu_i = \mu_i|_{P}+\sum_{j \in [m]} \mu_i|_{Q_j} \quad \text{for all } i\in \N.
		\end{equation}
		By Definition \ref{def_P_coconvex_sequence}\ref{def_P_coconvex_sequence_property_3}, there exists a compact set $C \subset \R^n$ such that $\mu_i$ is supported on $C$ for all $i \in \N$. Hence, Remark \ref{rem_bounded_subseq} implies that there exists a subsequence $(\mu_{k(i)})_{i \in \N}$ such that every part of the decomposition converges weakly, i.e., there exist signed measures $\nu \in \M^-(\R^n)$ and $\nu_1,\dots,\nu_m \in \M^+(\R^n)$ with $\mu_i|_P \warrow \nu$ and $\mu_i|_{Q_j} \warrow \nu_j$ for $j \in [m]$. Since $\mu=\nu + \sum_{j \in [m]} \nu_j$, the claim follows if we show that $\nu_j(S) = 0$ for $j \in [m]$.
		
		Let $\big[(\beta_i,K_i)\big]_{i \in \N}$ be the component sequence of $(\mu_i)_{i \in \N}$. By Lemma \ref{lemma_continuous_mapping}, we know that $(\pi_{\aff F})_*(\mu_{k(i)}|_{Q_1}) \warrow \nu_1$. Let $F$ be the facet $\{x \in P \mid \scpr{a_1,x}=b_1\}$, $u_F \in \sphere$ its outer normal vector and $z \in \relint F$. We define
		\begin{equation}
			F_{\frac{1}{2}} \coloneqq \tfrac{1}{2}(F-z)+z, \quad F_{\frac{2}{3}} \coloneqq \tfrac{2}{3}(F-z)+z \quad \text{and} \quad F_{2} \coloneqq 2(F-z)+z.
		\end{equation}
		Let $\alpha>0$ be such that $F_{\frac{2}{3}}-\alpha u_F \subset \interior P$. By Lemma \ref{lemma_basic_properties_of_P_coconvex_sequences}\ref{lemma_basic_properties_of_P_coconvex_sequences_2}, there exists an $N \in \N$ such that
		\begin{equation}\label{eq_lemma_m_face_negative_measure_choice_of_N}
			F_{\frac{2}{3}}-\alpha u_F \subset K_i \quad \text{and} \quad \pi_{\aff F}(K_i \cap Q_1) \subset F_2 \quad \text{for all } i \geq N.
		\end{equation}
		For $i \in \N$, we define $g_i\colon F_2 \rightarrow \R \cup \{- \infty\}$ by
		\begin{equation}
			g_i(x) \coloneqq \sup\{t \in \R \mid x+t\cdot u_F \in K_i\}.
		\end{equation}
		Because $K_i$ is convex, $g_i$ is a concave function $F_2 \rightarrow \R \cup \{- \infty\}$ for all $i \in \N$. Moreover, by \eqref{eq_lemma_m_face_negative_measure_choice_of_N}, the restriction $g_i|_{F_{\frac{2}{3}}}$ is real-valued for all $i \geq N$. Using Lemma \ref{lemma_uniform_concave_limit}\ref{lemma_uniform_concave_limit_boundedness}, we obtain that there exist constants $C_\ell, C_u \in \R$ such that %
		\begin{equation}
			\beta_i \cdot g_i \left(F_{\frac{1}{2}}\right) \subset \left[C_\ell,C_u\right] \quad \text{for all } i \geq N.
		\end{equation}
		Let $x \in F_2$. We set $y \coloneqq \frac{1}{4}x+\frac{3}{4}z=\frac{1}{4}(x-z)+z \in F_{\frac{1}{2}}$. Because $g_i \colon F_2 \rightarrow \R \cup \{- \infty\}$ is concave, we have $g_i(y) \geq \frac{1}{4}g_i(x)+\frac{3}{4}g_i(z)$ and hence
		\begin{equation}
			\beta_i \cdot g_i(x) \leq \beta_i \cdot \left(4 g_i(y)-3 g_i(z)\right) \leq 4C_u -3 C_\ell \eqqcolon \tilde{C} \quad \text{for all } i \geq N.
		\end{equation}	
		Since $\max\{\beta_i\cdot g_i,0\} \colon F_2 \rightarrow [0,\infty)$ is a density of $(\pi_{\aff F})_*(\mu_{i}|_{Q_1})$, this implies that
		\begin{equation}
			[(\pi_{\aff F})_*(\mu_{i}|_{Q_1})](A) \leq \tilde{C} \cdot\volt{n-1}(A) \quad \text{for all Borel sets } A \subset \aff F \text{ and }i \geq N.
		\end{equation}
		Let $\varepsilon>0$ and let $U \subset \aff F$ be (relatively) open with $\relbd F \subset U$ and $\volt{n-1}(U)<\varepsilon$. By Theorem \ref{thm_portemanteau}, we have
		\begin{equation}
			\nu_1(\relbd F)\leq\nu_1(U) \leq \liminf_{i \rightarrow \infty} [(\pi_{\aff F})_*(\mu_{k(i)}|_{Q_1})](U) \leq \tilde{C} \cdot\volt{n-1}(U) < \tilde{C} \varepsilon.
		\end{equation}
		Since $\varepsilon>0$ was arbitrary, this implies $\nu_1(S)=\nu_1(\relbd F)=0$.
		
		To complete the proof, we observe that $(\mu_i|_{P}+\sum_{j=2}^m \mu_i|_{Q_j})_{i \in \N}$ is again a $P$-contrast sequence, and we can iterate the argument above to obtain that $\nu_j(S) = 0$ for $j \in [m]$.
	\end{proof}

	\begin{remark} \label{rem_wlog_K_subset_P}
		Let $P$, $(\mu_i)_{i \in \N}$, $\mu$ and $S$ be as in Lemma \ref{lemma_m_face_negative_measure}. Moreover, let $F \subset P$ be an $m$-face for $0 \leq m < n-1$. Since $F \subset S$, Lemma \ref{lemma_m_face_negative_measure} implies that there exists a restricted subsequence $(\mu_{k(i)}|_P)_{i \in \N}$ that converges weakly to a signed measure $\nu$ with $\nu|_F=\mu|_F$. Therefore, for the purpose of determining $\mu|_F$, we can replace $(\mu_i)_{i \in \N}$ by $(\mu_{k(i)}|_P)_{i \in \N}$.
	\end{remark}
	
	The localization argument for the lower-dimensional faces of $P$ is somewhat more involved than the corresponding argument for the facets. Instead of restricting the sequence $(\mu_i)_{i \in \N}$ to a single compact set $Q \subset \R^n$, we consider a whole sequence of compact sets.
	
	\begin{lemma} \label{lemma_local_restriction_2}
		Let $(\mu_i)_{i \in \N}$ be a sequence in $\M^+(\R^n)$ that converges weakly to $\mu \in \M^+(\R^n)$. Let $Q \subset \R^n$ and let $(Q_j)_{j \in \N}$ be a sequence of compact sets with $Q_{j+1} \subset \interior Q_j$ for $j \in \N$ and $\bigcap_{j \in \N} \interior Q_j =Q$. Then there is a subsequence $(\mu_{k(i)})_{i \in \N}$ such that
		\begin{equation}
			\mu|_Q = \wlim_{j \rightarrow \infty} \wlim_{i \rightarrow \infty} (\mu_{k(i)}|_{Q_j}).
		\end{equation}
	\end{lemma}

	\begin{proof}
		Applying Remark \ref{rem_bounded_subseq} iteratively, we obtain a subsequence $(\mu_{k(i)})_{i \in \N}$ such that $(\mu_{k(i)}|_{Q_j})_{i \in \N}$ converges weakly to a signed measure $\nu_j$ for all $j \in \N$. By Lemma \ref{lemma_local_restriction_1}, we have $\nu_j = \nu_j|_{\interior Q_j} + \nu_j|_{\bd Q_j} = \mu|_{\interior Q_j} + \nu_j|_{\bd Q_j}$ for all $j \in \N$. By the continuity of measures from above, we have 
		\begin{equation}
			\norm{(\mu|_{\interior Q_j})-(\mu|_{Q})}_{\TV}=\mu(\interior Q_j \setminus Q)\textstyle\goestoinfty{j}0
		\end{equation}
		and hence $\mu|_{\interior Q_j} \warrow \mu|_Q$. It remains to show that $\nu_j|_{\bd Q_j} \warrow 0$. For this, we observe that
		\begin{equation}
			\mu|_{\bd Q_j} = \left(\wlim_{i \rightarrow \infty}\left.\left( \mu_{k(i)}|_{\R^n \setminus Q_j}\right)\right)\right|_{\bd Q_j}+\nu_j|_{\bd Q_j}.
		\end{equation}
		Since $\bd Q_{j} \cap \bd Q_{j'} = \emptyset$ for $j \neq j'$, it follows that 
		\begin{equation}
			\sum_{j \in \N} \nu_j(\bd Q_j)\leq \sum_{j \in \N}  \mu(\bd Q_j) = \mu \left(\bigcup_{j\in\N}\bd Q_j\right)< \infty,
		\end{equation}
		which implies $\norm{(\nu_j|_{\bd Q_j})}_{\TV}=\nu_j(\bd Q_j) \goestoinfty{j} 0$.
	\end{proof}
	
	The following fact is a straightforward consequence of Fubini's theorem, we state it as a lemma for later reference.

	\begin{lemma} \label{lemma_fubini_non_orthogonal_projection}
		Let $H \subset \R^n$ be an affine subspace and $\pi \colon \R^n \rightarrow H$ a \textit{projection}, i.e., an affine map with $\pi|_H = \id_H$. Then there exists a constant $\gamma(\pi)\in (0,1]$ such that
		\begin{equation}
			\vol(A) = \gamma(\pi)\cdot\int_{\pi(A)} \volt{n-\dim H}(A \cap \pi^{-1}(\{x\})) \dint \volt{\dim H}(x)
		\end{equation}
		for all measurable sets $A \subset \R^n$.
	\end{lemma}

	We now combine the insights of the previous lemmas into a result that allows us to identify the densities of limits of $P$-contrast sequences on lower-dimensional faces of $P$.
	
	\begin{lemma} \label{lemma_limit_shrinking_sets}
		Let $P \subset \R^n$ be a full-dimensional polytope and  $(\mu_i)_{i \in \N}$ an admissible $P$-contrast sequence with limit $\mu$ and component sequence  $\big[(\beta_i,K_i)\big]_{i \in \N}$. Let $F\subset P$ be an $m$-face for $0 \leq m < n-1$ and $G \subset \relint F$ a polytope with $\dim G= \dim F$. Let $\pi\colon \R^n \rightarrow \aff F$ be a projection. Moreover, let $(Q_j)_{j \in \N}$ be a sequence of compact sets with $Q_{j+1} \subset \interior Q_j$ for $j \in \N$ and 
		\begin{equation}
			\bd P \cap \bigcap_{j \in \N}Q_j = G \subset \relint F.
		\end{equation}
		For $i,j \in \N$, let $g_i^j\colon \relint F \rightarrow (-\infty,0]$ be given by
		\begin{equation}
			g_i^j(x)\coloneqq-\beta_i\cdot \gamma(\pi) \cdot\volt{n-\dim F}(Q_j\cap \pi^{-1}(\{x\})\cap(P \setminus K_i)),
		\end{equation}
		where $\gamma(\pi)$ is the constant from Lemma \ref{lemma_fubini_non_orthogonal_projection}.
		Finally, we assume that for every $j \in \N$, the sequence $(g_i^j)_{i \in \N}$ converges uniformly on $\relint F$ to a function $g_\infty^j \colon \relint F \rightarrow \R$. Then the sequence of limits $(g^j_\infty)_{j \in \N}$ converges pointwise to a function $g^*\colon \relint F \rightarrow  (-\infty,0]$ and we have
		\begin{equation}
			\mu|_G=[g^*]_{\M(F)}.
		\end{equation}
	\end{lemma}
	
	\begin{proof}
		By Remark \ref{rem_wlog_K_subset_P}, we can assume without loss of generality that $K_i \subset P$ for all $i \in \N$. Lemma \ref{lemma_continuous_mapping} and Lemma \ref{lemma_local_restriction_2} imply that there exists a subsequence $(\mu_{k(i)})_{i \in \N}$ such that
		\begin{equation}
			\mu|_G=\pi_* (\mu|_G) = \wlim_{j \rightarrow \infty} \wlim_{i \rightarrow \infty} \pi_* (\mu_{k(i)}|_{Q_j}).
		\end{equation}
		It follows from Lemma \ref{lemma_fubini_non_orthogonal_projection} that $g_{k(i)}^j$ is a density of $\pi_* (\mu_{k(i)}|_{Q_j})$. By assumption, the sequence $([g_{k(i)}^j]_{\M(F)})_{i \in \N}$ converges to $[g_\infty^j]_{\M(F)}$ in the total variation norm and hence also weakly. It remains to show that $\wlim_{j \rightarrow \infty}[g_\infty^j]_{\M(F)} = [g^*]_{\M(F)}$.
		
		Since the sequence $(Q_j)_{j \in \N}$ is decreasing, the sequence $(g^j_\infty)_{j \in \N}$ satisfies $g^j_\infty\leq g^{j+1}_\infty \leq 0$ for all $j\in \N$ and hence converges pointwise to a function $g^*\colon \relint F \rightarrow (-\infty,0]$. By the monotone convergence theorem \cite[Thm.~4.3.2]{dud02}, it follows that
		\begin{equation}
			\norm{[g^j_\infty]_{\M(F)}-[g^*]_{\M(F)}}_{\TV}=\int_{\relint F} g^*-g^j_\infty\dint\textstyle{\vol_m} \goestoinfty{j} 0. \qedhere
		\end{equation}
	\end{proof}

	To show that limits of $P$-contrast sequences have concave densities on lower-dimension\-al faces of $P$, we make use of the following sharpened version of the Brunn-Minkowski inequality for convex sets with identical projections. This result is classical, see for example \cite[Ch.~50]{bf34}.
	\begin{proposition} \label{prop_sharpenend_brunn_minkowski}
		Let $K,L \in \Knn$ and $H \subset \R^n$ a hyperplane. If $\pi_H(K)=\pi_H(L)$, then for every $\lambda \in [0,1]$, we have
		\begin{equation}
			\vol(\lambda K+ (1-\lambda)L) \geq \lambda \vol(K)+ (1-\lambda) \vol(L).
		\end{equation}
	\end{proposition}
	
	We now prove an analogous result to Lemma \ref{lemma_facet_concave_density} for faces $F \subset P$ with $\dim F < n-1$.
	
	\begin{lemma} \label{lemma_m_face_concave_density}
		Let $P \subset \R^n$ be a full-dimensional polytope and  $(\mu_i)_{i \in \N}$ an admissible $P$-contrast sequence with limit $\mu$. Let $F\subset P$ be an $m$-face for $0 \leq m < n-1$. Then $\mu|_{\relint F}$ is absolutely continuous with respect to $\vol_m \in \M(F)$ with a concave density $f\colon \relint F \rightarrow (-\infty,0]$.
	\end{lemma}

	The conclusion in Lemma \ref{lemma_m_face_concave_density} is weaker than the claim that $-f^{(n-\dim F)^{-1}}$ is convex, which is asserted in Theorem \ref{thm_main_result}. We will strengthen this conclusion below.
	
	\begin{proof}[Proof of Lemma \ref{lemma_m_face_concave_density}]
		Let $\big[(\beta_i,K_i)\big]_{i \in \N}$ be the component sequence of $(\mu_i)_{i \in \N}$. Again, by Remark \ref{rem_wlog_K_subset_P}, we can assume without loss of generality that $K_i \subset P$ for all $i \in \N$. Let $G \subset \relint F$ be a polytope with $\dim G = \dim F$ and set $G_j \coloneqq \frac{j}{j+1} G + \frac{1}{j+1}F$ for $j\in \N$.
		
		Let $u_1,\dots,u_\ell \in \sphere$ be the outer normal vectors of the facets $F_1,\dots,F_\ell$ that contain $F$. For every $v \in \vertices(G_1)$, there is an $\alpha_v>0$ such that
		$v - t \cdot \sum_{i=1}^\ell u_i \in \interior P$ for all $t \in (0,\alpha_v]$.
		Setting $w \coloneqq \left(\min_{v \in \vertices(G_1)} \alpha_v \right)\cdot \sum_{i=1}^\ell u_i$, we have 
		\begin{equation}
			G_1 - [\tfrac{w}{2},w]=\conv\left(\vertices G_1 - [\tfrac{w}{2},w]\right) \subset \interior P.
		\end{equation}
		For an arbitrary $z \in \aff F$, let $B \coloneqq \linh(F-z \cup \{w\})^\bot \cap B_n$. Because $G_1 - [\tfrac{w}{2},w]$ is compact, there exists an $N \in \N$ such that
		\begin{equation}
			G_1 - [\tfrac{w}{2},w]+\tfrac{1}{j} B \subset \interior P \quad \text{for all } j \geq N.
		\end{equation}
		By Lemma \ref{lemma_basic_properties_of_P_coconvex_sequences}\ref{lemma_basic_properties_of_P_coconvex_sequences_2}, we can assume without loss of generality that
		\begin{equation} \label{eq_lemma_m_face_concave_density_1}
			G_1 - [\tfrac{w}{2},w]+\tfrac{1}{j} B \subset \interior K_i \quad \text{for all } i,j \geq N.
		\end{equation}
		We define two sequences of convex bodies via
		\begin{equation} \label{eq_def_qlj}
			Q_j \coloneqq G_j+\tfrac{1+j}{2j}[-w,w]+\tfrac{1}{j}B \quad \text{and} \quad \tilde{Q}_j \coloneqq G_1+\tfrac{1+j}{2j}[-w,w]+\tfrac{1}{j}B \quad \text{for } j \in \N.
		\end{equation}
		By construction, the sequence $(Q_j)_{j \in \N}$ satisfies the condition  $Q_{j+1} \subset \interior Q_j$ from Lemma \ref{lemma_limit_shrinking_sets}; the sequence $(\tilde{Q}_j)_{j \in \N}$ will allow us to apply Lemma \ref{lemma_uniform_concave_limit}.
		In the following, we will use Lemma \ref{lemma_limit_shrinking_sets} to show that $\mu|_{G}$ is absolutely continuous with respect to $\vol_m \in \M(G)$ with a concave density $\relint G \rightarrow (-\infty,0]$. The claim of the lemma can then be deduced by an argument as in the proof of Proposition \ref{prop_facet_concave_density}.
		
		For $i,j \in \N$, let $g_i^j\colon \relint G_1 \rightarrow \R$ be given by 
		\begin{equation}
			g_i^j(x)\coloneqq-\beta_i\cdot \volt{n-m}(\tilde{Q}_j\cap \pi_{\aff F}^{-1}(x)\cap(P \setminus K_i)).
		\end{equation}
		We will first show that $g_i^j$ is concave if $i$ and $j$ are sufficiently large. For $x \in \relint G_1$, let $L_j(x) \coloneqq \tilde{Q}_j\cap \pi_{\aff F}^{-1}(x)$. We note that by \eqref{eq_lemma_m_face_concave_density_1} and \eqref{eq_def_qlj} the definition of $L_j(x)$ implies that $L_j(x) \cap K_i \supset x-\frac{w}{2}+\frac{1}{j}B$ and hence
		\begin{equation} \label{eq_lemma_m_face_concave_density_2}
			\pi_{w^\bot}(L_j(x) \cap K_i) = x+\frac{1}{j}B \quad \text{for all } x \in G_1 \text{ and } i,j \geq N.
		\end{equation}
		We consider the $(n-m)$-dimensional convex body $L_j(x)\cap P$ for a given $x \in \relint G_1$. Let $P_F \supset P$ be the polyhedron that results from $P$ by deleting all facet-defining inequalities of $P$ that are not active at relative interior points of $F$, i.e., $P_F$ is bounded by the facet-defining hyperplanes $\aff F_1, \dots, \aff F_\ell$. Then there exists an $N' \geq N$ such that
		\begin{equation} \label{eq_lemma_m_face_concave_density_3}
			L_j(x)\cap P = L_j(x)\cap P_F \quad \text{for all } j \geq N'.
		\end{equation}
		In the following, we assume that $i, j \geq N'$. Because $G_1 \subset F$, we have
		\begin{equation}
			[L_j(x)\cap P_F]-x=[L_j(y)\cap P_F]-y \quad \text{for all } y \in \relint G_1
		\end{equation}
		and hence, using \eqref{eq_lemma_m_face_concave_density_3},
		\begin{equation} \label{eq_volumes_are_equal}
			\volt{n-m}[L_j(x)\cap P]=\volt{n-m}[L_j(y)\cap P] \quad \text{for all } y \in \relint G_1.
		\end{equation}
		Let $x,y \in \relint G_1$ and $\lambda \in [0,1]$. By the convexity of $K_i$, we have
		\begin{equation} \label{eq_lemma_m_face_concave_density_4}
			\lambda(L_j(x)\cap K_i)+(1-\lambda)(L_j(y)\cap K_i) \subset L_j(\lambda x+(1-\lambda)y)\cap K_i.
		\end{equation}
		Since $i \geq N$, we obtain from \eqref{eq_lemma_m_face_concave_density_2} that $\pi_{w^\bot}([L_j(x)\cap K_i]-x)=\pi_{w^\bot}([L_j(y)\cap K_i]-y)= \tfrac{1}{j}B$ and Proposition \ref{prop_sharpenend_brunn_minkowski} yields
		\begin{align}
			\volt{n-m}(\lambda(L_j(x)\cap K_i)&+(1-\lambda)(L_j(y)\cap K_i))\\&\geq \lambda \cdot \volt{n-m}(L_j(x)\cap K_i)+(1-\lambda)\cdot \volt{n-m}(L_j(y)\cap K_i).
		\end{align}
		In combination with \eqref{eq_lemma_m_face_concave_density_4} and \eqref{eq_volumes_are_equal}, this shows that $g_i^j$ is concave for $i,j \geq N'$. We fix $j \geq \max\{N',2\}$. By Lemma \ref{lemma_basic_properties_of_P_coconvex_sequences}\ref{lemma_basic_properties_of_P_coconvex_sequences_1}, the sequence $(g_i^j)_{i \in \N}$ is bounded in the space $L^1(G_1)$, and Lemma \ref{lemma_uniform_concave_limit}\ref{lemma_uniform_concave_limit_sequence} implies that there exists a subsequence $(g_{k(i)}^j)_{i \in \N}$ that converges uniformly on $G_j$ to a concave function $g^j_\infty\colon G_j \rightarrow (-\infty,0]$. Using Lemma \ref{lemma_uniform_concave_limit}\ref{lemma_uniform_concave_limit_sequence} iteratively together with the usual diagonal argument, we obtain a strictly increasing function $k \colon \N \rightarrow \N$ such that $(g_{k(i)}^j)_{i \in \N}$ converges uniformly on $G_j$ to $g^j_\infty$ for every $j \geq \max\{N',2\}$.
		With regard to the assumptions of Lemma \ref{lemma_limit_shrinking_sets}, we observe that the sequence $(Q_j)_{j \in \N}$ satisfies $Q_j \subset \interior Q_{j-1}$ for $j \in \N$ and $\bd P \cap \bigcap_{j \in \N}Q_j = G$, and we have
		\begin{equation}
			g_i^j(x)=-\beta_i\cdot \volt{n-m}(Q_j\cap \pi_{\aff F}^{-1}(x)\cap(P \setminus K_i)) \quad \text{for }x \in G_j.
		\end{equation}
		Applying Lemma \ref{lemma_limit_shrinking_sets} to the sequence $(g_i^j \cdot \One_{G_j})_{j \in \N}$, we obtain that $(g^j_\infty)_{j \in \N}$ converges pointwise on $G$ to a function $g^*\colon G \rightarrow (-\infty,0]$ and $\mu|_{\relint G}=[g^*]_{\M(F)}$.
		Because $g^j_\infty$ is concave for every $j \geq N'$, the claim follows.
	\end{proof}

	Our strategy for strengthening Lemma \ref{lemma_m_face_concave_density} is inspired by Fillastre's proof of the reverse Brunn-Minkowski inequality for coconvex bodies. The main idea of this proof is contained in the following lemma, which is taken from \cite[Lem.~4]{fil17}.
	
	\begin{lemma} \label{lemma_fillastre}
		Let $C \subset \R^n$ be a convex cone and $f \colon C \rightarrow \R$ a convex function that is positively homogeneous of degree $m$, i.e., $f(\lambda x)=\lambda^m f(x)$ for all $x \in C \setminus \{0\}$ and $\lambda \geq 0$. Then $f$ is $\frac{1}{m}$-convex.
	\end{lemma}

	To prove the reverse Brunn-Minkowski inequality, Fillastre first shows that the volume functional is convex on a given cone of coconvex bodies, and then uses the homogeneity of the volume to deduce the desired $\tfrac{1}{n}$-convexity. Translated into our setting, this leads us to investigate how limits of $P$-contrast sequences behave under suitable projective transformations (in the sense of projective geometry).
	
	As a result, we obtain the following strengthening of Lemma \ref{lemma_m_face_concave_density}.
	
	\begin{proposition} \label{prop_m_face_alpha_convexity}
		Let $P \subset \R^n$ be a full-dimensional polytope and  $(\mu_i)_{i \in \N}$ an admissible $P$-contrast sequence with limit $\mu$. Let $F\subset P$ be an $m$-face for $0 \leq m < n-1$. Let $f\colon \relint F \rightarrow (-\infty,0]$ be the concave density of $\mu|_{\relint F}$, which exists by Lemma \ref{lemma_m_face_concave_density}. Then $-f$ is $(n-m)^{-1}$-convex
	\end{proposition}

	The outline of the proof of Proposition \ref{prop_m_face_alpha_convexity} is as follows. We replace $\relint F$ by its positive hull $C \coloneqq \pos(\relint F)$, i.e., the convex cone of all non-negative multiples of vectors in $\relint F$, and extend $-f$ on $C$ to a positively homogeneous function $\tilde{f}$ of degree $n-m$. By Lemma \ref{lemma_fillastre}, the claim of Proposition \ref{prop_m_face_alpha_convexity} follows if we show that $\tilde{f}$ is convex. For this, we consider two arbitrary points $\tilde{x}, \tilde{y} \in C \setminus \{0\}$ and show that $\tilde{f}|_{[\tilde{x}, \tilde{y}]}$ is, up to a negative constant, equal to a continuous density of the limit of a contrast sequence, which is concave by Lemma \ref{lemma_m_face_concave_density}. By definition, we have $\tilde{x}=r_1x$ and $\tilde{y}=r_2y$ for some $x,y \in \relint F$ and $r_1,r_2>0$. In order to apply Lemma \ref{lemma_m_face_concave_density}, our intermediate goal is to find a convexity-preserving transformation $\tau$ that satisfies $\tau(x)=\tilde{x}$, $\tau(y)=\tilde{y}$ and that restricts to a dilation on each of the parallel hyperplanes $z+\{y-x\}^\bot$ for $z \in [x,y]$. Such a convexity-preserving transformation $\tau$ can only be affine if $r_1=r_2$; in general, an additional ``degree of freedom'' is needed. This leads us to consider projective transformations, which have the desired properties. For a brief outline of the role of projective transformations in the study of polyhedra, we refer to \cite[Sect.~2.6]{zie07}. To apply Lemma \ref{lemma_m_face_concave_density}, we choose a projective transformation $\tau$ that dilates each hyperplane $z+\{y-x\}^\bot$ for $z \in [x,y]$ with a certain scaling factor and, intuitively speaking, moves the hyperplane $H_{\tilde{z}}= \tilde{z}+\{y-x\}^\bot$, $\tilde{z} \in \R^n$, where the scaling factor is infinite, to the ``hyperplane at infinity''. Truncating $P$ to a polytope $P^+$, we ensure that no point of $P^+$ is mapped ``to infinity'', which entails that $\tau(P^+)$ is a polytope whose boundary contains the segment $[\tilde{x}, \tilde{y}]$. Adapting some parts of the proof of Lemma \ref{lemma_m_face_concave_density}, we finally arrive at the conclusion that $\tilde{f}|_{[\tilde{x}, \tilde{y}]}$ is convex.
	
	\begin{proof}[Proof of Proposition \ref{prop_m_face_alpha_convexity}]
		Without loss of generality, we assume that $0 \in \interior P$. Let $C \coloneqq \pos(\relint F)$ and let $\tilde{f}$ be the unique function $C \rightarrow \R$ that is positively homogeneous of degree $n-m$ and satisfies $\tilde{f}|_{\relint F}=-f$. By Lemma \ref{lemma_fillastre}, it remains to show that $\tilde{f}$ is convex. For this, let $x, y \in \relint F$ with $x\neq y$ and let $r_1,r_2 >0$. We will show that $\tilde{f}|_{[r_1x,r_2y]}$ is convex. Because $\tilde{f}$ is continuous, it suffices to consider the case where $\scpr{x-y,r_1x-r_2y} \neq 0$, which can always be achieved by slightly varying $r_1$ or $r_2$. The general case then follows by approximation. Let $h\colon \R^n \rightarrow \R$ be given by
		\begin{equation}
			h(z) \coloneqq \frac{1}{r_1r_2\norm{x-y}^2} \SCPR{x-y,(r_2-r_1)z+r_1x-r_2y}.
		\end{equation}
		By construction, $h$ is the unique affine functional that satisfies $h(x)=\frac{1}{r_1}$ and $h(z)=\frac{1}{r_2}$ for all $z \in y + \{x-y\}^\bot$. Moreover, due to our assumptions on $r_1$ and $r_2$, we have
		\begin{equation}
			h(0)= \frac{1}{r_1r_2\norm{x-y}^2} \SCPR{x-y,r_1x-r_2y} \neq 0.
		\end{equation} Let $H^+ \coloneqq \{z \in \R^n \mid h(z) \geq \min\{\frac{1}{2r_1},\frac{1}{2r_2}\}\}$ and $P^+ \coloneqq P \cap H^+$. Now we have
		\begin{equation}
			\det \begin{bmatrix}
				I_n&0\\
				\frac{r_2-r_1}{r_1r_2\norm{x-y}^2}(x-y)^\transp&h(0)
			\end{bmatrix} = h(0) \neq 0
		\end{equation}
		and it follows that
		\begin{equation}
			\tau \colon H^+ \rightarrow \R^n, \quad z \mapsto \frac{z}{h(z)}
		\end{equation}
		is a projective transformation that maps $x$ to $r_1x$ and $y$ to $r_2y$ (see \cite[p.~68]{zie07}). We note in passing that if $r_1=r_2$, then $H^+=\R^n$ and $\tau$ is a dilation.
		
		Let $\big[(\beta_i,K_i)\big]_{i \in \N}$ be the component sequence of $(\mu_i)_{i\in\N}$. Without loss of generality, we assume that $K_i^+ \coloneqq K_i \cap H^+ \subset P^+$ holds for all $i\in\N$ and that the restricted sequence $(\mu_i|_{P^+})$ converges weakly. By a change of variables, we have
		\begin{align}
			\int_{\tau(P^+\setminus K_i^+)} f(x) \dint x= \int_{P^+\setminus K_i^+} f(\tau(x))\cdot|\!\det \tau(x)|\dint x,
		\end{align}
		for every $f \in C(\R^n)$. Since the function $(f \circ \tau)\cdot|\!\det \tau|$ is continuous on $H^+$, it follows that $\big[(\beta_i,\tau( K_i^+))\big]_{i \in \N}$ is the component sequence of an admissible $\tau(P^+)$-contrast sequence $(\nu_i)_{i \in \N}$ with limit $\nu$.
		
		Let $F^+$ be the face $F \cap H^+$ of $P^+$ and let $G \coloneqq \frac{1}{2} ([x,y]+F^+)$. For $j \in \N$, we set $G_j \coloneqq \frac{j}{j+1} G + \frac{1}{j+1}F^+$ and define $Q_j$ and $\tilde{Q}_j$ as in \eqref{eq_def_qlj}. %
		By construction, the affine functional $h$ is constant on the affine subspace $\pi^{-1}_{\aff F}(\{x\})$ for all $x \in \aff F$. Therefore, the map $\pi\coloneqq \tau \circ \pi_{\aff F} \circ\tau^{-1}$ is affine, in fact, a projection onto $\tau(\aff F)$. Let $\gamma(\pi)$ be as in Lemma \ref{lemma_fubini_non_orthogonal_projection}. Since $\tau$ is a homeomorphism onto its image, the sequence $(\tau(Q_j))_{j \in \N}$ satisfies the assumptions of Lemma \ref{lemma_limit_shrinking_sets}. For $i,j \in \N$, we define $\tilde{g}_i^j \colon \tau(\relint G_1) \rightarrow (-\infty, 0]$ by
		\begin{equation}
			\tilde{g}_i^j(z)\coloneqq -\beta_i \cdot \gamma(\pi)\cdot\volt{m}(\tau(Q_j)\cap \pi^{-1}(\{z\})\cap\tau(P^+ \setminus K_i^+)),
		\end{equation}
		and $g_i^j \colon \relint G_1 \rightarrow (-\infty, 0]$ by
		\begin{equation}
			g_i^j(z)\coloneqq -\beta_i \cdot\volt{m}(Q_j\cap \pi_{\aff  F}^{-1}(\{z\})\cap (P^+ \setminus K_i^+)).
		\end{equation}
		Using that $\tau$ is a dilation on $\pi^{-1}_{\aff F}(\{z\})$ for all $z \in \aff F$, we have
		\begin{equation} \label{eq_prop_m_face_alpha_convexity_g_star}
			\tilde{g}_i^j(\tau(z))=h(z)^{-n+m}\cdot \gamma(\pi) \cdot g_i^j(z) \quad \text{for all } z \in \relint G_1.
		\end{equation}
		Applying Lemma \ref{lemma_uniform_concave_limit}\ref{lemma_uniform_concave_limit_sequence} iteratively as in the proof of Lemma \ref{lemma_m_face_concave_density}, we obtain that there exists an $N' \in \N$ and a strictly increasing function $k \colon \N \rightarrow \N$ such that $(g_{k(i)}^j)_{i \in \N}$ converges uniformly on $\relint G_j$ for all $j \geq N'$. Since $h^{-n+m}$ is bounded on $\relint F^+$, it follows that $(\tilde{g}_{k(i)}^j)_{i \in \N}$ converges uniformly on $\tau(\relint G_j)$. %
		Applying Lemma \ref{lemma_limit_shrinking_sets} to the admissible $\tau(P^+)$-contrast sequence $(\nu_i)_{i \in \N}$, we obtain that the pointwise limit $\tilde{g}^*(z)\coloneqq \lim_{j \rightarrow \infty}\lim_{i \rightarrow \infty}\tilde{g}_i^j(z)$ exists and satisfies
		\begin{equation}
			\nu|_{\tau(G)} = [\tilde{g}^*]_{M(\tau(F^+))}.
		\end{equation}
		Because $\nu|_{\tau(\relint F^+)}$ has at most one continuous density and $\tilde{g}^*$ is continuous, Lemma \ref{lemma_m_face_concave_density} implies that $\tilde{g}^*|_{[r_1x,r_2y]}$ is concave. Finally, recalling the definition of $\tilde{f}$, we observe that \eqref{eq_prop_m_face_alpha_convexity_g_star} implies $\tilde{g}^*|_{[r_1x,r_2y]}=- \gamma(\pi) \cdot \tilde{f}|_{[r_1x,r_2y]}$, which completes the proof.
	\end{proof}
	
	The proof of our main result is now merely a matter of combining the previous results.
	
	\begin{proof}[Proof of Theorem \ref{thm_main_result}]
		In light of Proposition \ref{prop_main_result_first_inclusion}, it remains to show that every $\mu \in \W(P)$ can be written as a direct sum as in \eqref{eq_thm_main_result}. Let $(K_t)_{t \in [0,1]}$ be a family that realizes $\mu$. Then $[(i,K_{\frac{1}{i}})]_{i \in \N}$ is the component sequence of a $P$-contrast sequence, and the claim follows from Proposition \ref{prop_facet_concave_density} and Proposition \ref{prop_m_face_alpha_convexity}. The claim that $\W(P)$ is a convex cone follows from the fact that $\W(P)$ is the sequential closure of $\W^{\triangleleft}(P)$, which is clearly a convex cone.
	\end{proof}

	\begin{remark}
		Theorem \ref{thm_main_result} implies that, for $n \in \N$, the sum of two $\frac{1}{n}$-convex functions is again $\frac{1}{n}$-convex. This fact can also be deduced from the Hölder inequality.
	\end{remark}
	
	\section{A condition for polytopal maximizers of the isotropic constant} \label{sect_application}
	
	In this section, we use Theorem \ref{thm_main_result} to derive a necessary condition for polytopal maximizers of the isotropic constant and other geometric functionals.%
	
	Let $\phi \colon \Knn \rightarrow \R$ be a functional of the form \eqref{eq_functional_phi} and $P \in \Pnn$. We define the function $h_{\phi,P} \colon \R^n \rightarrow \R$ by
	\begin{equation} \label{eq_def_h_phi}
		h_{\phi,P} \coloneqq \sum_{i=1}^m \frac{\partial g}{\partial x_i} \left(\int_{P} f_1(x) \dint x, \dots, \int_{P} f_m(x) \dint x\right) \cdot f_i.
	\end{equation}
	If $P$ is a local maximizer of $\phi$, then it follows from Theorem \ref{thm_main_result} that
	\begin{equation} \label{eq_first_order_condition}
		 \int h_{\phi,P} \dint \mu \leq 0
	\end{equation}
	holds for all $\mu \in \W(P)$. In particular, \eqref{eq_first_order_condition} holds with equality for all $\mu \in \W^\pm(P)$. This leads us to the following definition.
	\begin{definition}
		Let $\phi$ be as above. A polytope $P \in \Pnn$ is called \textit{perturbation stable} (with respect to $\phi$) if $\W(P)$ is contained in the polar cone of $\{h_{\phi,P}\}$, i.e., if \eqref{eq_first_order_condition} holds for all $\mu \in \W(P)$. The polytope $P$ is called \textit{weakly perturbation stable} if $\W^\pm(P)$ is contained in the annihilator of $\{h_{\phi,P}\}$, i.e., if \eqref{eq_first_order_condition} holds for all $\mu \in \W^\pm(P)$.
	\end{definition}
	For $F \in \Phi(P)$, we denote the restriction of the Lebesgue measure on $\aff F$ to $F$ by $\vol_F$ and set
	\begin{equation}
		\volt{\Phi} \coloneqq \sum_{F \in \Phi(P)} \volt{F} \in \M^+(\bd P).
	\end{equation}
	We now embed the set of all elements of $\W(P)$ with square-integrable densities into the Hilbert space $L^2(\bd P, \vol_\Phi)$ by setting
	\begin{equation}
		\W^2(P) \coloneqq \left\{f \in L^2(\bd P, \volt{\Phi}) \MID \exists \mu \in \W(P)\colon \dint\mu=f \dint \volt{\Phi} \right\}.
	\end{equation}
	Since $\W^\triangleleft(P)$ is dense in $\W(P)$ and every element of $\W^\triangleleft(P)$ has a square-integrable density, a polytope $P$ is perturbation stable if and only if $\W^2(P)$ is contained in the polar cone of $\{h_{\phi,P}\}$. This fact entails the following result.
	\begin{proposition} \label{prop_necessary_condition}
		Let $P \in \Pnn$ and let $\phi$ and $h_{\phi,P}$ be as above. We denote the metric projection onto the closed convex cone $\W^2(P)$ in the Hilbert space $L^2(\bd P, \vol_\Phi)$ by $\pi_{\W}$. If $P$ is a local maximizer of $\phi$, then $\pi_{\W}(h_{\phi,P})=0 \in L^2(\bd P, \vol_\Phi)$.
	\end{proposition}
	\begin{proof}
		The claim that $\W^2(P)$ is closed follows from the fact that $L^2$-convergence implies $L^1$-convergence (because $\vol_\Phi$ is finite), and the $L^1$-norm coincides with the total variation norm. If there exists an $f \in \W^2(P)$ with $\norm{h_{\phi,P}-f}_{L^2} < \norm{h_{\phi,P}}_{L^2}$, then
		\begin{equation}
			-2\scpr{h_{\phi,P},f}<-2\scpr{h_{\phi,P},f}+\norm{f}_{L^2}^2<0,
		\end{equation}
		and $f$ is a witness for the fact that $P$ is not perturbation stable.
	\end{proof}
	The necessary condition from Proposition \ref{prop_necessary_condition} resembles the problem of convex regression in statistics. Since $\W^2(P)$ is a direct sum, the condition can be restated as follows: For each face $F\subset P$, the zero function is the unique minimizer of the risk minimization problem $\min_{f \in \mathcal{F}}\E[(-f(X)-h_{\phi,P}(X))^2]$, where $X$ is a random vector that is uniformly distributed on $\relint F$ and $\mathcal{F}$ denotes the class of square-integrable $(n-\dim F)^{-1}$-convex functions on $\relint F$.
	
	Rephrased for the isotropic constant, Proposition \ref{prop_necessary_condition} reads as follows.
	\begin{corollary}
		Let $P \in \Pnn$ be an isotropic local maximizer of the isotropic constant and let $\pi_{\W}$ be as in Proposition \ref{prop_necessary_condition}. Then $\pi_{\W}(x \mapsto \norm{x}^2_2-n-2)=0 \in L^2(\bd P, \vol_\Phi)$. Let $F \subset P$ be a facet and $X$ be uniformly distributed on $F$. Using that \eqref{eq_first_order_condition} holds with equality for the signed measures $\mu_i \in \M(F)$ with densities $x \mapsto x_i$, $i \in [n]$, we have
		\begin{equation}
			\E[\norm{X}^2_2X]=(n+2)\E[X].
		\end{equation}
	\end{corollary}

	\subsection*{Acknowledgments}
	
	This work was supported by the Deutsche Forschungsgemeinschaft (DFG), Graduiertenkolleg ``Facets of Complexity'' (GRK 2434). Moreover, I thank Ansgar Freyer and the referee for their very helpful comments and suggestions.
	
	\printbibliography
\end{document}